\documentclass[11pt,leqno,british]{amsart}
\usepackage{ifpdf}

\usepackage[T1]{fontenc}

\usepackage{eucal}
\usepackage{url}
\usepackage[pagebackref]{hyperref}
\usepackage{amsfonts}
\usepackage{amsthm}
\usepackage{amsmath}
\usepackage{amscd}
\usepackage{amssymb}
\usepackage{graphicx}
\usepackage{rotating}

\usepackage[small,nohug,heads=LaTeX]{diagrams}
\diagramstyle[labelstyle=\scriptstyle]

\def\CC {{\mathbb C}}     
\def\HH {{\mathbb H}}     
\def\NN {{\mathbb N}}     
\def\RR {{\mathbb R}}     
\def\SS {{\mathbb S}}     
\def\ZZ {{\mathbb Z}}     

\def\ring#1{\ifmmode \mathaccent'027 #1\else \rm\accent'027 #1\fi}

\newcommand{\ri}{{\mathrm i}}

\def\ul  {\underline}

\def\im {\mathfrak{im}}

\def\mc {\mathcal}
\def\mk {\mathfrak}

\def\st {\mathrm{Stab}}

\def \bd {\begin{diagram}}
\def \ed {\end{diagram}} 
\def\be  {\begin{eqnarray}}
\def\ee  {\end{eqnarray}}
\def\ben {\begin{eqnarray*}}
\def\een {\end{eqnarray*}}

\def\bpr {\begin{proof}[Proof]}
\def\epr {\end{proof}}
\def\bsp {\begin{split}}
\def\esp {\end{split}}
\def\bprr {\begin{proof}[solution]}
\def\bpru {\begin{proof}[hint]}
\def\bpro {\begin{proof}[answer]}
\def\bcd {\begin{CD}}
\def\ecd {\end{CD}}

\newcommand{\abs}[1]{\left\vert#1\right\vert}
\newcommand{\scal}[1]{\left\langle#1\right\rangle}

\def\Hom {\mathrm{Hom}}

\newtheorem{theorem}{Theorem}[section]
\newtheorem{lemma}[theorem]{Lemma}
\newtheorem{prop}[theorem]{Proposition}
\newtheorem{coro}[theorem]{Corollary}
\newtheorem{remark}[theorem]{Remark}
\newtheorem{df}[theorem]{Definition}

\begin{document}

\title[Bridgeland stability conditions on the  acyclic triangular  quiver]%
{Bridgeland stability conditions on the  acyclic  triangular  quiver}
 
\author{George Dimitrov}
\address[Dimitrov]{Universit\"at Wien\\
Oskar-Morgenstern-Platz 1, 1090 Wien\\
\"Osterreich}
\email{gkid@abv.bg}

\author{Ludmil Katzarkov}
\address[Katzarkov]{Universit\"at Wien\\
Oskar-Morgenstern-Platz 1, 1090 Wien\\
\"Osterreich }
\email{lkatzark@math.uci.edu}
\begin{abstract} Using  results  in a previous paper  ``Non-semistable exceptional objects in hereditary categories'', we focus here  on studying  the topology of the space of Bridgeland stability conditions on $D^b(Rep_k(Q ))$, where   $Q =\begin{diagram}[size=0.8em]
        &       &  \circ  &       &    \\
        & \ruTo &         & \luTo &       \\
 \circ  & \rTo  &         &       &  \circ
\end{diagram}$.  In particular,  we prove that this space is contractible (in the previous paper it was shown that it is connected).  \end{abstract}

\maketitle
\setcounter{tocdepth}{2}
\tableofcontents

\section{Introduction}
  
In 1994 Maxim Kontsevich interpreted a duality coming from physics in a 
consistent, powerful mathematical framework called Homological Mirror 
Symmetry (HMS).
HMS is now the foundation of a wide range of contemporary mathematical 
research.
Numerous works by many authors have demonstrated the interaction of 
mirror symmetry and HMS with a wide range of new and subtle mathematical 
structures.
One of this structures is the moduli space of stability conditions.

The study of stability in triangulated categories was initiated by M. 
Douglas
and  mathematically by T. Bridgeland.
The majority of the activity since then
has focused on categories of algebro-geometric origin.
Signicant work in this direction is due to T. Bridgeland, A. King, E. Macr\'i, S. Okada,  Y. Toda, A. 
Bayer, J. Woolf,  J. Collins, A. Polishchuck et. al.

In previous works \cite{DK1}, \cite[joint with F. Haiden and M. Kontsevich]{DHKK}    we developed   results and ideas by T. Bridgeland \cite{Bridg1}, A. King \cite{King}, E. Macr\'i \cite{Macri},  J. Collins and A. Polishchuck \cite{CP}.

Recently in \cite{Woolf} J.  Woolf showed  classes   of categories with  contractible component in the space of stability conditions. His  paper generalizes  and unifies various known results for stability spaces of specific categories, and settles some conjectures about the stability spaces associated to Dynkin quivers, and to their Calabi-Yau-N Ginzburg algebras.  However the  results in \cite{Woolf} do not cover  tame representation type quivers, these quivers  are beyond the scope of \cite{Woolf}.

In the  present paper we give a new example of a tame representation type quiver with contractible space of stability conditions. This paper  is  a natural consequence of our previous  paper   \cite{DK1}. Both are based  on   ideas of  E. Macr\'i, which he gave  in \cite{Macri}  studying $ \st(D^b(K(l))$, where $K(l)$ is the $l$-Kronecker quiver.
\vspace{3mm}

1.1. T. Bridgeland   defined  in  \cite{Bridg1} the  space  of stability conditions on  a triangulated category  $\mc T$, denoted by $\st(\mc T)$,  and proved that it is a  a complex manifold on which act the groups $\widetilde{GL}^+(2,\RR)$ and ${\rm Aut}(\mc T)$. To any bounded t-structure of $\mc T$ he assigned a family of stability conditions.

E. Macri constructed in  \cite{Macri}   stability conditions using exceptional collections and the action of $\widetilde{GL}^+(2,\RR)$ on $\st(\mc T)$. Applying results in \cite{BBD}, he showed  that the extension closure of   a full Ext-exceptional collection\footnote{An exceptional collection $\mc E$ is said to be \textit{Ext-exceptional} if $\Hom^{\leq 0}(E_i,E_j)=0$ for $0\leq j<j\leq n$.} $\mc E=(E_0,E_1,\dots, E_n)$ in  $\mc T$ is a bounded t-structure.
The stability conditions obtained from this t-structure together with their  translations by the right action of $\widetilde{GL}^+(2,\RR)$ will be  referred to as \textit{generated by $\mc E$}.  

E. Macr\`i, studying   $ \st(D^b(K(l))$ in  \cite{Macri},     gave  an idea for producing an  exceptional pair generating a given stability condition  $\sigma$ on  $D^b(K(l))$, where $K(l)$ is the $l$-Kronecker quiver. 

 We defined in \cite{DK1}  the notion of a \textit{$\sigma$-exceptional collection} (\cite[Definition 3.19]{DK1}), so that the full $\sigma$-exceptional collections are exactly the exceptional collections which generate $\sigma$, and  we   focused  on  constructing $\sigma$-exceptional collections from a given $\sigma \in \st(D^b(\mc A))$, where $\mc A$ is a hereditary, $\hom$-finite, abelian category.  We 
 developed tools   for constructing $\sigma$-exceptional collections of length at least three in  $D^b(\mc A)$.  These tools are based on  
   the notion of       \textit{regularity-preserving hereditary category},  introduced in \cite{DK1} to avoid difficulties related to the Ext-nontrivial couples (couples of exceptional objects in $\mc A$ with ${\rm Ext}^1(X,Y)\neq 0$ and ${\rm Ext}^1(Y,X)\neq 0$).
	
 After a detailed study of  the exceptional objects of the affine quiver $Q$ (see figure \eqref{Q1} below)  it was shown in \cite{DK1} that $Rep_k(Q)$ is regularity preserving and  the newly obtained methods for constructing $\sigma$-triples were applied to the case   $\mc A = Rep_k(Q)$.   As a result we  obtained        the following theorem:

 \begin{theorem}[\cite{DK1}] \label{main theorem for Q in intro} Let $k$ be an algebraically closed field.
 For each $\sigma \in \st(D^b(Rep_k(Q)))$ there exists a full
$\sigma$-exceptional collection.
\end{theorem}

 In other words,  all stability conditions on $D^b(Q)$ are generated by exceptional collections (in this case exceptional triples).  This theorem implies
 that $\st(D^b(Q))$ is connected \cite[Corollary 10.2]{DK1}.

 Using  Theorem \ref{main theorem for Q in intro} and the data about the exceptional collections  
 given in \cite[Section 2]{DK1},   we prove here the following:
\begin{theorem} \label{main theo} Let $k$ be an algebraically closed field.  Let  $Q $ be the following quiver:
\be \label{Q1} Q = \begin{diagram}[1em]
   &       &  \circ  &       &    \\
   & \ruTo &    & \luTo &       \\
\circ  & \rTo  &    &       &  \circ
\end{diagram}. \ee The space of Bridgeland stability conditions  $\st(D^b(Rep_k(Q))$ is a contractible (and connected) manifold, where $D^b(Rep_k(Q))$ is the derived category of representations of $Q$.  
\end{theorem}

 1.2. We give now more details  about  the structure of $\st(D^b(Rep_k(Q)))$ and about  the proof of Theorem \ref{main theo}.

	We call an exceptional pair $(E,F)$ in $D^b(Rep_k(Q))$ a 2-Kronecker pair if $\hom^{\leq 0}(E,F)=0$, and $\hom^{1}(E,F)=2$. 
	Recall that the Braid group on two strings $B_2\cong \ZZ$ acts on the set of equivalence classes of exceptional pairs in $\mc T$. \footnote{Here we take the equivalence $\sim$ explained in \textbf{Some notations} and it is clear when a given equivalence class w.r. $\sim$  will be called a 2-Kronecker pair} The set of  equivalence classes  of 2-Kronecker pairs is invariant  under  the action of $B_2$.   In Subsection \ref{the two orbits}   are described    the orbits of this action on the  2-Kronecker pairs (using \cite[Corollary 2.9]{DK1}). There are two such orbits and 
	in terms of our notations   they are  $\{(a^m,a^{m+1}[-1])\}_{m\in \ZZ}$ and $\{(b^m,b^{m+1}[-1])\}_{m\in \ZZ}$  (see Remark \ref{T_1234}). 
	
	It  turns out that the exceptional objects of
$D^b(Rep_k(Q))$ can be grouped as follows $\{a^m\}_{m\in \ZZ}  \cup \{M,M'\}\cup \{b^m\}_{m\in \ZZ}$, where $\{M,M'\}\subset \ Rep_k(Q)$ is the unique Ext-nontrivial couple of $D^b(Rep_k(Q))$.
	
	Let $\mk{T}_a^{st}$ and  $\mk{T}_b^{st}$ be the stability conditions generated by the exceptional triples containing a subsequence  of the  from   $(a^m[p],a^{m+1}[q])$ and $(b^m[p],b^{m+1}[q])$ for some $m,p,q \in \ZZ$, respectively. Using Theorem \ref{main theorem for Q in intro} we show in Section \ref{the union} that $\st(D^b(Rep_k(Q)))=\mk{T}_a^{st}\cup (\_,M,\_) \cup (\_,M',\_)\cup \mk{T}_b^{st}$, where $(\_,M,\_) \cup (\_,M',\_)$ denotes the set of stability conditions generated by triples of the form $(A,M[p],C)$ or $(A,M'[p],C)$ with $p\in \ZZ$ (these turn out to be the  triples $(A,B,C)$ for which $\dim(\Hom^i(A,B))\leq 1 $, $\dim(\Hom^i(A,C))\leq 1 $, $\dim(\Hom^i(B,C))\leq 1 $ for all $i\in \ZZ$).  
	
	The main steps are as follows. 
	In Section \ref{mk T_12 cap mk T_43}  we how that $\mk{T}_a^{st} \cap \mk{T}_b^{st}=\emptyset$. In Section \ref{T_a T_b cont} we show that  $\mk{T}_a^{st}$ and $\mk{T}_b^{st}$ are contractible. In Section \ref{connecting} we connect  $\mk{T}_a^{st}$ and $\mk{T}_b^{st}$  by $(\_,M,\_)\cup (\_,M',\_)$ and  show that in this procedure the contractibility is preserved. 
		
The  theorem from topology  which we use to glue stability conditions generated by different exceptional triples is the  Seifert-van Kampen    theorem, modified about contractile subsets in manifolds  (see Remark \ref{VK}).    In Section \ref{general remarks}  are  given several important  tools, which we use throughout to analyze the intersection of the sets of  stability conditions generated by different exceptional collections. These tools  are extensions of results and ideas  in \cite{King}, \cite{Macri}, \cite{DHKK}, \cite{DK1}.    In the final step (Section \ref{connecting}) we utilize as such a tool also the relation   $\bd  R & \rDotsto & (S,E)  \ed $  between a $\sigma$-regular object  $R$ and an exceptional pair generated by it (introduced in \cite{DK1}).  

In Section \ref{the exceptional objects} we organize in a better way the obtained  in \cite[Section 2]{DK1} data about $\Hom(X,Y)$ and ${\rm Ext}^1(X,Y)$, where $X,Y$ vary throughout the exceptional objects of $Rep_k(Q)$, and we add some observations about the behavior of the central charges of the exceptional objects, which  are   very essential for the proof of Theorem \ref{main theo} as well.  

Today, in view of the parallel between dynamical systems and categories \cite{DHKK}, \cite{BS} and  in view of the Motivic Donaldson Thomas invariants \cite{KS} the importance of studying the topology of the  space  of Bridgeland  stability conditions is even bigger.

We still do not understand  the meaning of the obtained picture about  $\st(D^b(Rep_k(Q)))$. We  hope that an understanding  of this  meaning  will open a way to analyzing more cases.

\textit{{\bf Acknowledgements:}}
The authors wish to express their gratitude to   Tom Bridgeland,  Dragos Deliu, Fabian Haiden,  Umut Isik,  Maxim Kontsevich,  Alexander Noll,  Tony Pantev, Pranav Pandit for their   interest in this paper. 

Both the authors were funded by NSF DMS 0854977 FRG, NSF DMS 0600800, NSF DMS 0652633
FRG, NSF DMS 0854977, NSF DMS 0901330, FWF P 24572 N25, by FWF P20778 and by an
ERC Grant. 

The second author was funded   by DMS-1265230 Wall Crossings in Geometry and Physics,
DMS-1201475 Spectra Gaps Degenerations and Cycles, and by  
OISE-1242272 PASI On Wall Crossings Stability Hodge Structures \& TQFT Grant.

 \textit{\textbf{Some notations.}} In these notes the letters ${\mathcal T}$ and $\mc A$ denote always  a triangulated category and an abelian category, respectively, linear over an  algebraically closed  field  $k$. T

he shift functor  in ${\mathcal T}$ is designated by $[1]$.   We write $\Hom^i(X,Y)$ for  $\Hom(X,Y[i])$ and  $\hom^i(X,Y)$ for  $\dim_k(\Hom(X,Y[i]))$, where $X,Y\in \mc T$.  For $X,Y\in\mc  A$,  writing $\Hom^i(X,Y)$, we consider $X,Y$ as elements in   $\mc T=D^b(\mc A)$, i.e.  $\Hom^i(X,Y)={\rm Ext}^i(X,Y)$.

  We write $\langle  S \rangle  \subset \mc T$ for  the triangulated subcategory of $\mc T$ 
   generated by $S$, when $S \subset Ob(\mc T)$.

  An \textit{exceptional object}  is an object $E\in \mc T$ satisfying $\Hom^i(E,E)=0$ for $i\neq 0$ and  $\Hom(E,E)=k $. We denote by ${\mc A}_{exc}$, resp. $D^b(\mc A)_{exc}$,  the set of all
    exceptional objects of  $\mc A$, resp. of  $D^b(\mc A)$.

An \textit{exceptional collection} is a sequence $\mc E = (E_0,E_1,\dots,E_n)\subset \mc T_{exc}$ satisfying $\hom^*(E_i,E_j)=0$ for $i>j$.    If  in addition we have $\langle \mc E \rangle = \mc T$, then $\mc E$ will be called a full exceptional collection.  For a vector $\textbf{p}=(p_0,p_1,\dots,p_n)\in \ZZ^{n+1}$ we denote $\mc E[\textbf{p}]=(E_0[p_0], E_1[p_1],\dots, E_n[p_n])$. Obviously  $\mc E[\textbf{p}]$ is also an exceptional collection. The exceptional collections of the form   $\{\mc E[\textbf{p}]: \textbf{p} \in \ZZ^{n+1} \}$ will be said to be shifts of $\mc E$. 

For two exceptional collections $\mc E_1$, $\mc E_2$ of equal length we  write $\mc E_1 \sim \mc E_2$ if $\mc E_2 \cong \mc E_1[\textbf{p}]$ for some $\textbf{p} \in \ZZ^{n+1}$.

 An abelian category $\mc A$ is said to be hereditary, if ${\rm Ext}^i(X,Y)=0$ for any  $X,Y \in \mc A$ and $i\geq 2$,  it is said to be of finite length, if it is Artinian and Noterian.

  For any quiver $Q$ we denote by  $D^b(Rep_k(Q))$ or just by $D^b(Q)$  the derived category of the category of representations of $Q$.
	
	For any $a\in \RR$ and any  complex number $z \in {\rm e}^{\ri \pi a} \cdot (\RR + \ri \RR_{>0})$, respectively  $z \in {\rm e}^{\ri \pi a} \cdot \left (\RR_{<0} \cup (\RR + \ri \RR_{>0}) \right )$,  we denote by $\arg_{(a,a+1)}(z)$, resp. $\arg_{(a,a+1]}(z)$,  the unique $\phi \in (a,a+1)$, resp. $\phi \in (a,a+1]$, satisfying $z=\abs{z} \exp(\ri \pi \phi)$.   
	
For a non-zero complex number $v \in \CC$ we denote   the two connected components of  $\CC \setminus \RR v $ by:
\begin{gather} \label{complement of a line} v^c_+ = v \cdot  (\RR + \ri \RR_{>0}) \qquad  v^c_- = v \cdot (\RR - \ri \RR_{>0}) \qquad \qquad v \in \CC \setminus \{0\}.\end{gather}

For $b\in (a,a+1)$, $c\in(a-1,a)$ $r_1>0$, $r_2>0$ we have
\begin{gather} 
\arg_{(a,a+1)}(r_1 \exp(\ri \pi a)+r_2 \exp(\ri \pi b))=a+\arg_{(0,1)}(r_1+r_2 \exp(\ri \pi (b-a))) \nonumber\\[-3mm] 
 \label{arg1}   \\[-3mm]  \arg_{(a-1,a)}(r_1 \exp(\ri \pi a)+r_2 \exp(\ri \pi c))=a+\arg_{(-1,0)}(r_1+r_2 \exp(\ri \pi (c-a))). \nonumber \end{gather}
These formulas imply that for  $c\in(a-1,a)$, $r_1>0$, $r_2>0$ we have
\begin{gather} \label{arg2}
\arg_{(a-1,a)} \left ( r_1 \exp(\ri \pi a) + r_2 \exp(\ri \pi c)\right )=-\arg_{(-a,-a+1)} \left ( r_1 \exp(-\ri \pi a) + r_2 \exp(-\ri \pi c)\right ).
\end{gather}

\section{Some general remarks} \label{general remarks} 
Here we  give tools which will be used throughout  to analyze the intersection of the  sets of stability conditions generated by different exceptional collections (Propositions \ref{all exc semistable},  \ref{phi_1>phi_2}, \ref{mutations} and Lemmas \ref{three comp factors}, \ref{two comp factors}). 
 A description of the set of stability conditions   generated by all shifts of a fixed exceptional triple is given in Proposition \ref{lemma for f_E(Theta_E)}, which is also important for the rest of the paper. Due to Remark \ref{for n geq 3} it seems that Proposition \ref{lemma for f_E(Theta_E)} can not be  generalized straightforwardly to the case of exceptional collections of length bigger than $3$.   
\subsection{Basic facts and notations related to Bridgeland stability conditions}  We use freely the axioms and notations  on stability conditions introduced by Bridgeland in \cite{Bridg1} and some additional notations  used in \cite[Subsection 3.2]{DK1}.  In particular, for $\sigma = ({\mathcal P}, Z) \in \st(\mathcal T)$ we  denote by
$\sigma^{ss}$ the  set of  $\sigma$-semistable
objects, i. e. \be
\label{sigma^{ss}} \sigma^{ss}=\cup_{t \in \RR} {\mathcal P}(t)\setminus \{0\}.
\ee 
 For any interval $I\subset \RR$    the extension closure of the slices $\{\mc P(x)\}_{x\in I}$ is denoted by $\mc P(I)$ in \cite{Bridg1}.    The nonzero objects in the subcategory $\mc P(I)$  are exactly those $X\in \mc T\setminus \{0\}$, which  satisfy $\phi_\pm(X)\in I$, i. e. whose HN factors have phases in $I$. In particular,  if  $X \in \mc P(a-1,a]\setminus \{0\}$ then  $Z(X) \in \exp(\ri \pi a)_-^c\cup \RR_{>0} \exp(\ri \pi a) $.

 From \cite{Bridg1} we know that for any  $\sigma = (\mc P, Z) \in \st(\mc T)$ and any $t \in \RR$ the subcategory $\mc P(t,t+1]$ is a heart of a bounded t-structure. In particular  $\mc P(t,t+1]$  is an abelian category, whose short exact sequences  are exactly these sequences
  $\bd A & \rTo ^{\alpha}& B & \rTo^{\beta} & C \ed $ with $A,B,C \in \mc P(t,t+1]$, s. t.  for some $\gamma : C \rightarrow A[1]$ the sequence
$\bd A & \rTo ^{\alpha}& B & \rTo^{\beta} & C & \rTo^{\gamma} A[1]\ed $ is a  triangle in $\mc T$. Using these remarks,  the HH filtration and by drawing pictures one easily shows the following properties:

\begin{remark} \label{arg remark} 
Let $t\in \RR$ and  $X \in  \mc P(a-1,a]$. Then: 
\begin{itemize}
	\item[(a)] If  $ X \not \in \sigma^{ss}$ then $\phi_-(X) < \arg_{(a-1,a]}(Z(X)) < \phi_+(X) $.
	\item[(b)] $X \not \in \sigma^{ss} $ iff there exists a monic arrow  $X'\rightarrow X$  in the abelian category  $ P(a-1,a]$ satisfying  $\arg_{(a-1,a]}(Z(X')) >\arg_{(a-1,a]}(Z(X))$. 
	\item[(c)] If $Z(X)\in v_+^c$ for some $v\in \CC^*$ with    $v=\abs{v}\exp(\ri \pi t)$ and $a-1 \in(t,t+1)$ or  $a \in(t,t+1)$, then $\arg_{(a-1,a]}(Z(X))=\arg_{(t,t+1)}(Z(X))$. In particular,   when $X \in \sigma^{ss}$, we have: \\  $\phi(X)=\arg_{(t,t+1)}(Z(X))$.
\end{itemize}
\end{remark}

\subsection{ Some remarks on \texorpdfstring{$\sigma$}{\space}-exceptional collections}

E. Macr\`i proved in \cite[Lemma 3.14]{Macri} that the extension closure  ${\mc A}_{\mc
E}$   of a full Ext-exceptional collection ${\mc
E}=(E_0,E_1,\dots,E_n)$ in $\mc T$ is
  a heart of a bounded t-structure. Furthermore,
 ${\mc A}_{\mc E}$ is of finite length and $E_0, E_1,\dots, E_n$ are the simple objects in it.  
  By Bridgeland's \cite[Proposition 5.3]{Bridg1} from the bounded t-structure ${\mc A}_{\mc E}$  is produced a family of stability conditions, which we denote by  $\HH^{\mc A_{\mc E}}\subset \st(\mc T)$ or sometimes just $\HH^{\mc E}\subset \st(\mc T)$.  

 For a given $\sigma \in \st(\mc T)$ we define a
    \textit{$\sigma$-exceptional collection} (\cite[Definition 3.19]{DK1})  as an  Ext-exceptional collection ${\mc E } = (E_0,E_1,\dots,E_n)$,     s. t. the
     objects $\{E_i\}_{i=0}^n$ are $\sigma$-semistable, and  $\{\phi(E_i)\}_{i=0}^n \subset (t,t+1)$
 for some $t\in \RR$. The following Proposition is basic for this paper:   

\begin{prop} \label{all exc semistable} Let  $\mc T$ be a $k$-linear triangulated category and $\sigma=(\mc P, Z) \in \st(\mc T)$. Let $\mc E=(E_0,E_1,\dots, E_n)$ be a full $\sigma$-exceptional collection such that  $\phi(E_i)\geq \phi(E_{i+1})$ and  $\hom^1(E_i,E_{i+1})$ $\neq $ $ 0$  for some $i\in \{0,1,\dots,n-1\}$. Let $\mc A_{i,i+1}$  be the extension closure of  $E_i, E_{i+1}$  in $\mc T$.  Then each element in ${\mc T}_{exc}\cap \mc A_{i,i+1}$ is semistable. 
\end{prop}
\bpr  If $\phi(E_i)=\phi(E_{i+1})=t$, then $\mc A_{i,i+1}\subset \mc P(t)$ and hence  all non-zero objects in $\mc A_{i,i+1}$  are semistable, therefore we can assume that $\phi(E_i)>\phi(E_{i+1})$. 

 By \cite[Corollary 3.20]{DK1} we have $\sigma \in \Theta_{\mc E}'=\HH^{\mc E}\cdot\widetilde{GL}^+(2,\RR)$.  Since the action of $\widetilde{GL}^+(2,\RR)$ does not change the order of the phases, we can assume that $\sigma=(\mc P, Z) \in \HH^{\mc E}$, which means that the extension closure of $\mc E$ is the t-structure $\mc P(0,1]$ and \be \label{all exc semistable1}  \phi(E_j)=\arg_{(0,1]}(Z(E_j)) \qquad j=1,\dots,n. \ee 

Let us denote  $ \mc       T_{i,i+1} =\left \langle E_i,E_{i+1}\right \rangle $.  From \cite[Proposition 3.17]{DK1} we have a projection map $\HH^{\mc E} \rightarrow \HH^{{\mc A}_{i,i+1}}\subset \st(\mc T_{i,i+1})$ and it maps $\sigma =(\mc P,Z)$ to a stability condition $\sigma' =(\mc P',Z')\in  \HH^{{\mc A}_{i,i+1}}$ with $Z'(E_i)=Z(E_i)$, $Z'(E_{i+1})=Z(E_{i+1})$ and $\{{\mc P}'(t)=\mc P(t)\cap \mc       T_{i,i+1}\}_{t\in\RR}$. Therefore it remains to show that the objects in ${\mc T}_{exc}\cap \mc A_{i,i+1}$ are $\sigma'$-semistable.

From \cite[Lemma 3.22]{DHKK} we have that $\mc A_{i,i+1}$ is a bounded t-structure in $\mc T_{i,i+1}$ and an equivalence of abelian categories $F:\mc A_{i,i+1} \rightarrow Rep_k(K(l))$ with $F(E_i)=s_1$, $F(E_{i+1})=s_2$, where  $l=\hom^1(E_i,E_{i+1})$ and $s_1$, $s_2$ are the simple representations of $K(l)$ with $k$ at the source, sink, respectively. This equivalence maps   $\sigma'\in \HH^{\mc A_{i,i+1}}$  to a stability condition   
 \begin{gather} \sigma'' = (\mc P'', Z'')\in \HH^{Rep_k(K(l))} \subset \st(D^b(K(l))) 
 \quad Z(E_i)=Z''(s_1), Z(E_{i+1})=Z''(s_2). \nonumber \end{gather}  If $E\in \mc T_{exc}\cap \mc A_{i,i+1}$, then by the fact that $F$ is an equivalence of abelian categories it follows that $F(E)\in Rep_k(K(l))$ is an exceptional representation. Since $\{F(\mc P'(t))=\mc P''(t) \}_{t\in (0,1]}$, it remains to show that each exceptional representation of   $Rep_k(K(l))$ is $\sigma''$-semistable. 

Let $\rho\in Rep_k(K(l))_{exc}$. Then the dimension vector $\ul{\dim}(\rho)=(n,m)\in (n,m)$ is a real  root of $K(l)$, furthermore it is a Schur root. From \eqref{all exc semistable1} we have $\arg(Z''(s_1))>\arg(Z''(s_2))$.  By the arguments in the proof of \cite[Lemma 3.19]{DHKK} using  a theorem by King (\cite[Proposition 4.4]{King} ) and  $\arg(Z''(s_1))>\arg(Z''(s_2))$  we obtain a $\sigma''$-stable representation $X\in Rep_k(K(l))$ with  $\ul{\dim}(X)=(n,m)$. Since $X$ is stable, it is simple in $\mc P''(t)$, where $t=\phi''(X)$, in particular it is indecomposable in $\mc P''(t)$. Since  $\mc P''(t)$  is a thick subcategory (see \cite[Lemma 3.7]{DK1}), it follows that   $X$ is indecomposable in  $Rep_k(K(l))$. Since $\ul{\dim}(\rho)$ is a real root and both $X$, $\rho$ are indecompsable representations, the equality  $\ul{\dim}(\rho)=\ul{\dim}(X)$ implies $\rho \cong  X$(see \cite[Theorem 2, c)]{Kac}). The proposition  follows. 
\epr
Other statements, which will be widely  used in the next sections   are Propositions  \ref{lemma for f_E(Theta_E)}, \ref{phi_1>phi_2} and \ref{mutations}.  For the proof of Proposition  \ref{lemma for f_E(Theta_E)}   it is useful to define:  
\begin{df} \label{S(I)} Let $n\geq 1$ be an integer.  Let  $\mc I = \{I_{ij}=(l_{ij},r_{ij})\subset \RR\}_{ 0 \leq i <j \leq n}$  be a  family of non-empty open  intervals, and let  $ {\mk l}=\{ l_{ij} \in \{-\infty\}\cup \RR \}_{0 \leq i <j \leq n}$,  ${\mk r}= \{ r_{ij} \in  \RR \cup \{+\infty\} \}_{0 \leq i <j \leq n}$  be the corresponding families of left and right endpoints. 

 We will denote the following open convex set $ \{(y_0,y_1,\dots,y_n)\in \RR^{n+1} : y_i - y_j \in I_{ij} \ i<j\}\subset \RR^{n+1}$ 
  by $S^n(\mc I)$ or $S^n({\mk l}, {\mk r})$.
\end{df}
For a full Ext-exceptional collection $\mc E=(E_0, E_1, \dots, E_n)$ in $\mc T$  we denote  $\Theta_{\mc E}'=\HH^{\mc E}\cdot\widetilde{GL}^+(2,\RR)$. If  $\mc E$ is a full Ext-exceptional collection, then  we have (see \cite[Remark 3.21]{DK1}): 
\begin{gather}  \label{Theta_mc E'} \Theta_{\mc E}'=\HH^{\mc E}\cdot\widetilde{GL}^+(2,\RR)=\left \{\sigma  : {\mc E}\subset \sigma^{ss} \ \mbox{and} \  \abs{\phi^\sigma(E_i)-\phi^\sigma(E_j)}<1   \ \mbox{for} \  i < j  \right \} \end{gather}
and the assignment:
\begin{gather} \label{the map} \bd  \{ \sigma \in \st(\mc T): \mc E \subset \sigma^{ss}\} \ni (\mc P,Z) & \rMapsto^{f_{\mc E}}  & \left (\{\abs{Z(E_i)} \}_{i=0}^n, \{ \phi^\sigma(E_i) \}_{i=0}^n\right )  \in \RR^{2(n+1)} \ed \end{gather}
restricted to $\Theta_{\mc E}'$ defines   a homeomorphism between $\Theta_{\mc E}'$ and $\RR_{>0}^{n+1}\times S^n(-\textbf{1},+\textbf{1})$ (as defined in Definition \ref{S(I)}). 

Assume now that   $\mc E=(E_0, E_1,\dots, E_n)$ is any  full exceptional collection in $\mc T$ (not restricted to be Ext).  If $\mc T$ is a triangulated category of finite type, then there are infinitely many choices of $\textbf{p}\in \ZZ^{n+1}$ such that $\mc E[\textbf{p}]=(E_0[p_0], E_1[p_1],\dots, E_n[p_n])$ is an  Ext-exceptional collection.  \cite[Lemma 3.19]{Macri} says that   the following  open subset of stability conditions is    connected and simply conected:
 \begin{gather}\label{theta_{mc E}}  \Theta_{{\mc E}}= \bigcup_{\left \{\textbf{p} \in \ZZ^{n+1} : {\mc E}[\textbf{p}] \ \mbox{is Ext} \right \}} \Theta_{{\mc E}[\textbf{p}]}'\subset \st(\mc T). \end{gather} 

 For the sake of completeness we will comment on  this set as well (compare with \cite[proof of Lemma 3.19]{Macri}). 

 By \cite[Corollary 3.20]{DK1} $\Theta_{{\mc E}}$ is the set of stability conditions $\sigma \in \st(\mc T)$ for which a shift of ${\mc E}$ is a $\sigma$-exceptional collection, in particular for each $\sigma \in \Theta_{{\mc E}}$ we have $\mc E \subset \sigma^{ss}$. Hence the assignment \eqref{the map} is well defined on   $\Theta_{\mc E}$. Furthermore, this defines a homeomorphism between $\Theta_{\mc E}$ and $f_{\mc E}(\Theta_{\mc E})$.  Indeed,  if $\mc E[\textbf{p}]$ is an Ext-collection for some $\textbf{p}\in \ZZ^{n+1}$, then  $f_{\mc E[\textbf{p}]}$ maps $\Theta_{{\mc E}[\textbf{p}]}'$ homeomorphically   to  $\RR_{>0}^{n+1}\times S^n(-\textbf{1},+\textbf{1})$ (see after \eqref{the map}) and  due to   $f_{\mc E[\textbf{p}]}-(0,\textbf{p})=f_{\mc E}$ we see that  ${f_{\mc E}}_{\vert \Theta_{{\mc E}[\textbf{p}]}'}$ is homeomorphism onto its image $\RR_{>0}^{n+1}\times( S^n(-\textbf{1},+\textbf{1}) - \textbf{p})$. Therefore, provided that  $f_{\mc E}$ is injective on $\Theta_{\mc E}$, \textbf{ the following restriction  is a  homeomorphism:}
\begin{gather} \label{A and homeo} {f_{\mc E}}_{\vert \Theta_{\mc E}}: \Theta_{{\mc E}} \ \ \rightarrow \ \     \RR_{>0}^{n+1} \times \left (  \bigcup_{ \textbf{p} \in A }  S^n(-\textbf{1},+\textbf{1}) - \textbf{p} \right ), \ \ \mbox{where}  \ \ A=\left \{\textbf{p} \in \ZZ^{n+1} : {\mc E}[\textbf{p}] \ \mbox{is Ext} \right \}. \end{gather}
 To show that the obtained function is injective,    assume that $\sigma_i=(\mc P_i, Z_i) \in \Theta_{\mc E}$, $i=1,2$ and $f_{\mc E}(\sigma_1)=f_{\mc E}(\sigma_2)$, i. e. $ \abs{Z_1(E_j)}=\abs{Z_2(E_j)}, \phi^{\sigma_1}(E_j)=\phi^{\sigma_2}(E_j)$  for all $j$, then by \eqref{Theta_mc E'} and the axiom $\phi^{\sigma}(E_j[p_j])=\phi^{\sigma}(E_j)+p_j$  we see that for any $\textbf{p}$ the incidence   $\sigma_1 \in \Theta_{\mc E[\textbf{p}]}'$ is equivaaent to  $\sigma_2 \in \Theta_{\mc E[\textbf{p}]}'$, hence by the injectivity of $f_{\mc E[\textbf{p}]}$ and $f_{\mc E[\textbf{p}]}-(0,\textbf{p})=f_{\mc E}$  we obtain $\sigma_1=\sigma_2$. Thus, we see that \eqref{A and homeo}  is a  homeomorphism.

 Finally, note that by \eqref{Theta_mc E'}, \eqref{the map}, \eqref{theta_{mc E}}, and $f_{\mc E[\textbf{p}]}=f_{\mc E}+(0,\textbf{p})$  one easily shows that 
\begin{gather} \label{Theta_E} \Theta_{\mc E} =\left  \{\sigma\in \st(\mc T): \mc E \subset \sigma^{ss}  \ \mbox{and} \  \phi^\sigma(\mc E) \in   \bigcup_{ \textbf{p} \in A }  S^n(-\textbf{1},+\textbf{1}) - \textbf{p}  \right \}. \end{gather}

\subsection{The set \texorpdfstring{$f_{\mc E}(\Theta_{\mc E})$}{\space} when \texorpdfstring{$n=2$}{\space}} \label{f_E(Theta_E)}   
Here is given  an explicit representation of $f_{\mc E}(\Theta_{\mc E})$, when $n=2$.  Remark \ref{for n geq 3} shows that the case  $n\geq 3$ is not completely analogous.  The only statement of this subsection, which will be used later is Proposition  \ref{lemma for f_E(Theta_E)}, the rest is its proof. 

Let us denote first:  \be B^n=\{(0,q_1 ,q_2,\dots, q_n)\subset \NN^{n+1} : 0\leq q_1 \leq q_2 \leq\dots\leq q_n \}\ee 
The following properties 		are clear from the definitions of $S^n(\mc J)$(Definition \ref{S(I)}) and of  $A\subset \ZZ^{n+1}$(formula \eqref{A and homeo})
\begin{gather} \label{prop1}  \forall \textbf{v}\in diag(\RR^{n+1}) \ \   \ \ S^n(\mc J) - \textbf{v} = S^n(\mc J) \\
                 \forall \textbf{v}\in diag(\ZZ^{n+1}) \ \   \ \  A-\textbf{v} =  A \\
                 \forall \textbf{v} \in B^n  \ \  \ \  A-\textbf{v} \subset  A.
 \end{gather}
Any $\textbf{p}=(p_0,p_1,\dots,p_n)\in A$ can be represented as  $\textbf{p}-(p_0,p_0,\dots,p_0) + (p_0,p_0,\dots,p_0)$, hence if we denote 
\be   \label{A_0}  A_0=\left \{\textbf{p} \in \ZZ^{n+1} : p_0=0, {\mc E}[\textbf{p}] \ \mbox{is Ext} \right \} \ee 
by the properties above we can write
\begin{gather} \label{f_E(Theta_E)1} \bigcup_{ \textbf{p} \in A }  S^n(-\textbf{1},+\textbf{1}) - \textbf{p}= \bigcup_{ \textbf{p} \in A_0 }  S^n(-\textbf{1},+\textbf{1}) - \textbf{p} =  \bigcup_{ \textbf{p} \in A_0 } \left  ( \bigcup_{ \textbf{v} \in B^n }  S^n(-\textbf{1},+\textbf{1}) + \textbf{v} \right ) -  \textbf{p}.
\end{gather}
For the cases $n=1,2$ we have the following simple form of the expression in the brackets:

\begin{lemma} \label{lemma about S^2(-1,1)} The following equalities hold:
\begin{gather}  \label{lemma about S^2(-1,1)eq}  \bigcup_{ \textbf{v} \in B^1 }  S^1(-1,+1) + \textbf{v}  = S^1(-\infty,1) \qquad   \bigcup_{ \textbf{v} \in B^2 }  S^2(-\textbf{1},+\textbf{1}) + \textbf{v}  = S^2(-\infty,\textbf{1}).\end{gather}
Recall  that  $S^n(-\infty,\textbf{1})=\{(y_0,y_1,\dots,y_n)\in \RR^{n+1} : y_i-y_j<1,  i<j  \}$ (see Definition \ref{S(I)}).
\end{lemma}
\begin{remark} \label{for n geq 3} For $n\geq 3$ we have not such an equality. For example, we  have $(0,-\frac{1}{2},\frac{1}{2}, 0,\dots,0) \in  S^n(-\infty,\textbf{1})$ but  $(0,-\frac{1}{2},\frac{1}{2}, 0,\dots,0) \not \in  \bigcup_{ \textbf{v} \in B^n }  S^n(-\textbf{1},+\textbf{1}) + \textbf{v} $ for $n\geq 3$.

 More precisely, it holds 
$\bigcup_{ \textbf{v} \in B^n }  S^n(-\textbf{1},+\textbf{1}) + \textbf{v} \subset  \not =  S^n(-\infty,\textbf{1})$ for $n\geq 3$.

\end{remark}
	\bpr(of Lemma \ref{lemma about S^2(-1,1)})   Note first that for any $\mc I = \{I_{ij}:  i <j \}$ as in Definition \ref{S(I)} and any $\textbf{p} \in \ZZ^{n+1}$  we have \be \label{shift of S(I)} S^n(\{I_{ij}:  i <j \})-\textbf{p}=S^n(\{I_{ij}-(p_i-p_j):  i <j \}). \ee  
	In particular for $n=1$ we have(now the index set of $\mc I$  has only one element: $(0,1)$):  \begin{gather} \bigcup_{ \textbf{v} \in B^1 }  S^1(-1,+1) + \textbf{v} = \bigcup_{ (0,k) \in \NN^2 }  S^1(-1,+1) + (0,k) = \bigcup_{k \in \NN } S^1(-1-k,1-k) \nonumber \\ = \bigcup_{k \in \NN } \{-1-k < y_0-y_1 < 1-k\}= \{ y_0-y_1 < 1\}=  S^1(-\infty ,+1) \nonumber \end{gather}
Using \eqref{prop1} and \eqref{shift of S(I)} one easily shows that:  
\begin{gather}   \bigcup_{ \textbf{v} \in B^2 }  S^2(-\textbf{1},+\textbf{1}) + \textbf{v} =  diag(\RR^{n+1}) \oplus  \{y_2=0\} \cap \left (\bigcup_{ \textbf{v} \in B^2 }  S^2(-\textbf{1},+\textbf{1}) + \textbf{v} \right ) \nonumber \\
 S^2(-\infty,\textbf{1}) = diag(\RR^{n+1}) \oplus  \{y_2=0\} \cap  S^2(-\infty,\textbf{1}). \nonumber \end{gather}

Obviously we have 
\begin{gather} \{y_2=0\} \cap  S^2(-\infty,\textbf{1}) = \{y_2=0\}\cap \left \{ \begin{array}{c} y_0 - y_1 <1 \\ y_0 - y_2 <1\\  y_1 - y_2 <1\end{array} \right \} = \left \{ \begin{array}{c} y_0 - y_1 <1 \\ y_0  <1\\  y_1  <1\end{array} \right \}.\nonumber \end{gather}
We will prove the second equality in \eqref{lemma about S^2(-1,1)eq} by showing that:
\begin{gather}\label{lemma about S^2(-1,1)eq1} \{y_2=0\} \cap \left (\bigcup_{ \textbf{v} \in B^2 }  S^2(-\textbf{1},+\textbf{1}) + \textbf{v} \right )= \left \{ \begin{array}{c} y_0 - y_1 <1 \\ y_0  <1\\  y_1  <1\end{array} \right \}. \end{gather}
Let $(0,k,k+l) \in B^2$, $k$, $l\in \NN$ be a vector in  $B^2$. By \eqref{shift of S(I)} we have: 
\begin{gather} S^2(-\textbf{1},\textbf{1})+(0,k,k+l)=\left \{ \begin{array}{c}-1-k < y_0-y_1 < 1-k \\ 
-1-k-l < y_0-y_2 < 1-k - l \\ -1-l < y_1-y_2 < 1- l  \end{array} \right \}\subset \left \{ \begin{array}{c} y_0-y_1 < 1 \\ 
 y_0-y_2 < 1 \\  y_1-y_2 < 1  \end{array} \right \}. \end{gather}
Denoting  the unit open square by $C(-1,+1)=\{\abs{y_i}<1;i=0,1\}\subset \RR^2$, we can write: 
\begin{gather} \{y_2=0 \}\cap \left (S^2(-\textbf{1},\textbf{1})+(0,k,k+l) \right )=\left \{ \begin{array}{c}-1-k < y_0-y_1 < 1-k\nonumber  \\ 
-1-k-l < y_0 < 1-k - l \\ -1-l < y_1 < 1- l  \end{array} \right \}\nonumber \\ 
=S^1(-1-k,+1-k)\cap \left ( C(-1,+1) -(k+l,l)\right )\nonumber \\
=\left ( S^1(-1,+1)+(0,k) \right ) \cap \left ( C(-1,+1) -(k+l,l)\right ) \nonumber \\
=\left ( S^1(-1,+1)-(k+l,k+l)+(0,k) \right ) \cap \left ( C(-1,+1)  -(k+l,l)\right ) \nonumber \\ 
= \left ( S^1(-1,+1) \cap C(-1,+1)\right )-(k+l,l). \nonumber\end{gather}	
Therefore:
\begin{gather}\label{lemma about S^2(-1,1)1} \{y_2=0 \}\cap \left (\bigcup_{\textbf{v}\in B^2}S^2(-\textbf{1},\textbf{1})+\textbf{v} \right )= \bigcup_{k\in \NN} \left ( \bigcup_{l\in \NN}   \left ( S^1(-1,1) \cap C(-1,1)\right )-(l,l) \right ) - (k,0).\end{gather}

Before we continue with the proof of Lemma \ref{lemma about S^2(-1,1)}, we prove:
\begin{lemma} For any $k\in \ZZ\cup\{+\infty\}$  we have the following equality: \begin{gather} \label{eq for S(-1,1)} 
\bigcup_{l\leq k}    S^1(-1,+1) \cap C(-1,+1) +(l,l) =  S^1(-1,+1) \cap \left \{\begin{array}{c} y_0 <1+k \\ y_1 <1+k  \end{array}\right \} . \end{gather}
\end{lemma}
\bpr We show first the equality for $k=+\infty$. Let $(a_0, a_1) \in S^1(-1,+1)$, i. e. $\abs{a_0-a_1}<1$.   Since $\RR=\bigcup_{l\in \ZZ} [2l-1,2l+1)$, there exists $l\in \ZZ$ such that $a_0+a_1 \in [2 l -1,  2 l +1)$, i. e. $-1 \leq a_0+a_1 - 2 l \leq 1$.  We have also     $-1<a_0-a_1<+1$ and  due to the equalities:
\begin{gather} a_0 -l  =\frac{ a_0+a_1-2l }{2} + \frac{ a_0-a_1}{2}; \ \  a_1 -l  =\frac{ a_0+a_1-2l }{2} + \frac{ a_1-a_0}{2}\nonumber\end{gather}
 we obtain $-1 =-\frac{1}{2}-\frac{1}{2}< a_i -l  < \frac{1}{2}+\frac{1}{2}=1$ for  $i=0,1$.  Hence $(a_0,a_1) -(l,l)\in C(-1,+1) \cap S(-1,+1)$, and we proved the equality \eqref{eq for S(-1,1)}  with $k=+\infty$. By \eqref{prop1}  and since the translation in $\RR^2$ is bijective  we  rewrite this equality as follows   $  S^1(-1,+1)= \bigcup_{l\in \ZZ}    S^1(-1,+1) \cap C(-1,+1) +(l,l) = \bigcup  \left (  S^1(-1,+1) +(l,l)\right ) \cap\left ( C(-1,+1) +(l,l) \right )=  S^1(-1,+1) \cap \left ( \bigcup_{l\in \ZZ}    C(-1,+1) +(l,l) \right )$.  Hence
\begin{gather} \label{some shifts of square} S^1(-1,1) \cap \left \{\begin{array}{c} y_0 <1+k \\ y_1 <1+k  \end{array}\right \} =  S^1(-1,1) \cap \left ( \bigcup_{l\in \ZZ}    C(-1,1) +(l,l) \right ) \cap \left \{\begin{array}{c} y_0 <1+k \\ y_1 <1+k  \end{array}\right \}. \end{gather}
Due to the equalities   
\begin{gather}  \left \{\begin{array}{c} y_0 <1+k \\ y_1 <1+k  \end{array}\right \} \cap \left ( C(-1,1) +(l,l) \right )  = \left \{ \begin{array}{c c}  \emptyset & \ \mbox{if} \ \ l\geq k+2  \\ 
\left \{\begin{array}{c} k< y_0 <1+k \\ k < y_1 <1+k  \end{array}\right \} \subset  C(-1,1) +(k,k) & \ \mbox{if} \ \ l= k+1 \\ 
 C(-1,+1) +(l,l) & \ \mbox{if} \  \ l\leq  k  \end{array} \right. \nonumber \end{gather}
we obtain  $ \left ( \bigcup_{l\in \ZZ}    C(-1,1) +(l,l) \right ) \cap \left \{\begin{array}{c} y_0 <1+k \\ y_1 <1+k  \end{array}\right \} =  \bigcup_{l\leq k}    C(-1,1) +(l,l)  $. By \eqref{some shifts of square} and applying again \eqref{prop1} we obtain the equality \eqref{eq for S(-1,1)}   for  $k\in \ZZ$.
\epr
Now we put \eqref{eq for S(-1,1)} with $k=0$  in \eqref{lemma about S^2(-1,1)1} and obtain 
\begin{gather}\label{lemma about S^2(-1,1)2} \{y_2=0 \}\cap \left (\bigcup_{\textbf{v}\in B^2}S^2(-\textbf{1},\textbf{1})+\textbf{v} \right )= \bigcup_{k\in \NN} \left ( S^1(-1,+1) \cap \left \{\begin{array}{c} y_0 <1 \\ y_1 <1  \end{array}\right \} \right  )- (k,0).\end{gather}
The next step is to show that 
\begin{gather}\label{lemma about S^2(-1,1)3} \bigcup_{k\in \NN} \left ( S^1(-1,+1) \cap \left \{\begin{array}{c} y_0 <1 \\ y_1 <1  \end{array}\right \} \right  )- (k,0)= \bigcup_{k\in \NN} \left ( S^1(-1,+1) - (k,0) \right  ) \cap \left \{\begin{array}{c} y_0 <1 \\ y_1 <1  \end{array}\right \}.\end{gather}
The inclusion $\subset$ is clear.  Assume now that $a_0, a_1 \in \RR$,  $k\in \NN$ and $\abs{a_0-a_1} <1$ and $a_0-k <1$, $a_1 <1$.  We have to find $a'_0 \in \RR$, and $k' \in \NN$ such that \be \abs{a'_0-a_1}<1 \qquad a'_0 <1 \qquad  a'_0 - k'= a_0-k.  \ee
First note that $a_0 = a_0- a_1 + a_1 < \abs{a_0- a_1 } +  a_1 < 2$.
If $k=0$ or $a_0 <1$, then we put $a'_0 = a_0$, $k'=k$.  
Thus, we can assume that  $k\geq 1$ and $1 \leq a_0 < 2$. Now $a_1 <1$ and  $\abs{a_0-a_1}<1$ imply $0\leq a_1 < 1$.  It follows  that $ -1 < -a_1\leq  a_0-1-a_1 <1 $,  therefore we can put  $a'_0=a_0-1$, $k'=k-1$. Hence we obtain \eqref{lemma about S^2(-1,1)3}. 

On the other hand by \eqref{prop1} and the already proven first equality in \eqref{lemma about S^2(-1,1)eq}  we have   
$$\bigcup_{k\in \NN} S^1(-1,1)-(k,0) = \bigcup_{k\in \NN}S^1(-1,1)+(0,k) =S^1(-\infty,1).$$ The latter equality and equalities   \eqref{lemma about S^2(-1,1)2},  \eqref{lemma about S^2(-1,1)3} imply  \eqref{lemma about S^2(-1,1)eq1} and the lemma follows.
 \epr

Putting \eqref{lemma about S^2(-1,1)eq}  in \eqref{f_E(Theta_E)1} and then using \eqref{shift of S(I)} we obtain for the case $n=2$:
\begin{gather} \label{f_E(Theta_E)2} \bigcup_{ \textbf{p} \in A }  S^2(-\textbf{1},+\textbf{1}) - \textbf{p}= \bigcup_{ \textbf{p} \in A_0 } S^2(-\infty, \textbf{1}) -  \textbf{p}= \bigcup_{ (0,p_1,p_2) \in A_0 } \left \{ \begin{array}{c} y_0 - y_1 < 1+p_1 \\  y_0 - y_2 < 1+p_2 \\ y_1 - y_2 < 1+p_2-p_1\end{array}\right \}.
\end{gather}
Using  the equality \eqref{f_E(Theta_E)2},  the homeomorphism \eqref{A and homeo}, and  \eqref{Theta_E}  we will prove the main result of this subsection:
\begin{prop} \label{lemma for f_E(Theta_E)} Let  $\mc T$ be a $k$-linear triangulated category. Let $\mc E = (A_0,A_1,A_2)$ be a full exceptional collection, such  that: \begin{gather} 1+\alpha = \min \{i:\hom^i(A_0,A_1)\neq 0\}\in \ZZ \nonumber \\ \label{alpha,beta,gamma} 1+\beta = \min \{i:\hom^i(A_0,A_2)\neq 0\}\in \ZZ  \\ 1+\gamma = \min \{i:\hom^i(A_1,A_2)\neq 0\}\in \ZZ . \nonumber  \end{gather}
Then the subset $ \Theta_{\mc E}\subset \st(\mc T)$ defined in  \eqref{theta_{mc E}} has the following description:  \begin{gather} \label{Theta_E n=2} \Theta_{\mc E} =\left  \{\sigma\in \st(\mc T): \mc E \subset \sigma^{ss}  \ \mbox{and} \  \begin{array}{l} \phi^\sigma(A_0) -\phi^\sigma(A_1) < 1+\alpha \\ \phi^\sigma(A_0) -\phi^\sigma(A_2) < 1+\min\{\beta, \alpha+\gamma\} \\\phi^\sigma(A_1) -\phi^\sigma(A_2) < 1+\gamma\end{array} \right \}  \end{gather} 
 and  $\Theta_{\mc E}$ is homeomorphic with  the set  $ \RR_{>0}^3 \times   \left \{  \begin{array}{c} y_0 - y_1 < 1+\alpha \\  y_0 - y_2 < 1+\min\{\beta, \alpha+\gamma\} \\ y_1 - y_2 < 1+\gamma\end{array} \right \}$ by the map  $f_{\mc E}$ in  \eqref{the map}  restricted to $\Theta_{\mc E}$. In particular $\Theta_{\mc E}$ is contractible. 
\end{prop}
\bpr Given a family  $\mc I$ of  the form: $\mc I$ $=$ $ \{ I_{01}=(-\infty,u), I_{02}=(-\infty,v),   I_{12}=(-\infty,w) \}$,  we  write $S\left (\begin{array}{c} -\infty,u \\  -\infty,v \\  -\infty,w  \end{array} \right ) $  for $S^2(\mc I) $ throughout the proof. By \eqref{f_E(Theta_E)2},  \eqref{Theta_E}, and \eqref{A and homeo}  the proof is reduced to showing that:    
\begin{gather} \label{lemma for f_E(Theta_E)1} \bigcup_{ (0,p_1,p_2) \in A_0 }S\left (\begin{array}{c} -\infty,1+p_1 \\  -\infty,1+p_2 \\  -\infty,1+p_2-p_1  \end{array} \right ) =  S\left (\begin{array}{c} -\infty,1+\alpha \\  -\infty,1+\min\{\beta, \alpha+\gamma\} \\  -\infty,1+\gamma  \end{array} \right ) .     \end{gather}
From the definition of $A_0$ in \eqref{A_0} and the definition of $\alpha$, $\beta$, $\gamma$ one easily obtains:
\begin{gather} \label{p_1 leq alpha} (0,p_1,p_2) \in  A_0 \ \ \Rightarrow \ \ p_1 \leq \alpha, p_2 \leq \min\{\beta,\alpha+\gamma \}; \qquad (0,\alpha, \min\{\beta, \alpha+\gamma\}) \in A_0. \end{gather}
If $u\leq u'$,  $v\leq v'$, $w\leq w'$, then $S\left ( -\infty, (u,v ,w )  \right )  \subset S\left ( -\infty, (u',v' ,w' )  \right )$, hence by \eqref{p_1 leq alpha}  we have:
\begin{gather} \bigcup_{ (0,p_1,p_2) \in A_0 }S\left (\begin{array}{c} -\infty,1+p_1 \\  -\infty,1+p_2 \\  -\infty,1+p_2-p_1  \end{array} \right )=S\left (\begin{array}{c} -\infty,1+\alpha \\  -\infty,1+\min\{\beta, \alpha+\gamma\} \\ -\infty,1+\min\{\beta, \alpha+\gamma\}-\alpha  \end{array} \right ) \cup \nonumber \\[-2mm] \label{lemma for f_E(Theta_E)2} \\[-2mm] \bigcup_{\left \{ \begin{array}{c} (0,p_1,p_2)\in A_0: \\ p_2-p_1 > \min\{\beta, \alpha+\gamma\}-\alpha \end{array} \right \}} S\left (\begin{array}{c} -\infty,1+p_1 \\  -\infty,1+p_2 \\  -\infty,1+p_2-p_1  \end{array} \right ). \nonumber  \end{gather}
Now we consider two cases.

\ul{If $\min\{\beta, \alpha+\gamma\}= \alpha + \gamma$,} then $\min\{\beta, \alpha+\gamma\}-\alpha=\gamma$ and  $(A_0,A_1[p_1], A_2[p_2])$ is not an Ext-collection for $p_2-p_1 >\gamma$ (since $\hom^{p_1+\gamma+1-p_2}(A_1[p_1],A_2[p_2]) \neq 0 $, $p_2-p_1-\gamma-1\geq 0$), hence the equality \eqref{lemma for f_E(Theta_E)2} reduces to \eqref{lemma for f_E(Theta_E)1}.

\ul{If $\min\{\beta, \alpha+\gamma\}= \beta < \alpha + \gamma$,} then   $\beta \leq \alpha -i + \gamma$ for $i \leq \alpha + \gamma - \beta$  and hence
\be \{(0,\alpha-i, \beta): 0\leq i \leq \alpha + \gamma - \beta \} \subset A_0. \ee
 Furthermore, we claim that the equality \eqref{lemma for f_E(Theta_E)2} reduces to 
\begin{gather} \label{lemma for f_E(Theta_E)3} \bigcup_{ (0,p_1,p_2) \in A_0 }S\left (\begin{array}{c} -\infty,1+p_1 \\  -\infty,1+p_2 \\  -\infty,1+p_2-p_1  \end{array} \right )=   \bigcup_{i=0}^{\alpha + \gamma - \beta} S\left (\begin{array}{c} -\infty,1+\alpha -i \\  -\infty,1+\beta \\  -\infty,1+\beta-\alpha+i  \end{array} \right ).   \end{gather}
Indeed,  the first set of the union in \eqref{lemma for f_E(Theta_E)2} is the same as the first set of the union  \eqref{lemma for f_E(Theta_E)3}. Now assume that $(0,p_1,p_2) \in A_0$ and $p_2-p_1>\beta -\alpha$, then $\beta -\alpha < p_2-p_1 \leq \gamma$. Therefore for some $1\leq i \leq \gamma+\alpha-\beta$ we have $ p_2-p_1 =\beta -\alpha + i $.
From  \eqref{p_1 leq alpha} we have also  $ p_2 \leq \beta $, therefore $ p_1 =p_2-\beta +\alpha - i\leq \alpha - i$, and then $S\left (\begin{array}{c} -\infty,1+p_1 \\  -\infty,1+p_2 \\  -\infty,1+p_2-p_1  \end{array} \right )\subset S\left (\begin{array}{c} -\infty,1+\alpha -i \\  -\infty,1+\beta \\  -\infty,1+\beta-\alpha+i  \end{array} \right )$ and we showed \eqref{lemma for f_E(Theta_E)3}. The last step  of the proof is to show that
\begin{gather} \label{lemma for f_E(Theta_E)4}  \bigcup_{i=0}^{\alpha + \gamma - \beta}S\left (\begin{array}{c} -\infty,1+\alpha -i \\  -\infty,1+\beta \\  -\infty,1+\beta-\alpha+i  \end{array} \right )=  S\left (\begin{array}{c} -\infty,1+\alpha  \\  -\infty,1+\beta \\  -\infty,1+\gamma  \end{array} \right ).   \end{gather}
The inclusion $\subset $ is clear. To show the inclusion $\supset$,  assume that $(a_0,a_1,a_2) \in \RR^3$ and \\ $a_0 - a_1 < 1 +\alpha$, $a_0 - a_2 < 1 +\beta$,   $a_1 - a_2 < 1 +\gamma$. 

 If $a_0 - a_1 < 1+\alpha -(\alpha + \gamma -\beta) = 1+\beta - \gamma$, then by   $a_1 - a_2 < 1 +\gamma$ it follows that  $(a_0,a_1,a_2)$ is in the set with index  $i=\alpha + \gamma -\beta$ on the right-hand side. 

It remains to consider the case, when  $ 1+\alpha -i > a_0-a_1 \geq 1 + \alpha -i -1  $ for some $0\leq i<\alpha + \gamma -\beta$.  Now  $(a_0,a_1,a_2)$ is in the set  indexed by the given    $i$. Indeed, now  $   a_1-a_0 \leq i-\alpha $ and by $a_0-a_2 < 1 + \beta $ we have $a_1-a_2 = a_1-a_0 + a_0-a_2<1+\beta  +i -\alpha$. 
\epr 
\subsection{More propositions used  for gluing} Since we will often use the notion of a $\sigma$-triple, for the sake of completeness we rewrite here  \cite[Definition 3.19]{DK1} for triples (see also \cite[Remark 3.31]{DK1}):
\begin{df} \label{sigma triple} An exceptional triple $(A_0, A_1,A_2)$ is a $\sigma$-triple iff  the following conditions hold: {\rm \textbf{(a)}} $\hom^{\leq 0}(A_i, A_j)=0$ for $i\neq j$;  {\rm \textbf{(b)}} $ \{A_i\}_{i=0}^2\subset \sigma^{ss} $ ;  {\rm \textbf{(c)}} $ \{\phi(A_i)\}_{i=0}^2\subset (t,t+1) $ for some $t\in \RR$.
\end{df}
We   enhance now  Proposition \ref{all exc semistable} for the case $n=2$: 
\begin{prop} \label{phi_1>phi_2} Let  $\mc T$ be a $k$-linear triangulated category. Let $\mc E = (A_0,A_1,A_2)$,  $\alpha$,  $\beta$, $\gamma$  be as in Proposition \ref{lemma for f_E(Theta_E)}. Let $\sigma \in \Theta_{\mc E}$ (hence we have the inequalities in \eqref{Theta_E n=2}).  

{\rm (a)} If   $\phi^\sigma(A_0)\geq  \phi^\sigma(A_1[\alpha])$, then $\mc A \cap \mc T_{exc}\subset \sigma^{ss}$, where $\mc A$ is the extension closure of $(A_0,A_1[\alpha])$.

{\rm (b) } If  $\phi^\sigma(A_1)\geq  \phi^\sigma(A_2[\gamma])$, then $\mc A \cap \mc T_{exc}\subset \sigma^{ss}$, where $\mc A$ is the extension closure of $(A_1,A_2[\gamma])$.
\end{prop}
\bpr If an equality holds in (a) or (b), then we have $\mc A \subset \mc P(t)$ for some $t\in \RR$ and the Proposition follows. Hence we can assume that we have a proper inequality in both the cases.
 
 (a)  By the definition of $\Theta_{\mc E}$ in \eqref{theta_{mc E}} and \cite[Corollary 3.20]{DK1} we see that $(A_0[l], A_1[i], A_2[j])$ is a $\sigma$-triple for some $l,i,j \in \ZZ$. We can assume\footnote{note that $(A_0, A_1[i], A_2[j])$ is a $\sigma$-triple iff  $(A_0[k], A_1[i+k], A_2[j+k])$ is a $\sigma$-triple}  $l=0$ and then   $\hom^{\leq 0}(A_0,A_1[i])=0$ and $\abs{\phi(A_0)- \phi(A_1[i])}<1$. From the definition of $\alpha$ we see that $i\leq \alpha$. Actually we must have $i=\alpha$, otherwise  the given  inequality $\phi(A_0)- \phi(A_1[\alpha])>0$ implies $\phi(A_0)- \phi(A_1[i])>1$, which is a contradiction. Thus $(A_0,A_1[\alpha],A_2[j])$ is a $\sigma$-triple for some $j\in \ZZ$. Now we apply Proposition \ref{all exc semistable}. 

(b)  In this case we shift the given triple to a $\sigma$-triple of the form  $(A_0[l], A_1, A_2[j])$  for some $l,j \in \ZZ$, in particular we have  $\hom^{\leq 0}(A_1,A_2[j])=0$ and $\abs{\phi(A_1)- \phi(A_2[j])}<1$. From the definition of $\gamma$ and the given  inequality $\phi(A_1)- \phi(A_2[\gamma])>0$  it follows  that $j= \gamma$. Thus $(A_0[l],A_1,A_2[\gamma])$ is a $\sigma$-triple for some $l\in \ZZ$. Now we apply Proposition \ref{all exc semistable}.
\epr
 
\begin{prop} \label{mutations} Let $\mc T$ has the property that  for each exceptional triple $(A_0,A_1,A_2)$ and any two $0\leq i<j\leq 2$ there exists unique $k\in \ZZ$ satisfying $\hom^k(A_i,A_j)\neq 0$.    Let $\mc E = (A_0,A_1,A_2)$ be a full exceptional collection in $\mk{\mc T}$. 

 Let $R_0(\mc E) = (A_1,R_{A_1}(A_0),A_2)$, $L_0(\mc E) = (L_{A_0}(A_1), A_0, A_2)$, $R_1(\mc E) = (A_0, A_2, R_{A_2}(A_1))$, $L_1(\mc E) = (A_0, L_{A_1}(A_2), A_1)$    be the triples obtained by a single   mutation applied to $\mc E$.\footnote{ Recall that for any exceptional pair $(A,B)$  the exceptional objects $L_A(B)$ and $R_B(A)$ are determined by the triangles $\bd L_A(B)&\rTo &\Hom^*(A,B)\otimes A & \rTo^{ev^*_{A,B}} & B \ed$; \  $\bd A &\rTo^{coev^*_{A,B}} &\Hom^*(A,B)^{\check{}}\otimes B &\rTo & R_B(A) \ed$ and that $(L_A(B),A)$, $(B,R_B(A))$ are exceptional pairs.}   
Then the four intersections   $\Theta_{\mc E}\cap \Theta_{R_0(\mc E)}$, $\Theta_{\mc E}\cap \Theta_{L_0(\mc E)}$, $\Theta_{\mc E}\cap \Theta_{R_1(\mc E)}$, $\Theta_{\mc E}\cap \Theta_{L_1(\mc E)}$ are all contractible and non-empty. 
\end{prop}
\bpr    Since  $\mc E \sim \mc E'$ implies  $\Theta_{\mc E}= \Theta_{\mc E'}$, $R_i(\mc E)\sim R_i(\mc E')$, $L_i(\mc E)\sim L_i(\mc E')$,  we can assume that $l=\hom^1(A_0,A_1)> 0$, $p = \hom^1(A_1,A_2)> 0$.  By the assumptions on $\mc T$  the other degrees are zero  and  it follows that the integers  $\alpha$, $\gamma$ defined in \eqref{alpha,beta,gamma}  vanish and  from Proposition \ref{lemma for f_E(Theta_E)}  we  get: 
\begin{gather} \label{Theta_E n=2 1} \Theta_{\mc E} =\left  \{\sigma\in \st(\mc T): \mc E \subset \sigma^{ss}  \ \mbox{and} \  \begin{array}{l} \phi^\sigma(A_0) -\phi^\sigma(A_1) < 1 \\ \phi^\sigma(A_0) -\phi^\sigma(A_2) < 1+\min\{\beta,0\} \\\phi^\sigma(A_1) -\phi^\sigma(A_2) < 1 \end{array} \right \}.  \end{gather}

We start with the  intersection $\Theta_{\mc E}\cap \Theta_{R_0(\mc E)}$. 
Let us denote $X = R_{A_1}(A_0)[-1]$. Let  $\alpha'$, $\beta'$, $\gamma'$ be the integers corresponding  to the triple $(A_1,X,A_2)$ used in Proposition  \ref{lemma for f_E(Theta_E)}.  We have $1+\beta'=\min\{k : \hom^k(A_1,A_2)\neq 0\}=1+\gamma=1$, hence $\beta'=0$.   On the other hand     from the definition of $R_{A_1}(A_0)$ we have a triangle
 \begin{gather} \label{one triangle} A_1^{\oplus l} \rightarrow   X  \rightarrow  A_0 \rightarrow  A_1^{\oplus l}[1]\end{gather}
and it follows that   $\hom(A_1,X)\neq 0$, hence  $\alpha'=-1$.  We apply Proposition \ref{lemma for f_E(Theta_E)} to the triple  $(A_1,X,A_2)$ and obtain (note that $1+\min\{\beta', \alpha'+\gamma'\}=1+\min\{0, \gamma'-1\}=\min\{1, \gamma'\}$)  \begin{gather} \label{Theta_E n=2 2}  \Theta_{R_0(\mc E)}=\Theta_{(A_1,X,A_2)} =\left  \{\sigma\in \st(\mc T): \begin{array}{c}A_1\in \sigma^{ss} \\ X \in \sigma^{ss}\\ A_2\in \sigma^{ss}\end{array} \ \mbox{and} \  \begin{array}{l} \phi^\sigma(A_1) -\phi^\sigma(X) < 0 \\ \phi^\sigma(A_1) -\phi^\sigma(A_2)  < \min\{1, \gamma'\}\\ \phi^\sigma(X) -\phi^\sigma(A_2) < 1+\gamma'\end{array} \right \}.  \end{gather} 
From the defintion of $\beta, \gamma$ we have  $0=\hom^{\leq \min\{\beta,\gamma\}}(A_0,A_2)=\hom^{\leq \min\{\beta,\gamma\}}(A_1,A_2)$, and  then the triangle \eqref{one triangle} implies  that $\hom^{\leq \min\{\beta,\gamma\}}(X,A_2)=0$, it follows that
\begin{gather} \label{min beta, gamma leq gamma'} \min\{\beta, \gamma\}= \min\{\beta, 0\} \leq \gamma'. \end{gather}

Assume that $\sigma\in  \Theta_{(A_1,X,A_2)}\cap \Theta_{\mc E}$. Then   $A_0,A_1,A_2,X$ are all semi-stable and $\phi(A_1)<\phi(X)$.\footnote{We omit  sometimes  the superscript $\sigma$ in expressions like  $\phi^\sigma(X)$  and write  just $\phi(X)$.}  It is easy to show that $\hom(X,A_0)\neq 0$ (using the triangle \eqref{one triangle}), hence $\phi(X)\leq \phi(A_0)$ and therefore  $\phi(A_1)<\phi(A_0)$, and  we obtain the inclusion $\subset$ in the following formula (the third inequality in this formula is  the second in \eqref{Theta_E n=2 2}, the other inequalities are in \eqref{Theta_E n=2 1} together with $\phi(A_1)<\phi(A_0)$) 
 \begin{gather} \label{Theta_E n=2 3} \Theta_{\mc E}\cap  \Theta_{R_0(\mc E)} =\left  \{\sigma\in \st(\mc T): \mc E \subset \sigma^{ss}  \ \mbox{and} \  \begin{array}{l} 0<\phi^\sigma(A_0) -\phi^\sigma(A_1) < 1 \\ \phi^\sigma(A_0) -\phi^\sigma(A_2) < 1+\min\{\beta, 0\} \\ \phi^\sigma(A_1) -\phi^\sigma(A_2) < \min\{\gamma',1\}\end{array} \right \}.  \end{gather} 
We  show now the inclusion $\supset$. Assume that $\mc E \subset \sigma^{ss}$ and that the inequalities on the right hand side of \eqref{Theta_E n=2 3} hold. In particular the inequalities in \eqref{Theta_E n=2 1} hold, hence we have     $\sigma \in \Theta_{\mc E}$ and $\phi^\sigma(A_0) >\phi^\sigma(A_1)$.  Proposition \ref{phi_1>phi_2} (a) ensures $X\in \sigma^{ss}$ and by \eqref{one triangle} we get $\hom(A_1,X)\neq 0$, $\hom(X,A_2)\neq 0$, hence 
\be \label{Z(X)=l Z(A_1) + Z(A_0)2} X \in \sigma^{ss} \quad \phi(A_1)\leq \phi(X)\leq \phi(A_0) \quad  Z(X)=l Z(A_1) + Z(A_0). \ee
Using \eqref{phase formula} and  $ 0<\phi^\sigma(A_0) -\phi^\sigma(A_1) < 1$  we see that  $Z(A_1),  Z(A_0)$ are not collinear(see Definition \ref{collinear}), therefore  $Z(X)=l Z(A_1) + Z(A_0)$ is collinear neither with $ Z(A_1)$ nor with   $Z(A_0)$.  Now we apply \eqref{phase formula} again  and by \eqref{Z(X)=l Z(A_1) + Z(A_0)2}  we obtain $ \phi(A_1)< \phi(X)< \phi(A_0) $. In particular, we obtain the first inequality in \eqref{Theta_E n=2 2}. The second inequality in \eqref{Theta_E n=2 2} is the same as the third inequality of \eqref{Theta_E n=2 3}.     From $\phi^\sigma(A_0) -\phi^\sigma(A_2) < 1+\min\{\beta, 0\} $   and \eqref{min beta, gamma leq gamma'} we get  $\phi(X)-\phi^\sigma(A_2) < \phi^\sigma(A_0) -\phi^\sigma(A_2) < 1+\gamma' $, hence the third  inequality in \eqref{Theta_E n=2 2} is verified also. Thus we showed \eqref{Theta_E n=2 3}. This equality implies that the set $\Theta_{\mc E}\cap  \Theta_{R_0(\mc E)}$ is contractible. Indeed,we have a homeomorphism ${f_{\mc E}}_{\vert \Theta_{\mc E}}:\Theta_{\mc E} \rightarrow f_{\mc E}(\Theta_{\mc E})$ (see \eqref{A and homeo}, \eqref{the map}). The proved equality \eqref{Theta_E n=2 3} shows that: 
$ f_{\mc E}\left (  \Theta_{\mc E}\cap  \Theta_{R_0(\mc E)}\right )= \RR_{>0}^3 \times \left \{   \begin{array}{c} 0 <\phi_0 - \phi_1 <1 \\  \phi_0 -\phi_2 < 1+\min\{\beta, 0\} \\ \phi_1 -\phi_2 < \min\{\gamma',1\}\end{array}\right\} $,
hence  $\Theta_{\mc E}\cap  \Theta_{R_0(\mc E)}$ is contractible.

 Next, we consider the  intersection $\Theta_{\mc E}\cap \Theta_{L_1(\mc E)}$, where $L_1(\mc E)=(A_0, L_{A_1}(A_2), A_1)$.

Let us denote $Y = L_{A_1}(A_2)[1]$. Let  $\alpha'$, $\beta'$, $\gamma'$ be the integers corresponding  to the triple $(A_0,Y,A_1)$.  Obviously $\beta'=\alpha=0$.   From the definition of $L_{A_1}(A_2)$ we have a triangle  
\begin{gather} \label{another triangle} A_2 \rightarrow   Y  \rightarrow  A_1^{\oplus p} \rightarrow A_2[1]\end{gather}
and it follows that   $\hom(Y,A_1)\neq 0$, hence  $\gamma'=-1$.   Proposition \ref{lemma for f_E(Theta_E)} applied  to the triple  $(A_0,Y,A_1)$ results in the equality(note that $1+\min\{0',\alpha'-1\}=\min\{1,\alpha'\}$)  \begin{gather} \label{Theta_E n=2 5} \Theta_{L_1(\mc E)}=\Theta_{(A_0,Y,A_1)} =\left  \{\sigma\in \st(\mc T): \begin{array}{c}A_0\in \sigma^{ss} \\ Y \in \sigma^{ss}\\ A_1\in \sigma^{ss}\end{array} \ \mbox{and} \  \begin{array}{l} \phi^\sigma(A_0) -\phi^\sigma(Y) < 1+\alpha' \\ \phi^\sigma(A_0) -\phi^\sigma(A_1)  < \min\{1,\alpha'\} \\ \phi^\sigma(Y) -\phi^\sigma(A_1) < 0 \end{array} \right \}.  \end{gather} 
From the defintion of $\alpha, \beta$ for the initial sequence $\mc E$ we have  $0=\hom^{\leq \min\{\alpha,\beta\}}(A_0,A_1)$,\\ $0=\hom^{\leq \min\{\alpha,\beta\}}(A_0,A_2)$, and  then the triangle \eqref{another triangle} implies  that $\hom^{\leq \min\{\alpha,\beta\}}(A_0,Y)=0$, it follows that
\begin{gather} \label{min beta, gamma leq gamma'1} \min\{\alpha,\beta\}= \min\{0,\beta\} \leq \alpha'. \end{gather}

Assume that $\sigma\in  \Theta_{(A_0,Y,A_1)}\cap \Theta_{\mc E}$. Then   $A_0,A_1,A_2,Y$ are all semi-stable and by \eqref{Theta_E n=2 5} $\phi(Y)<\phi(A_1)$.  The triangle \eqref{another triangle} implies  $\hom(A_2,Y)\neq 0$, hence  $\phi(A_2)\leq \phi(Y) <\phi(A_1)$. Combining this inequality with the inequalities in \eqref{Theta_E n=2 5}, \eqref{Theta_E n=2 1}  we obtain the inclusion $\subset$ in the following formula: 
 \begin{gather} \label{Theta_E n=2 6} \Theta_{\mc E}\cap  \Theta_{L_1(\mc E)} =\left  \{\sigma\in \st(\mc T): \mc E \subset \sigma^{ss}  \ \mbox{and} \  \begin{array}{l} \phi^\sigma(A_0) -\phi^\sigma(A_1) < \min\{1,\alpha'\} \\ \phi^\sigma(A_0) -\phi^\sigma(A_2) < 1+\min\{\beta,0\} \\0<\phi^\sigma(A_1) -\phi^\sigma(A_2) < 1 \end{array} \right \}.  \end{gather} 
To show the inclusion $\supset$, assume that $\mc E \subset \sigma^{ss}$ and that the inequalities on the right hand side of \eqref{Theta_E n=2 6} hold. In particular, we have     $\sigma \in \Theta_{\mc E}$(see \eqref{Theta_E n=2 1}). It remains to show that  $Y\in\sigma^{ss}$ and that the inequalities in  \eqref{Theta_E n=2 5} hold.  From Proposition \ref{phi_1>phi_2} (b) and   $\sigma \in \Theta_{\mc E}$, $\phi^\sigma(A_1) >\phi^\sigma(A_2)$ we obtain  $Y \in \sigma^{ss}$.   The triangle \eqref{another triangle} implies 
\be \label{Z(Y)= Z(A_1) +p  Z(A_2)} Y \in \sigma^{ss} \quad \phi(A_2)\leq \phi(Y)\leq \phi(A_1) \quad  Z(Y)=p Z(A_1) +  Z(A_2). \ee
By similar arguments as in the previous case, using   \eqref{phase formula},   $ 0<\phi^\sigma(A_1) -\phi^\sigma(A_2) < 1$ and \eqref{Z(Y)= Z(A_1) +p  Z(A_2)} one shows that  $ \phi(A_2) < \phi(Y) < \phi(A_1) $. In particular, we obtain the third inequality in \eqref{Theta_E n=2 5}. The second inequality in \eqref{Theta_E n=2 5} is the same as the first inequality of \eqref{Theta_E n=2 6}.     From $\phi^\sigma(A_0) -\phi^\sigma(A_2) < 1+\min\{\beta, 0\} $   and \eqref{min beta, gamma leq gamma'1}  we get  $\phi(A_0)-\phi^\sigma(Y) < \phi^\sigma(A_0) -\phi^\sigma(A_2) < 1+\alpha' $ and the first   inequality in \eqref{Theta_E n=2 6} is verified also. Thus we showed \eqref{Theta_E n=2 6}. As in the previous case this implies that $ \Theta_{\mc E}\cap  \Theta_{L_1(\mc E)}$ is contractible.

Finally, recall that $\mc E \sim R_0(L_0(\mc E))$, therefore   $\Theta_{\mc E}\cap \Theta_{L_0(\mc E)}$ = $\Theta_{R_0(L_0(\mc E))}\cap \Theta_{L_0(\mc E)}$  and by the already proved first  case we see that  $\Theta_{\mc E}\cap \Theta_{L_0(\mc E)}$ is contractible.  For the  case $\Theta_{\mc E}\cap \Theta_{R_1(\mc E)}$  we have $\Theta_{\mc E}\cap \Theta_{R_1(\mc E)}=\Theta_{L_1(R_1\mc E)}\cap \Theta_{R_1(\mc E)}$ and contractibillity follows from a previous case.  The Proposition is proved.
\epr

Propositions \ref{phi_1>phi_2} and \ref{all exc semistable} ensure semi-stability of certain exceptional objects. The following two  lemmas are similar in that respect and will be used later, when we analyze the intersections of the form  $\Theta_{\mc E_1}\cap \Theta_{\mc E_2}$, when  $\mc E_2$ is obtained from $\mc E_1$  by more than one and different mutations.

\begin{lemma}\label{three comp factors} Let $\mc T=D^b(\mc A)$, where $\mc A$ is a hereditary abelian hom-finite category, and let for any two exceptional objects $E,F \in \mc T_{exc}$ there exists at most one $k\in \ZZ$ satisfying $\hom^k(E,F)\neq 0$.  

 Let $(A_0,A_1,A_2)$ be a full  Ext-exceptional (``Ext-'' means that it satisfies (a) in Definition \ref{sigma triple}) collection in $\mc T$, such that $\hom^1(A_0,A_2)= 0$ and  $A_0, A_1, A_2$ are semistable.    Let $X,Y$ be  exceptional objects in $\mc T$ for which we have a diagram of distinguished triangles, where all arrows are non-zero:  
\be
\begin{diagram}[size=1em] \label{JH filtration}
0 & \rTo      &        &       &   A_2 & \rTo      &            &        & X     &   \rTo      &           &       & Y \\
  & \luDashto &        & \ldTo &         &\luDashto  &          &  \ldTo &       &   \luDashto &           & \ldTo &       \\
  &           &  A_2   &       &         &           &  A_1     &        &       &             & A_0       &
\end{diagram}.
\ee
{\rm (a) } If we have the following system of inequalities: 
$$ \begin{array}{c}   \phi(A_0)-1 < \phi(A_1) < \phi(A_0),  \ \ \ \phi(A_0)-1 < \phi(A_2) < \phi(A_0) \\ 
  \arg_{(\phi(A_0)-1,\phi(A_0))}(Z(A_0)+Z(A_1)) > \phi(A_2)
\end{array}, $$

then $Y \in \sigma^{ss}$ and $\phi(Y)<\phi(A_0)$.

{\rm (b) } If we have 
  $   \phi(A_2)< \phi(A_1)\leq  \phi(A_0)  < \phi(A_2)+1$,  
then $Y \in \sigma^{ss}$ and $\phi(Y)<\phi(A_0)$.
 \end{lemma}
\bpr  We note first some vanishings. From the given diagram it follows that $\hom(Y,A_0)\neq 0$  and $\hom(X,A_1)\neq 0$.  Since $X,Y$ are also exceptional objects, from \cite[Lemma 9.1]{DK1} and the hereditariness of $\mc A$ it follows that $\hom(A_0,Y)=\hom(A_1,X)= 0$. On the other hand  $\hom(A_1,Y)=\hom(A_1,X)$ (follows by applying $\hom(A_1,\_)$ to the last triangle and using $\hom^*(A_1,A_0)=0$). 
Thus, we obtain 
\begin{gather} \label{three comp factors1}\hom(A_0,Y)=\hom(A_1,Y)= 0. \end{gather}
 Since    $(A_0,A_1,A_2)$ is an Ext-exceptional collection, its extension closure is a heart of a bounded t-structure(\cite[Lemma 3.14]{Macri}), furthermore this heart is of finite length and $(A_0,A_1,A_2)$ are the simple objects in it.  Let us denote for simplicity $t=\phi(A_0)$(in case (a)) or   $t-1=\phi(A_2)$(in case (b)).   In both the cases  from the given inequalities and since  $\mc P(t-1,t]$ is also a heart, it follows that     the extension closure of  $(A_0,A_1,A_2)$  is exactly  $\mc P(t-1,t]$.  Now \eqref{JH filtration} can be considered as the Jordan-H\"older filtration of $Y$ in the abelian category $\mc P(t-1,t]$ and the composition factors of $Y$ are $\{A_0,A_1,A_2\}$.

Suppose that $Y\not \in  \sigma^{ss}$.  From Remark \ref{arg remark} (b) there exists $Y'\in \mc P(t-1,t]$ and a non-trivial monic arrow $Y'\rightarrow Y$, s. t. $\arg_{(t-1,t]}(Z(Y'))>\arg_{(t-1,t]}(Z(Y))$. 

 We have  $ Z(Y)=Z(A_0)+Z(A_1)+Z(A_2) $ and one can show that  the given inequalities in either  case  (a) or (b) imply that  \begin{gather}  \label{arg condition1} \arg_{(t-1,t]}(Z(Y)) >  \arg_{(t-1,t]}(Z(A_2)), \quad \arg_{(t-1,t]}(Z(Y)) >  \arg_{(t-1,t]}(Z(A_2)+Z(A_1)). \end{gather}
 Since $Y'$ is a subobject of $Y$, the composition factors of $Y'$ in $\mc P(t-1,t]$ are subset of $\{A_0,A_1,A_2\}$. The cases $Y'\cong A_0$, $Y'\cong A_1$ are excluded by \eqref{three comp factors1}.  The case  $Y'\cong A_2 $ is excluded by the first inequality in  \eqref{arg condition1} and the condition  $\arg_{(t-1,t]}(Z(Y'))>\arg_{(t-1,t]}(Z(Y))$. Since $Y'$ is a proper subobject of $Y$ we reduce to the case when $Y'$ has two composition factors (two different elements of the set $\{A_0,A_1,A_2\}$). Using \eqref{three comp factors1} again we reduce to the following options for a  Jordan H\"older filtration \begin{gather} \bd[1em] 0 &\rTo & A_2 &\rTo & Y'  &\rTo & A_1 &\rTo & 0 \ed  \qquad   \bd[1em] 0 &\rTo & A_2 &\rTo & Y'  &\rTo & A_0 &\rTo & 0 \ed. \end{gather}
In the first case we have $Z(Y')=Z(A_2)+Z(A_1)$ which contradicts the second inequality on \eqref{arg condition1}. In the second case  we have a distinguished triangle$ \bd[1em]  A_2 &\rTo & Y'  &\rTo & A_0 &\rTo & A_2[1] \ed $ in  $\mc T$, and form the given vanishing $\hom^1(A_0,A_2)=0$ it follows $Y'\cong A_0\oplus A_2$, which contradicts \eqref{three comp factors1}. 
So we proved $Y\in \sigma^{ss}$. 
The inequality  $\phi(Y)<\phi(A_0)$(in either case (a) or (b)) follows from the given inequalities and  $ Z(Y)=Z(A_0)+Z(A_1)+Z(A_2) $. 
 The lemma is proved.
\epr
\begin{lemma}\label{two comp factors} Let    $\mc T=D^b(\mc A)$ be as in Lemma \ref{three comp factors}.  

  Let $(A_0,A_1,A_2)$ be a full  Ext-exceptional collection\footnote{ ``Ext-'' means that it satisfies (a) in Definition \ref{sigma triple}} in $\mc T$ such that   $A_0, A_1, A_2$ are  semistable.    Let $Y$ be an exceptional object in $\mc T$ for which we have a  triangle, where all arrows are non-zero:  
\begin{gather}
\begin{diagram}[size=1em] \label{JH filtration1}
 A_2   &   \rTo      &           &       & Y \\
       &   \luDashto &           & \ldTo &   \\
       &             & A_0       &
\end{diagram}.
\end{gather}
 If one of the two systems: 
$ \begin{array}{c}   \phi(A_2) < \phi(A_0) < \phi(A_2)+1  \\ 
   \phi(A_2) < \phi(A_1) < \phi(A_2)+1 
\end{array}$ or   $ \begin{array}{c}   \phi(A_0)-1 < \phi(A_1) < \phi(A_0)  \\ 
   \phi(A_0)-1 < \phi(A_2) < \phi(A_0) 
\end{array}$ holds, 
then we have:  
 $Y \in \sigma^{ss}$,    $\phi(Y)=\arg_{(\phi(A_2), \phi(A_2)+1)}(Z(A_0)+Z(A_2))=\arg_{( \phi(A_0)-1 , \phi(A_0))}(Z(A_0)+Z(A_2))  $ and  $\phi(A_2) < \phi(Y)<\phi(A_0)$.
\end{lemma}
\bpr  Since $Y, A_0$ are exceptional objects and $\hom(Y,A_0)\neq 0$, from \cite[Lemma 9.1]{DK1} and the hereditarines of $\mc A$ it follows that $\hom(A_0,Y)= 0$.   Due to the given inequalities, in both the cases  we can choose $t$ so that $\phi(A_0), \phi(A_1), \phi(A_2) \in (t-1,t]$. By the same arguments as in the previous lemma, one sees that   the extension closure of  $(A_0,A_1,A_2)$  is   $\mc P(t-1,t]$ and that this is an abelian category of finite length with simple objects $A_0,A_1,A_2 $.  Now \eqref{JH filtration1} can be considered as the Jordan-H\"older filtration of $Y$ in the abelian category $\mc P(t-1,t]$ and the composition factors of $Y$ are $\{A_0,A_2\}$. 

 We have  $ Z(Y)=Z(A_0)+Z(A_2) $ and  the given inequalities   (in either  case)
 imply that:   \begin{gather}  \phi(A_2)=\arg_{(t-1,t]}(Z(A_2)) < \arg_{(t-1,t]}(Z(Y))=\arg_{( \phi(A_0)-1 , \phi(A_0))}(Z(A_0)+Z(A_2)) \nonumber  \\[-2mm]  \label{arg condition2}  \\[-2mm]=\arg_{(\phi(A_2), \phi(A_2)+1)}(Z(A_0)+Z(A_2)) < \arg_{(t-1,t]}(Z(A_0))= \phi(A_0).  \nonumber \end{gather}

Suppose that $Y\not \in  \sigma^{ss}$.  From Remark \ref{arg remark} (b) it follows that  there exists $Y'\in \mc P(t-1,t]$ and a non-trivial monic arrow $Y'\rightarrow Y$ in $\mc P(t-1,t]$, s. t. $\arg_{(t-1,t]}(Z(Y'))>\arg_{(t-1,t]}(Z(Y))$. 
 Since  the composition factors of $Y$ are $\{A_0,A_2\}$ and $Y'$ is a non-zero proper sub-object of $Y$, then we have  $Y'\cong A_2$ or $Y'\cong A_0$. The case $Y'\cong A_0$,  is excluded by $\hom(A_0,Y)= 0$.   The case  $Y'\cong A_2 $ is excluded by   \eqref{arg condition2} and the condition  $\arg_{(t-1,t]}(Z(Y'))>\arg_{(t-1,t]}(Z(Y))$. 
 So we proved $Y\in \sigma^{ss}$. 

 By $Y\in \mc P(t,t+1]$ it follows that $\phi(Y)=\arg_{(t-1,t]}(Z(Y))$. Now the lemma  follows from \eqref{arg condition2}. \epr

\section{The exceptional objects  in \texorpdfstring{$D^b(Q )$}{\space}} \label{the exceptional objects}
From now on we fix $\mc T=D^b(Rep_k(Q ))$, where $Q$ is the affine quiver in figure \eqref{Q1}.

In this Section  we organize in a better way the  data about $\{ \Hom(X,Y), {\rm Ext}^1(X,Y)\}_{X\in \mc T_{exc}}$ obtained in \cite[Section 2]{DK1}. In Subsection \ref{central charges} are  given  some observations about the behavior of the vectors  $\{ Z(X)\}_{X\in \mc T_{exc}}\subset \CC$, which  will be helpful when we analyze the intersections of the form $\Theta_{\mc E_1}\cap \Theta_{\mc E_2}$ in the next sections.  

We start by recalling the classification  of the exceptional objects in $Rep_k(Q )$ obtained in \cite{DK1}.
Let us denote    for any $m\geq 1$:
\begin{gather} \nonumber \pi_+^m: k^{m+1} \rightarrow k^{m}, \quad  \pi_-^m: k^{m+1} \rightarrow k^{m}, \quad j_+^m: k^{m} \rightarrow k^{m+1}, \quad  j_-^m: k^{m} \rightarrow k^{m+1}  \\
\nonumber  \pi_+^m(a_1,a_2,\dots, a_m, a_{m+1}) =(a_1,a_2,\dots, a_m) \qquad   \pi_-^m(a_1,a_2,\dots, a_m, a_{m+1})=(a_2,\dots, a_m, a_{m+1}) \\
\nonumber  j_+^m(a_1,a_2,\dots, a_m) =(a_1,a_2,\dots, a_m,0)  \qquad
  j_-^m(a_1,a_2,\dots, a_m)=(0,a_1,\dots, a_m).
\end{gather}

In \cite{DK1} was shown that: 

\begin{prop}\cite[Proposition 2.2]{DK1} \label{exceptional objects in Q1} The exceptional objects up to isomorphism in  $ Rep_{k}(Q ) $ are ($m=0,1,2,\dots$)
\begin{gather} E_1^{m} = \begin{diagram}[1em]
   &       &  k^m &       &    \\
   & \ruTo^{\pi_+^m} &    & \luTo^{Id} &       \\
k^{m+1}  & \rTo^{\pi_-^m}  &    &       &  k^m
\end{diagram} \ \ \ \  E_2^m = \begin{diagram}[1em]
   &       &  k^{m+1} &       &    \\
   & \ruTo^{j_+^m} &    & \luTo^{Id} &       \\
k^{m}  & \rTo^{j_-^m}  &    &       &  k^{m+1}
\end{diagram} \ \ \ \  E_3^m = \begin{diagram}[1em]
   &       &  k^{m+1} &       &    \\
   & \ruTo^{j_+^m} &    & \luTo^{j_-^m} &       \\
k^{m}  & \rTo^{Id}  &    &       &  k^{m}
\end{diagram} \nonumber  \\
E_4^m = \begin{diagram}[1em]
   &       &  k^m &       &    \\
   & \ruTo^{\pi_+^m} &    & \luTo^{\pi_-^m} &       \\
k^{m+1}  & \rTo^{Id}  &    &       &  k^{m+1}
\end{diagram} \ \ \ \ M = \begin{diagram}[1em]
   &       &  0 &       &    \\
   & \ruTo &    & \luTo &       \\
0  & \rTo  &    &       &  k
\end{diagram} \ \ \ \  M'= \begin{diagram}[1em]
   &       &  k &       &    \\
   & \ruTo^{Id} &   &  \luTo  &       \\
k  &  \rTo &    &       &  0
\end{diagram}.\nonumber \end{gather}
\end{prop}

We denote by $K(\mc T)$ the Grothendieck group of $\mc T$.  For $X\in \mc T$ we denote by $[X]\in K(T)$   the corresponding equivalence class in $K(\mc T)$.  From Proposition \ref{exceptional objects in Q1} it follows:
\begin{coro} \label{[X]} Let us denote $\delta =[E_1^0]+[E_3^0]+[M]\in K(\mc T)$. We  have the following equalities in $K(\mc T)$: 
\begin{gather} \label{[delta]} \delta =[E_1^0]+[E_3^0]+[M] = [E_1^0]+[E_2^0] =[E_3^0]+[E_4^0] = [M]+ [M']\\
\label{[E_1^m]} [E_1^m]= m \delta + [E_1^0] = (m+1) \delta - [E_2^0]  \qquad [E_2^m]= m \delta + [E_2^0] = (m+1) \delta - [E_1^0] \\ \label{[E_3^m]}
[E_3^m]= m \delta + [E_3^0] = (m+1) \delta - [E_4^0]  \qquad [E_4^m]= m \delta + [E_4^0] = (m+1) \delta - [E_3^0] \\ 
\label{E M M'} [E_1^m]+[M]=[E_4^m] \quad [E_3^m]+[M]=[E_2^m] \quad  [E_4^m]+[M']=[E_1^{m+1}] \quad [E_2^m]+[M']=[E_3^{m+1}]. 
\end{gather}
\end{coro}

	\subsection{The two orbits of 2-Kronecker pairs in \texorpdfstring{$D^b(Q)$}{\space}} \label{the two orbits}\mbox{}\\
	
	In Propositions 	 \ref{all exc semistable} and \ref{phi_1>phi_2} were discussed  exceptional pairs  $(E,F)$  with $\hom^{\leq 0}(E,F)=0$ and $\hom^{1}(E,F)=l\neq 0$, and their extension closures. We call such a pair \textit{$l$-Kronecker pair}.  Kronecker  pairs were used in \cite{DHKK} for studying   the density of the set of phases of  Bridgeland stability conditions. In \cite[Corollary 3.31]{DHKK}  was shown that for any affine acyclic quiver $A$ (like the quiver $Q$ in figure \eqref{Q1}) only 1- and 2-Kronecker pairs can appear in $D^b(A)$.   In this subsection we give some comments  on the 1- and 2-Kronecker pairs in $D^b(Q)$, which will be useful later when we apply Propositions \ref{all exc semistable}, \ref{phi_1>phi_2} and Lemmas \ref{three comp factors},  \ref{two comp factors}.
	
 From \cite[Remark 2.11]{DK1} we see that the 2-Kronecker pairs in $D^b(Q )$ up to shifts are:	
	\begin{gather} \label{Kron12}\mk{P}_{12}=\{ (E_1^{m+1},E_1^m[-1]),   (E_1^0,E_2^0), (E_2^m,E_2^{m+1}[-1]): m\in \NN \} \\
		\label{Kron34} \mk{P}_{43}=\{ (E_4^{m+1},E_4^{m}[-1]),   (E_4^0,E_3^0), (E_3^m,E_3^{m+1}[-1]): m\in \NN \}.
		\end{gather}
	Recall that the Braid group on two strings $B_2\cong \ZZ$ acts on the set of equivalence classes of exceptional pairs in $\mc T$ (here we take the equivalence $\sim$ explained in \textbf{Some notations} and it is clear when a given equivalence class w.r. $\sim$  will be called a 2-Kronecker pair). Using \cite[Corollary 2.10]{DK1} and the list of triples  in \cite[Corollary 2.9]{DK1} one can show that the set of 2-Kronecker pairs is invariant  under  this action of $B_2$  and this action on the  2-Kronecker pairs has two orbits. They are \eqref{Kron12} and \eqref{Kron34}. 
	
We will describe now the sets ${\mc T}_{exc}\cap \mc A$, up to isomorphism,  where $\mc A$ is the extension closure in $\mc T$ of a 2-Kronecker pair. This will be helpful  later (e. g. when we apply Propositions \ref{all exc semistable} and \ref{phi_1>phi_2}).  
We note first a simple lemma (in which $Rep_k(Q)$ can be any hereditary category):
\begin{lemma} \label{simple lemma} Let $A$, $B \in Rep_k(Q)$, let $\mc C$ be the extension closure of $A, B[-1]$ in $D^b(Rep_k(Q))=\mc T$. Then any $X \in \mc C \cap \mc T_{exc}$ has the form $X'[i]$, where $X'\in Rep_k(Q )_{exc}$ and $i\in \{0,-1\}$. 
\end{lemma}
\bpr Since $Rep_k(Q)$ is hereditary,   any object $X \in D^b(Rep_k(Q))$  decomposes as follows
$X \cong \bigoplus_{i \in \ZZ} H^i(X)[-i]$, where $H^i: \mc T \rightarrow Rep_k(Q)$ are the cohomology functors.   

Since $A,B\in  Rep_k(Q)$, it follows that  $H^i(A) = H^i(B[-1]) = 0$ for each  $i \neq \{0,1\}$. The  functors $H^i: \mc T \rightarrow Rep_k(Q)$ map triangles to short exact sequences (see e.g. \cite{GM}), therefore  $H^i(X)=0$ for any  $X \in \mc C$ and  any  $i \neq \{0,1\}$. By the first paragraph of the proof we see that each  $X \in \mc C$ has the form $X' \oplus X''[-1]$ with  $X'$, $X'' \in Rep_k(Q)$. Finally, if $X \in \mc C \cap \mc T_{exc}$ , then $X$ is indecomposable in $\mc T$, hence  either $X\cong X'$ or   $X\cong X'[-1]$ for some $X'\in Rep_k(Q)$ and obviously $X'$ is also exceptional, i. e.  $X'\in Rep_k(Q)_{exc}$. 
\epr
\begin{lemma} \label{closures of 2-Kronecker}  Let $(U,V)$  be one  of the  2-Kronecker pairs given in \eqref{Kron12} or \eqref{Kron34}. Let $\mc A$ be its extension closure in $\mc T$. Then representatives of the iso-classes of objects in $\mc A \cap \mc T_{exc}$ are:

	\begin{gather} \begin{array}{| c | c | c | c |}
  \hline
   (U,V)     =               &   (E_{1/4}^{m+1},E_{1/4}^{m}[-1] )         &   (E_{1/4}^0,E_{2/3}^0 )  &     (E_{2/3}^{m},E_{2/3}^{m+1}[-1] )         \\ \hline
 \mc A \cap {\mc T}_{exc}= & \left \{\begin{array}{c c} E_{1/4}^n[-1] &  0\leq n\leq m  \\ 
	                                                       E_{1/4}^n     &   n\geq m+1  \\ 
																												  E_{2/3}^n    & n\in \NN     \end{array}\right  \} & \left \{\begin{array}{c c} E_{1/4}^n & n\in \NN \\ E_{2/3}^n & n\in \NN  \end{array} \right \}    & \left \{\begin{array}{c c} E_{2/3}^n &  0\leq n\leq m  \\ 
	                                                       E_{2/3}^n[-1]     &   n\geq m+1  \\ 
																												  E_{1/4}^n[-1]    & n\in \NN     \end{array}\right  \}        \\ \hline
	\end{array} \nonumber
	\end{gather}
 where the  subscript in the table is either  everywhere the first or everywhere the second. 
 \end{lemma}
\bpr We show the case when the subscript is everywhere the first (i. e. the pairs in \eqref{Kron12}), the other case is analogous.  From \cite[Lemma 3.22]{DHKK} we have that $\mc A$ is a bounded t-structure in $\scal{U,V}$ and we have also an  equivalence of abelian categories \begin{gather} \label{closures of 2-Kronecker1} F:\mc A \rightarrow Rep_k(K(2)) \qquad   F(U)=\bd[1em]k & \pile{\rTo \\ \rTo }& 0  \ed, \quad F(V)=\bd[1em]0 & \pile{\rTo \\ \rTo }& k  \ed.    \end{gather}
Using the facts that $\mc A$ is a bounded t-structure and that $F$ is equivalence, one can  show that if $X \in \mc A \cap \mc T_{exc}$, then $F(X) \in Rep_k(K(2))_{exc}$. Furthermore, since $\mc T = D^b(Rep_k(Q ))$ and $Rep_k(Q )$  is a  hereditary abelian category, it is easy to show that: 
\begin{gather} \label{closures of 2-Kronecker2} X \in \mc A \cap {\mc T}_{exc}  \ \   \Leftrightarrow \ \  F(X) \in Rep_k(K(2))_{exc}.\end{gather}
As in the proof of \cite[Proposition 2.2]{DK1} one can classify $ Rep_k(K(2))_{exc} $  and the result is:

	\begin{gather}\label{closures of 2-Kronecker3} \forall X \in  Rep_k(K(2))_{exc}     \qquad \quad  X \cong  \bd k^{n+1} & \pile{\rTo^{\pi_+^n} \\ \rTo_{\pi_-^n}} & k^n \ed    \ \ \mbox{or}   \ \  X \cong    \bd k^{n} & \pile{\rTo^{j_+^n} \\ \rTo_{j_-^n}} & k^{n+1} \ed \ \ \mbox{for some } \   n \in \NN.  \end{gather}

Since $\mc A$ is a bounded t-structure in $\scal{U,V}$, the inclusion functor $\mc A \rightarrow \mc T$ induces an embedding of groups $K(\mc A)\rightarrow K(\mc T)$. Now from \eqref{closures of 2-Kronecker1},  \eqref{closures of 2-Kronecker2}, \eqref{closures of 2-Kronecker3}  it follows that: 
\begin{gather} \label{closures of 2-Kronecker4} \left \{[X] \in K(\mc T): X \in \mc A \cap \mc T_{exc}\right \}=\left \{(n+1) [U]+n [V], \ \  n [U]+(n+1) [V]: n\in \NN\ \right \}.\end{gather}
\ul{If $(U,V)=(E_1^{m+1},E_1^{m}[-1] )$,} then using \eqref{[E_1^m]} we obtain: 
\begin{gather}  (n+1) [U]+n [V]=(n+1)\left [E_1^{m+1}\right ]-n \left [E_1^{m}\right ]=(n+1) \left ((m+1)\delta + \left [E_1^0 \right ] \right )-n \left  (m \delta +\left  [E_1^0 \right ] \right )  \nonumber \\
= (n +m+1) \delta +\left [E_1^0\right ] =\left [E_1^{n+m+1} \right ] \nonumber \\
n [U]+(n+1) [V]=n\left  ((m+1)\delta + \left [E_1^0\right ]\right )-(n+1) \left  (m \delta +\left [E_1^0\right ] \right ) 
= (n -m) \delta -\left  [E_1^0\right ] \nonumber  \\  = \left \{ \begin{array}{c c} \left [E_2^{n-m-1} \right ] \  & \  n \geq m+1  \\ -\left [ E_1^{m-n} \right ] =\left [ E_1^{m-n}[-1] \right ]  \  & \  n \leq m  \end{array}    
 \right. . \nonumber  
\end{gather}
Hence \eqref{closures of 2-Kronecker4} in this case is 
$  \left \{[X] \in K(\mc T): X \in \mc A \cap \mc T_{exc}\right \}=  \left \{\begin{array}{c c}\left [ E_1^n[-1] \right ] &  0\leq n\leq m  \\ 
	                                                      \left [ E_1^n \right ]    &   n\geq m+1  \\ 
																												\left [  E_2^n \right ]   & n\in \NN     \end{array}\right  \} 
$.  
Now the second column in the table follows easily from Lemma \ref{simple lemma} and the fact that     there is at most one, up to isomorphism, exceptional  representation in $Rep_k(Q)$  of a given dimension vector (\cite[p. 13]{WCB2}).  

\ul{If $(U,V)=(E_1^0,E_2^0 )$,} then using \eqref{[E_1^m]} and \eqref{[delta]} we reduce   \eqref{closures of 2-Kronecker4} to $\left \{[X] \in K(\mc T): X \in \mc A \cap \mc T_{exc}\right \}=\{[E_1^n],[E_2^n]: n\in \NN \}$ and  the third column of the table follows.

\ul{If $(U,V)=(E_2^{m},E_2^{m+1}[-1] )$,} then using \eqref{[E_1^m]} we obtain: 
\begin{gather}  (n+1) [U]+n [V]=(n+1) \left (m\delta + \left [E_2^0 \right ] \right )-n \left  ((m+1) \delta +\left  [E_2^0 \right ] \right )  
= (m-n) \delta +\left [E_2^0\right ]  \nonumber \\  
= \left \{ \begin{array}{c c} \left [E_2^{m-n} \right ] \  & \  n \leq m  \\ -\left [ E_1^{n-m-1} \right ] =\left [ E_1^{n-m-1}[-1] \right ]  \  & \  n \geq m+1  \end{array}  \right.  \nonumber \\[2mm]
n [U]+(n+1) [V]=n\left  (m\delta + \left [E_2^0\right ]\right )-(n+1) \left  ((m+1)\delta +\left [E_2^0\right ] \right ) 
\nonumber  \\= -(n +m+1) \delta -\left  [E_2^0 \right ]   =  \left [ E_2^{n+m+1}[-1] \right ]  
\nonumber  
\end{gather}
and now \eqref{closures of 2-Kronecker4} and similar arguments as in the first case give the fourth column of the table.  
 
  The case when the subscript is everywhere the second (i. e. the pairs in \eqref{Kron34})
 is obtained by substituting $E_1$ with $E_4$,  $E_2$ with $E_3$, and using \eqref{[E_3^m]} instead of \eqref{[E_1^m]}. 
\epr

Some 1-Kronecker pairs in $D^b(Q)$ are  (see \cite[table (4)]{DK1}):
\begin{gather} \label{1-kronecker} (M', E_1^m[-1]), (M', E_2^m), (M, E_3^m), (M, E_4^m[-1]).\end{gather}
 In the following lemma are listed several short exact sequences in $Rep_k(Q)$. On one hand, these sequences determine the set $\mc A\cap \mc T_{exc}$, where $\mc A$ is the extension closure of some of the 1-Kronecker pairs in \eqref{1-kronecker}, so they will be helpful  when we apply Propositions \ref{all exc semistable} and \ref{phi_1>phi_2}. On the other hand, they (and their combinations) will play the role of the triangles \eqref{JH filtration} and \eqref{JH filtration1} when we apply Lemmas \ref{three comp factors} and \ref{two comp factors}.
\begin{lemma} There exist arrows in $Rep_k(Q )$ as shown below, so  that  the resulting  sequences are  exact($m\in \NN$):
\begin{gather}
\label{ses2} \bd 0 & \rTo & E_3^{m} & \rTo & E_2^m & \rTo & M & \rTo & 0 \ed \\
\label{ses4} \bd 0 & \rTo & M & \rTo & E_4^m & \rTo & E_1^{m} & \rTo & 0 \ed\\
\label{ses7} \bd 0 & \rTo & M' & \rTo &E_1^{m+1} & \rTo & E_4^m  & \rTo & 0 \ed \\
\label{ses8} \bd 0 & \rTo & E_2^m & \rTo & E_3^{m+1} & \rTo & M' & \rTo & 0 \ed\\
\label{ses10} \bd 0 & \rTo & E_3^0 & \rTo & M' & \rTo & E_1^{0} & \rTo & 0. \ed
\end{gather}
\end{lemma}
\bpr The proof is an exercise using  Proposition \ref{exceptional objects in Q1} .\epr 
	\subsection{Recollection of some results of \cite{DK1}  with new notations} \mbox{}\\

	It is useful to introduce some notations (see Proposition \ref{exceptional objects in Q1} for the notations $E_i^j$, $M$, $M'$):
\begin{gather} \label{a^m b^m} a^m=\left \{  \begin{array}{c c} E_1^{-m}  & m\leq 0 \\  E_2^{m-1}[1]  & m\geq 1 \end{array} \right. ; \qquad \qquad \qquad b^m= \left \{ \begin{array}{c c} E_4^{-m}  & m\leq 0 \\  E_3^{m-1}[1]  & m\geq 1  \end{array} \right . . \end{gather}

\begin{remark} \label{exceptional objects}  The objects in $\mc T_{exc}$ up to isomorphism  are   $\{a^j[k], b^j[k], M[k],M'[k]: j\in \ZZ, k\in \ZZ\}$. 
\end{remark}

  Using   \cite[table (4) in Proposition 2.4]{DK1}), one verifies that:
\begin{coro}(of \cite[Proposition 2.4]{DK1})\label{nonvanishings} For each $m\in \ZZ$ we have: 
\begin{gather} \label{nonvanishing1} \hom(M',a^m)\neq 0; \quad \hom(M,b^m)\neq 0; \quad  \hom^*(a^m,M')= 0;  \\ 
   \label{nonvanishing2} \hom^1(a^m,M)\neq 0; \quad \hom^1(b^m,M')\neq 0; \quad  \hom^*(b^m, M)=0 \\ 
 \label{nonvanishing3}  \hom^1(b^{m+1},a^n) \neq 0 \ \mbox{for} \ m> n;   \quad  \hom(b^{m},a^{n}) \neq 0   \ \mbox{for} \ m\leq n; \quad  \hom^*(b^{m+1},a^m) = 0 \\
  \label{nonvanishing4}   \hom^1(a^{m},b^{n}) \neq 0 \ \mbox{for} \ m>n;   \quad   \hom(a^m,b^{n+1}) \neq 0    \ \mbox{for} \ m\leq n;  \quad   \hom^*(a^m,b^{m}) = 0; \\
    \label{nonvanishing5}  \hom(a^{m},a^n) \neq 0 \ \mbox{for} \ m\leq  n;   \quad  \hom^1(a^{m},a^{n}) \neq 0   \ \mbox{for} \ m> n+1; \quad \hom^*(a^{m},a^{m-1})=0\\ 
		 \label{nonvanishing6}  \hom(b^{m},b^n) \neq 0 \ \mbox{for} \ m\leq  n;   \quad  \hom^1(b^{m},b^{n}) \neq 0   \ \mbox{for} \ m> n+1; \quad \hom^*(b^{m},b^{m-1})=0 \\ 
	\label{nonvanishing7}	\hom^1(M,M') \neq 0 \qquad \hom^1(M',M)\neq 0 .\end{gather}
		\end{coro}
 It is useful to keep in mind  the following remarks:
\begin{remark}  \label{inc dec seq} Recall that $\phi_-(A)> \phi_+(B)$ implies $\hom(A,B) = 0$  and in particular $\hom(A,B) \neq  0$ implies  $\phi_-(A)\leq \phi_+(B)$ (for each stability condition).  

 Let $\{x^i\}_{i\in \ZZ}$ be either $\{a^i\}_{i\in \ZZ}$ or  $\{b^i\}_{i\in \ZZ}$.  From  \eqref{nonvanishing5}  and \eqref{nonvanishing6}  we see that:

{\rm (a)} For $m\leq n$  we have  $\hom(x^{m}, x^n) \neq 0  $.   In particular,  if $x^{m}, x^n \in \sigma^{ss}$ and $m\leq n$ then  $\phi(x^m)\leq \phi(x^n)$.

 {\rm (b)} For $m+1<n$ we have $\hom^1(x^n, x^{m}) \neq 0$.  In particular,  if $x^{m}, x^n \in \sigma^{ss}$ and $m+1<n$ then $\phi(x^n)\leq \phi(x^{m})+1$.
\end{remark}

\begin{remark} \label{ext closure ab} Let $\{x^i\}_{i\in \ZZ}$ be either $\{a^i\}_{i\in \ZZ}$ or  $\{b^i\}_{i\in \ZZ}$.   Lemma \ref{closures of 2-Kronecker} in terms of the notations \eqref{a^m b^m} is equivalent to saying that for any three integers  $i\leq p$, $ p+1\leq j$ we have that $x^i$ and $x^j[-1]$ are in the extension closure of $\{x^p,x^{p+1}[-1]\}$.
\end{remark}

Keeping these remarks in mind one easily proves:
\begin{lemma}\label{from phases of big distance}   Let $\{x^i\}_{i\in \ZZ}$ be either $\{a^i\}_{i\in \ZZ}$ or  $\{b^i\}_{i\in \ZZ}$.  
 If there exists $m\in \ZZ$  such that     $\{x^m, x^{m+1}\}\subset \sigma^{ss}$ and $\phi(x^m)+1<\phi(x^{m+1})$, then for $i\not \in \{m,m+1\}$ we have  $x^i\not \in \sigma^{ss}$. 
\end{lemma}
\bpr 
Suppose $ x^{i} \in \sigma^{ss}$ with $i< m$,  then by Remark (a) we have \ref{inc dec seq} $\phi(x^{i})+1\leq \phi(x^{m})+1<\phi(x^{m+1})$, hence $\hom^1(x^{m+1},x^{i})=0$, which contradicts the second part of Remark \ref{inc dec seq} (b).  If $ x^{i} \in \sigma^{ss}$ with $i> m+1$, then by Remark \ref{inc dec seq} (a) we obtain  $\hom^1(x^i,x^m)=0$, which  again contradicts  Remark \ref{inc dec seq} (b).
\epr

We will  use often the following result obtained in \cite{DK1}:
\begin{coro} \label{one nonvanishing degree} \cite[Corollary 2.6 (b)]{DK1}  For any two exceptional objects $X, Y \in D^b(Q )$ at most one  element of the family $\{ \hom^p(X,Y) \}_{p\in \ZZ}$ is nonzero.
\end{coro}
Due to Corollary \ref{one nonvanishing degree} we can apply Lemmas \ref{three comp factors}, \ref{two comp factors} to $D^b(Q)$. Furthermore, we have:
	\begin{coro}(\cite{DK1}) \label{exceptional colleections}  The full  exceptional collections in $D^b(Q )$ up to isomorphism and schifts are in the set of triples $\mk T$ given below. Propositions \ref{lemma for f_E(Theta_E)}, \ref{phi_1>phi_2}, \ref{mutations} can be applied to any of these triples.  
\ben {\mk T} = \left \{  \begin{array}{c  c  c} (M',a^{m},a^{m+1})  & (a^m, b^{m+1},a^{m+1}) & (a^m,a^{m+1},M) \\
(M,b^{m},b^{m+1}) &(b^m,a^m,b^{m+1})  & (b^m,b^{m+1},M')\\
 (b^{m},M',a^{m}) &(a^m , M, b^{m+1}) &  .\end{array}: m\in \NN \right \}.
\een
\end{coro}
\bpr The list $\mk T$ follows straightforwardly from \cite[Corollary 2.9]{DK1}. By Corollary \ref{nonvanishings}  $\hom^*(X,Y)\neq 0$  for any  exceptional pair $(X,Y)$, therefore  Propositions \ref{lemma for f_E(Theta_E)}, \ref{phi_1>phi_2} can be applied to any of the triples. Proposition \ref{mutations} can be applied due to Corollary \ref{one nonvanishing degree}. \epr

\begin{remark}	 It is known  \cite{WCB1} that the Braid group on three strings $B_3$ acts transitively on the exceptional triples  of $Rep_k(Q)$. This action is not free (see \cite[Remark 2.12]{DK1}).
	\end{remark} 
\begin{remark} \label{T_1234} With the notations \eqref{a^m b^m} the two orbits of 2-Kronecker pairs (see \eqref{Kron12} and \eqref{Kron34}) are $\{(a^m,a^{m+1}[-1])\}_{m\in \ZZ}$ and $\{(b^m,b^{m+1}[-1])\}_{m\in \ZZ}$. Each of these pairs can be extended to three non-equivalent triples, so we obtain two sets of triples. Having the list $\mk{T}$ above, we see that  these  two sets of triples are: 
\begin{gather} \label{mk{T}_12} \mk{T}_{a}=\{ (M',a^{m},a^{m+1}) ,  (a^m, b^{m+1},a^{m+1}),  (a^m,a^{m+1},M): m\in \ZZ\} \\
\label{mk{T}_43} \mk{T}_{b}=\{(M,b^{m},b^{m+1}), (b^m,a^m,b^{m+1}),  (b^m,b^{m+1},M'): m \in \ZZ \}.
\end{gather}
Furthermore we have:
\begin{gather} \label{one union for mkT} \mk{T} =\mk{T}_{a} \cup \{ (b^{m},M',a^{m}), (a^m , M, b^{m+1}):m\in \ZZ \} \cup \mk{T}_{b} \qquad  \mk{T}_{a}\cap \mk{T}_{b}=\emptyset. \end{gather}
\end{remark}

\begin{remark}  \label{ses ab}
The short exact sequences \eqref{ses7}, \eqref{ses8}, \eqref{ses10} in terms of the notations \eqref{a^m b^m} become a sequence of distinguished triangles (for each $p$): 
\be \label{short filtration 1} 
\begin{diagram}[size=1em] 
b^{p+1}[-1] & \rTo      &        &       &   M'   \\
            & \luDashto &        & \ldTo &        \\
            &           &  a^p   &       &        
\end{diagram}
\ee 
The short exact sequences \eqref{ses2} and \eqref{ses4}  become the following   distinguished triangles  ($q \in \ZZ$): 
\be \label{short filtration 2} 
\begin{diagram}[size=1em] 
a^{q}[-1] & \rTo      &        &       &   M   \\
            & \luDashto &        & \ldTo &        \\
            &           &  b^q   &       &        
\end{diagram}.
\ee 

\end{remark}

\subsection{Comments on the vectors \texorpdfstring{$\{Z(X):X \in \mc T_{exc}\}$}{\space}.} \label{central charges} \mbox{} \\
Corollary \ref{[X]} shows  that for each $\sigma=(\mc P, Z) \in \st(\mc T)$ we have $ Z(\delta) = Z(E_1^0)+ Z(E_2^0)= Z(E_4^0)+Z(E_3^0)$ and  $ Z(E_k^m)= m Z(\delta) + Z(E_k^0) $   for each $m\in \NN$  and each $ k=1,2,3,4$.  Due to these equalities,  with the notations \eqref{a^m b^m}   we can write (recall that $Z(X[j]) = (-1)^j Z(X)$ for any $j\in \ZZ$, $X \in \mc T$):
\begin{gather} 
 \label{Z(delta)}\forall j \in \ZZ \qquad  Z(a^{j+1})=Z(a^{j})-Z(\delta) \ \  \mbox{and} \ \ Z(b^{j+1})=Z(b^{j})-Z(\delta).\end{gather} 
Therefore  for any two integers $m,n$ we have:  
\begin{gather} \label{Z(E_k^m)}  Z(a^m)=Z(a^n)-(m-n) Z(\delta) \qquad  \mbox{and} \qquad  Z(b^m)=Z(b^n)-(m-n) Z(\delta) . \end{gather}
Next we discuss collinear vectors among  $\{Z(a^j)\}_{j\in \ZZ}$ and  $\{Z(b^j)\}_{j\in \ZZ}$.  We fix first the meaning of ``collinear'':
\begin{df} \label{collinear} We say that  a family  $\{A_i\}_{i \in I}$  of complex numbers  is collinear if $\{A_i\}_{i \in I} \subset \RR c$ for some $c \in \CC \setminus \{0\}$. In particular,  $0 \in \CC$ is collinear to any $a\in \CC$.  \end{df}
With this definition we have: 
 \begin{lemma}\label{vectors Z(E_k^m)} Let   $\sigma=(\mc P, Z) \in \st(\mc T)$.  Let $\{x^i\}_{i\in \ZZ}$ be either $\{a^i\}_{i\in \ZZ}$ or  $\{b^i\}_{i\in \ZZ}$.   Recall that $\delta$ is defined in \eqref{[delta]} and consider a  sequence in $\CC$(infinite in both directions) of the form:    \begin{gather} \label{the sequence of vectors} \dots,  Z(x^{-i}),\dots,  Z(x^{-2}),Z(x^{-1}),  Z(x^0),Z(\delta),  Z(x^1),  Z(x^2), Z(x^3), \dots, Z(x^j), \dots. \end{gather}
Then the following conditions are equivalent:  (a)  Two of the vectors  in   this sequence  are collinear; (b) The entire sequence   is collinear.
\end{lemma}  \bpr  
Recall that formula  \eqref{Z(E_k^m)} holds for any $m,n \in \ZZ$. 

If $Z(x^i)$ and $Z(\delta)$ are collinear for some $i\in \ZZ$, then $ Z(\delta)$ and $Z(x^0) $ are collinear  by $Z(x^0) =Z(x^i)+i Z(\delta) $  and (b) follows from the equalities   $Z(x^j)=Z(x^0) -j  Z(\delta)$, $j\in \ZZ$.  

If $Z(x^i)$ and $Z(x^j)$ are collinear for some $i\neq j$, then by  the equality$Z(\delta)=\frac{1}{j-i} (Z(x^i)-Z(x^j))$ we see that   $ Z(\delta)$ and $Z(x^i) $ are collinear and (b) follows from the considered case above
\epr

\begin{coro} \label{from equal phases to}   Let $\{x^i\}_{i\in \ZZ}$ be either $\{a^i\}_{i\in \ZZ}$ or  $\{b^i\}_{i\in \ZZ}$.   Let two of the vectors in the  sequence \eqref{the sequence of vectors} be non-collinear.  
Then:

{\rm (a) }  All the vectors in this sequence  are non-zero and no two of them are collinear. 
 
 {\rm (b) } If  for two  integers  $n\neq m$  holds  $\{x^n, x^m\}\subset \sigma^{ss}$,  then we have   $\phi(x^n)\not \in \phi(x^n)+\ZZ$.

{\rm  (c) }  The   numbers  $\{Z(x^j)\}_{j\in \ZZ} $ are contained in a common connected component of  $\CC \setminus \RR Z(\delta)$.

 {\rm  (d) }   If  for two  integers  $n< m$   we have $\{x^n, x^m\}\subset \sigma^{ss}$ and   $\phi(x^m)<\phi(x^n)+1$, then:\footnote{See \eqref{complement of a line} for the notations $Z(\delta)^c_\pm$.}   $$\{Z(x^j)\}_{j\in \NN} \subset  Z(\delta)^c_+.$$

 \end{coro}
\bpr (a) and (b)  follow from Lemma \ref{vectors Z(E_k^m)}, and   the following  axiom in \cite{Bridg1}:
\be
\label{phase formula}  X\in \sigma^{ss} \qquad
\Rightarrow \qquad  Z(X)=r(X)\ \exp(\ri \pi \phi(X)), \ r(X) >0.
\ee

 Since $Z(x^0)$ and  $Z(\delta)$ are non-colinear, it follows that either $Z(x^0)\in Z(\delta)^c_+$ or $Z(x^0)\in Z(\delta)^c_-$. From formula  \eqref{Z(E_k^m)}  we have   $Z(x^j)=Z(x^0) -j  Z(\delta)$ for any $j\in \ZZ$ therefore either  $\{Z(x^j)\}_{j\in\ZZ}\subset Z(\delta)^c_+$ or $\{Z(x^j)\}_{j\in\ZZ}\subset Z(\delta)^c_-$. Therefore we obtain (c). Now to show (d), it is enough to show that  $Z(x^m)\in Z(\delta)^c_+$.
 From (b)  and  Remark \ref{inc dec seq} (a) we get the inequalities  $\phi(x^n)<\phi(x^m)<\phi(x^n)+1$. By drawing a picture and taking into account formula \eqref{phase formula} and the equality    $Z(x^n)=Z(x^m)+(m-n) Z(\delta)$, one  sees that $\phi(x^n)<\phi(x^m)<\phi(x^n)+1$ is  impossible if  $Z(x^m)\in Z(\delta)^c_-$. 
\epr 
\begin{coro} \label{noncolinear ab}    Let $\{x^i\}_{i\in \ZZ}$ be  either the sequence $\{a^i\}_{i\in \ZZ}$ or   the sequence $\{b^i\}_{i\in \ZZ}$.

If $Z(\delta) \neq 0$ and  $Z(x^q) \in Z(\delta)_+^c$ for some $q\in \ZZ$, then   $\{Z(x^i)\}_{i\in \ZZ}\subset Z(\delta)_+^c$ and for any $t\in \RR$ with $Z(\delta)=\abs{Z(\delta)}\exp(\ri \pi t)$ we have: \begin{gather} \label{noncolinear ab1}  \forall p \in \ZZ \ \ \  \arg_{(t,t+1)}(Z(x^p))<\arg_{(t,t+1)}(Z(x^{p+1})) \\
 \label{noncolinear ab2}    \lim_{p\rightarrow -\infty } \arg_{(t,t+1)}(Z(x^p))=t; \quad  \lim_{p\rightarrow+ \infty } \arg_{(t,t+1)}(Z(x^p))=t+1.  \end{gather}
 \end{coro}
\bpr 

 Since $Z(\delta)$ and $Z(x^q)$ are not collinear by  Corollary \ref{from equal phases to} (c) and  $Z(x^q) \in Z(\delta)_+^c$ it follows that   $\{Z(x^i)\}_{i\in \ZZ}\subset Z(\delta)_+^c$.     The inequalities \eqref{noncolinear ab1} follow from  $Z(x^{p+1}), Z(x^p) \in  Z(\delta)_+^c$ and  $ Z(x^{p+1})=Z(x^p)-Z(\delta)$  (see  \eqref{Z(E_k^m)} ). The formulas in   \eqref{noncolinear ab2} follow also  from \eqref{Z(E_k^m)} and $\{Z(x^i)\}_{i\in \ZZ}\subset Z(\delta)_+^c$. 
\epr

\begin{coro} \label{noncolinear ab full1}  Let $Z(M)$ and $Z(M')$ be non-zero and  have the same direction.\footnote{We mean that $Z(M)= y Z(M')$ for some $y\in \RR_{>0}$. In particular, by  $Z(\delta)=Z(M)+Z(M')$  (recall \eqref{[delta]}) we see that $Z(\delta)$ is non-zero.}  Let  $Z(a^q), Z(b^p) \in Z(\delta)_+^c$ for some $p,q\in \ZZ$. 

Then  $\{Z(a^j), Z(b^j)\}_{j\in \ZZ}\subset Z(\delta)^c_+$ and  for any $t\in \RR$ with   $Z(\delta)=\abs{Z(\delta)}\exp(\ri \pi t)$ the formulas  \eqref{noncolinear ab2}, \eqref{noncolinear ab1} hold  for both the sequences  $\{Z(a^j)\}_{j\in \ZZ}$ and $\{Z(b^j)\}_{j\in \ZZ}$.

 Furthermore,  for any three integers $i,j,m$ we have:
\begin{gather} \label{j<mleqi} j<m\leq i \ \ \Rightarrow \ \  \arg_{(t,t+1)}(Z(a^j))< \arg_{(t,t+1)}(Z(b^m))<\arg_{(t,t+1)}(Z(a^i)).\end{gather}
\end{coro}
\bpr    Corollary  \ref{noncolinear ab}   shows the first part of the conclusion.  
To show \eqref{j<mleqi} we note first that   the equalities \eqref{E M M'} with the notations \eqref{a^m b^m} become the following (for any $m\in \ZZ$): 
\begin{gather} \label{a b M M'} Z(b^m)-Z(M)=Z(a^m) \qquad Z(b^m)+Z(M')=Z(a^{m-1}).  \end{gather}
Since  $Z(a^{m-1}),  Z(a^{m}), Z(b^{m}) \in Z(\delta)^c_+$ for any $m\in \ZZ$ and $Z(M)$, $Z(M')$ have the same direction as $Z(\delta)$ (recall \eqref{[delta]})    the equalities \eqref{a b M M'} imply that $\arg_{(t,t+1)}(Z(a^{m-1}))< \arg_{(t,t+1)}(Z(b^m))<\arg_{(t,t+1)}(Z(a^m))$ for any $m\in \ZZ$.  Now \eqref{j<mleqi} follows from \eqref{noncolinear ab1} (applied to the case $\{x^i\}_{i\in \ZZ }=\{a^i\}_{i\in \ZZ }$).   
\epr

\section{The union \texorpdfstring{$\st(D^b(Q )) = \mk{T}_{a}^{st} \cup (\_,M,\_) \cup (\_,M',\_) \cup \mk{T}_{b}^{st}$}{\space}} \label{the union}
In this Section we  distinguish some   building blocks of   $\st(D^b(Q ))$  and   organize them  in a manner  consistent with the order in which we will glue these blocks in the next sections.

Theorem \ref{main theorem for Q in intro} says   that for each $\sigma \in \st(D^b(Q ))$ there exists a $ \sigma$-triple.  This means that (see \cite[Corollary 3.20]{DK1})  for each $\sigma$ there exists an Ext-exceptional triple $\mc E$ with  $\sigma \in \Theta_{\mc E}'$. From Corollary \ref{exceptional colleections}  we see that  $\mc E$ is a shift of  some of the triples in $\mk{T}$. Recalling the notation    \eqref{theta_{mc E}} we get
$ \st(D^b(Q )) = \bigcup_{{\mc E} \in \mk{T}} \Theta_{\mc E}.$ Our basic  building blocks are $\{ \Theta_{\mc E}\}_{\mc E \in \mk{T}}$ and by Proposition \ref{lemma for f_E(Theta_E)} they are contractible.   

\textit{For a given triple $(A,B,C)\in \mk{T}$  we will denote the open subset $\Theta_{(A,B,C)} \subset \st(D^b(Q ))$ by   $(A,B,C)$, when (we believe that) no confusion may arise.} With this convention we can write 
\begin{gather} \label{st} \st(D^b(Q )) = \bigcup_{(A,B,C) \in \mk{T}} (A,B,C).
\end{gather}

\textit{For a given $X\in \{M,M'\}$ we denote by $ (X,\_,\_)$ the following  open subset of $\st(D^b(Q )) $:}
\begin{gather} \label{(A,_,_)} \st(D^b(Q )) \supset (X,\_,\_) = \bigcup_{\{(B_0,B_1,B_2) \in \mk{T} : B_0 = X \}} (X,B_1,B_2).
\end{gather}
\textit{Similarly we define $(\_,X,\_)$ and  $(\_,\_,X)$.} Looking at the list $\mk{T}$ and denoting (see \eqref{mk{T}_12}, \eqref{mk{T}_43}):
\begin{gather} \label{mk{T}} \mk{T}_{a}^{st}=\bigcup_{(A,B,C)\in \mk{T}_{a}} (A,B,C) \subset \st(D^b(Q )); \qquad \mk{T}_{b}^{st}=\bigcup_{(A,B,C)\in \mk{T}_{b}} (A,B,C) \subset \st(D^b(Q ))\end{gather}   we can regroup the union \eqref{st} using \eqref{one union for mkT} as follows:
\begin{gather} \label{st with T} \st(D^b(Q )) = \mk{T}_{a}^{st} \cup (\_,M,\_) \cup (\_,M',\_) \cup \mk{T}_{b}^{st}. \end{gather}

\begin{remark} \label{remark about shifts} From the very definition \eqref{theta_{mc E}} of $\Theta_{\mc E}$ it is clear that $\Theta_{\mc E[\textbf{p}]} = \Theta_{\mc E}$ for any triple $\mc E=(A,B,C)\in \mk{T}$ and any $\textbf{p}\in \ZZ^3$. Using the notations explained  here,   we have $(A,B,C)=(A[p_0],B[p_1],C[p_2])\subset \st(D^b(Q ))$ for any $p_0,p_1,p_2 \in \ZZ$.
\end{remark}

\section{Some contractible subsets of  \texorpdfstring{$\mk{T}_{a}^{st} \ \mbox{and}\   \mk{T}_{b}^{st}$}{\space}.  Proof that  \texorpdfstring{$\mk{T}_{a}^{st}\cap \mk{T}_{b}^{st}=\emptyset$}{\space}} \label{mk T_12 cap mk T_43} 

In this  Section is shown that $(X,\_,\_)$ and  $(\_,\_,X)$  are contractible subsets of $\st(\mc T)$ for any $X\in \{M,M'\}$ and that  {$\mk{T}_{a}^{st}\cap \mk{T}_{b}^{st}=\emptyset$. 

We will refer  often
to some of the formulas in Corollary \ref{nonvanishings}.
Whenever we discuss $\hom(A,B)$ or  $\hom^1(A,B)$ with $A,B$ varying in  the symbols $M,M', a^m, b^m$, $m\in \ZZ$, we refer to
Corollary \ref{nonvanishings}.

Putting  \eqref{mk{T}_12}, \eqref{mk{T}_43} in \eqref{mk{T}}  we obtain: \begin{gather} \label{T_12new} \mk{T}_{a}^{st} = (M',\_,\_) \cup (\_,\_,M) \cup \bigcup_{p\in \ZZ}(a^p,b^{p+1},a^{p+1})\\ \label{T_43new} \mk{T}_{b}^{st} =  (M,\_,\_) \cup (\_,\_,M') \cup\bigcup_{q\in \ZZ}(b^q,a^q,b^{q+1}) \\  \label{(M',_,_)} (M',\_,\_) =  \bigcup_{m \in \ZZ} (M',a^m,a^{m+1}) ;  \qquad  (M,\_,\_) =  \bigcup_{m \in \ZZ} (M, b^m,b^{m+1})\\ \label{(_,_,M)} (\_,\_,M) =  \bigcup_{m \in \ZZ} (a^m,a^{m+1},M)   ; \qquad  (\_,\_,M') =  \bigcup_{m \in \ZZ} (b^{m},b^{m+1},M').
\end{gather} 
 
 We apply Proposition \ref{lemma for f_E(Theta_E)} to the triples $(a^p,b^{p+1},a^{p+1})$ and $(b^q,a^q,b^{q+1})$.  Using  Corollary \ref{one nonvanishing degree} and  the formulas in Corollary \ref{nonvanishings} we see that in both the cases the coefficients $\alpha, \beta, \gamma$ defined in \eqref{alpha,beta,gamma} are $\alpha=\beta=\gamma=-1$.  Thus,  we obtain  the following formulas for the sets $(a^p,b^{p+1},a^{p+1})\subset \st(D^b(\mc T))$   and  $(b^q,a^q,b^{q+1})\subset \st(D^b(\mc T))$ in the first and the second column, respectively:
\begin{gather} \label{no M} \begin{array}{| c   | c  |}   \hline  
 (a^p,b^{p+1},a^{p+1})   & (b^q,a^q,b^{q+1})  \\
  \hline
 \left \{a^p,b^{p+1},a^{p+1} \in \sigma^{ss} :                      \begin{array}{c} \phi\left ( a^{p}\right) < \phi\left (b^{p+1}\right) \\ \phi\left ( a^{p}\right)+1 < \phi\left (a^{p+1}\right) \\ \phi\left ( b^{p+1}\right) < \phi\left (a^{p+1}\right)\end{array}  \right \}                            &       \left \{ b^q,a^q,b^{q+1} \in \sigma^{ss} :                      \begin{array}{c} \phi\left ( b^{q}\right) < \phi\left (a^{q}\right) \\ \phi\left ( b^{q}\right)+1 < \phi\left (b^{q+1}\right) \\ \phi\left ( a^{q}\right) < \phi\left (b^{q+1}\right)\end{array}  \right \}   \\  \hline
	\end{array}
\end{gather}

Similarly, applying Proposition \ref{lemma for f_E(Theta_E)} to the triples in   the unions  \eqref{(_,_,M)}, \eqref{(M',_,_)} (with the help of Corollary \ref{nonvanishings} and Corollary \ref{one nonvanishing degree})  we see that  $ (M',\_,\_) \cup  (\_,\_,M)   $ and $(M,\_,\_)\cup (\_,\_,M') $ are the unions of the sets in the first and the second column of the following table, respectively (where  $m,n,i,j$ vary throughout $\ZZ$):
\begin{gather} \label{left right M} \begin{array}{| c   | c  |}   \hline  
 (M',\_,\_) \cup  (\_,\_,M)      &  (M,\_,\_)\cup (\_,\_,M')   \\
  \hline
  \left \{  M',a^{j},a^{j+1} \in \sigma^{ss} :   \begin{array}{c} \phi\left (M'\right) < \phi\left (a^{j}\right) \\ \phi\left ( M'\right)+1 < \phi\left (a^{j+1}\right) \\ \phi\left ( a^{j}\right) < \phi\left (a^{j+1} \right)\end{array}\right \}                                                &    \left \{  M,b^{n},b^{n+1} \in \sigma^{ss} :                     \begin{array}{c} \phi\left (M\right) < \phi\left (b^{n}\right) \\ \phi\left ( M\right)+1 < \phi\left (b^{n+1}\right) \\ \phi\left ( b^{n}\right) < \phi\left (b^{n+1} \right)\end{array} \right \} \\  \hline
  \left \{ a^{m},a^{m+1}, M \in \sigma^{ss} :                      \begin{array}{c} \phi\left ( a^{m}\right) < \phi\left (a^{m+1}\right) \\ \phi\left ( a^{m}\right) < \phi\left (M\right) \\ \phi\left ( a^{m+1}\right) < \phi\left (M\right)+1\end{array}  \right \}                            &    \left \{ b^{i},b^{i+1}, M' \in \sigma^{ss}:                       \begin{array}{c} \phi\left ( b^{i}\right) < \phi\left (b^{i+1}\right) \\ \phi\left ( b^{i}\right) < \phi\left (M'\right) \\ \phi\left ( b^{i+1}\right) < \phi\left (M'\right)+1\end{array}    \right \}   .    \\ \hline
	\end{array}
\end{gather}
  For the triples on the first row of table \eqref{left right M} we have  $\alpha=\beta=\gamma-1$ and for the triples on the second row we have  $\alpha=-1$, $\beta=\gamma=0$ (one shows this using  Corollaries \ref{nonvanishings} and  \ref{one nonvanishing degree}).

\subsection{Proof that \texorpdfstring{$\mk{T}_{a}^{st}\cap \mk{T}_{b}^{st}=\emptyset$}{\space}} \label{empty intersection} \mbox{}\\
\textit{Recall the axioms of Bridgeland \cite{Bridg1}, that  $\phi(A[1])=\phi(A)+1$ for any $A\in \sigma^{ss}$, and that   $A, B \in \sigma^{ss}$ and $\phi(A)>\phi(B)$ imply $\hom(A,B)=0$.}  We will use these axiom  often implicitely. 
We start with:
\begin{lemma} \label{mk T_12 cap mk T_43 1} $\left ( (M',\_,\_) \cup  (\_,\_,M)   \right ) \cap \left ( (M,\_,\_)\cup (\_,\_,M')  \right )=\emptyset $. \end{lemma}
\bpr

Suppose  $\sigma \in  (a^{m},a^{m+1}, M) \cap  ( M,b^{n},b^{n+1})$, then  by  the  table \eqref{left right M}  we obtain $\hom^1(b^{n+1},a^m)=0$ and $\hom(b^{n+1},a^{m+1})=0$, which contradicts  \eqref{nonvanishing3}. 

Suppose  $\sigma \in  (a^{m},a^{m+1}, M) \cap  (b^{i},b^{i+1}, M')$, then  by  $\hom(M',a^{m})\neq 0$ (see \eqref{nonvanishing1}) and   table \eqref{left right M}   we obtain $ \phi\left ( b^{i}\right) <\phi(M')\leq\phi\left ( a^{m}\right)  <\phi(M)$, which contradicts $\hom(M,b^i)\neq 0$ (see \eqref{nonvanishing1}). 

Suppose  $\sigma \in  (M',a^{j},a^{j+1} ) \cap  ( M,b^{n},b^{n+1})$, then  by    $\hom^1(a^{j+1},M)\neq 0$ (see \eqref{nonvanishing2}) and   table \eqref{left right M}  we obtain $ \phi(M')+1  <\phi\left ( a^{j+1}\right)\leq \phi(M)+1<\phi\left ( b^{n+1}\right) $, which contradicts $\hom^1(b^{n+1},M')\neq 0$. 

Suppose  $\sigma \in  (M',a^{j},a^{j+1} ) \cap   (b^{i},b^{i+1}, M')$, then  by the table  we have $\hom^1(a^{j+1},b^i)=0$, $\hom(a^{j+1},b^{i+1})=0$,  which contradicts  \eqref{nonvanishing4}. The lemma is proved.
\epr

\begin{lemma}  \label{mk T_12 cap mk T_43 3} For any $p,q\in \ZZ$ we have $ (a^p,b^{p+1},a^{p+1})   \cap  (b^q,a^q,b^{q+1})=\emptyset$. 
\end{lemma}
\bpr Let $\sigma \in  (a^p,b^{p+1},a^{p+1})$, then in table \eqref{no M} we see that $a^p ,a^{p+1} \in \sigma^{ss}$ and $\phi\left ( a^{p}\right)+1 < \phi\left (a^{p+1}\right)$. Now by   Lemma \ref{from phases of big distance}  we see that $a^q \not \in \sigma^{ss}$ for $q \not \in \{p,p+1\}$, and therefore $\sigma \not \in  (b^q,a^q,b^{q+1})$ for  $q \not \in \{p,p+1\}$. 

Suppose that $\sigma \in  (b^p,a^p,b^{p+1})$, then from table \eqref{no M} we obtain  $\phi\left ( b^{p}\right)+1 <\phi\left ( a^{p}\right)+1 < \phi\left (a^{p+1}\right)$, hence $\hom^1(a^{p+1},b^p)=0$, which contradicts \eqref{nonvanishing4}. 

Suppose that $\sigma \in  (b^{p+1},a^{p+1},b^{p+2})$, then from table \eqref{no M} we obtain  $\phi\left ( a^{p}\right)+1 < \phi\left (a^{p+1}\right)<\phi\left (b^{p+2}\right)$, hence $\hom^1(b^{p+2},a^p)=0$, which contradicts \eqref{nonvanishing3}. The lemma is proved. 
\epr
\begin{lemma}  \label{mk T_12 cap mk T_43 4} For any $p,q \in \ZZ$  we have:
$
\left (  (M',\_,\_)  \cup (\_,\_,M)\right ) \cap  (b^q,a^q,b^{q+1})=\emptyset $ and \\ 
$ \left (  (M,\_,\_) \cup (\_,\_,M')\right ) \cap  (a^p,b^{p+1},a^{p+1}) =\emptyset $.
 \end{lemma}
\bpr Assume first  that $\sigma \in  (b^q,a^q,b^{q+1})$, then table \eqref{no M} shows that:
\begin{gather}\label{big distance1} \phi(b^q)+1 <\phi(b^{q+1}).\end{gather}
Suppose that $\sigma \in (\_,\_,M)$, then using \eqref{big distance1},  $\hom(M,b^q)\neq 0$ and  table \eqref{left right M} we see that $\phi(a^m)+1<\phi(b^{q+1})$ and $\phi(a^{m+1})<\phi(b^{q+1})$ for some $m\in \ZZ$, hence $\hom^1(b^{q+1},a^m)=\hom(b^{q+1},a^{m+1})=0$, which contradicts  \eqref{nonvanishing3}.  

Suppose that $\sigma \in (M',\_,\_)$, then using \eqref{big distance1},  $\hom^1(b^{q+1},M')\neq 0$ and  table \eqref{left right M} we see that $\phi(b^q)+1<\phi(a^{j+1})$ and $\phi(b^{q})<\phi(a^{j})$ for some $j\in \ZZ$, hence $\hom^1(a^{j+1},b^q)=\hom(a^{j},b^{q})=0$, which contradicts  \eqref{nonvanishing4}.  So far we proved that $
\left (  (\_,\_,M) \cup (M',\_,\_) \right ) \cap  (b^q,a^q,b^{q+1})=\emptyset $. 

 Assume now   that $\sigma \in  (a^p,b^{p+1},a^{p+1})$, then table \eqref{no M} shows that:
\begin{gather}\label{big distance2} \phi(a^p)+1 <\phi(a^{p+1}).\end{gather}
Suppose that $\sigma \in (\_,\_,M')$, then using \eqref{big distance2},  $\hom(M',a^p)\neq 0$ and  table \eqref{left right M} we see that $\phi(b^i)+1<\phi(a^{p+1})$ and $\phi(b^{i+1})<\phi(a^{p+1})$ for some $i\in \ZZ$, hence $\hom^1(a^{p+1},b^i)=\hom(a^{p+1},b^{i+1})=0$, which contradicts  \eqref{nonvanishing4}.  

Suppose that $\sigma \in (M,\_,\_)$, then using \eqref{big distance2},  $\hom^1(a^{p+1},M)\neq 0$ and  table \eqref{left right M} we see that $\phi(a^p)+1<\phi(b^{n+1})$ and $\phi(a^{p})<\phi(b^{n})$ for some $n\in \ZZ$, hence $\hom^1(b^{n+1},a^p)=\hom(b^{n},a^{p})=0$, which contradicts  \eqref{nonvanishing3}.  Thus, we proved the second equality as well. 
 \epr

Lemmas \ref{mk T_12 cap mk T_43 1}, \ref{mk T_12 cap mk T_43 3}, \ref{mk T_12 cap mk T_43 4}, and formulas \eqref{T_12new}, \eqref{T_43new}  imply that $\mk{T}_{a}^{st}\cap \mk{T}_{b}^{st}=\emptyset$.

\subsection{The subsets \texorpdfstring{$(\_,\_,M), (\_,\_,M')$}{\space}, \texorpdfstring{$(M,\_,\_) \ \mbox{and} \ (M',\_,\_)$}{\space} are contractible
} \mbox{}\\

We start with:
\begin{lemma} \label{(_,_,X)0}  Let $\{x^i\}_{i\in \ZZ}$ be  either the sequence $\{a^i\}_{i\in \ZZ}$ or   the sequence $\{b^i\}_{i\in \ZZ}$.  If $m>j$ then: \begin{gather}(x^m, x^{m+1},X)\cap (x^j, x^{j+1},X) =\left \{ \sigma: \begin{array}{c} x^m \in \sigma^{ss}\\  x^{m+1} \in \sigma^{ss} \\ X \in \sigma^{ss}\end{array} , \  \begin{array}{c} 0 <\phi(x^{m+1}) - \phi(x^m) <1 \\ \phi(x^m)  < \phi(X) \\ \phi(x^{m+1})  < \phi(X)+1  \end{array}  \right \}, \nonumber \end{gather}
where $X=M$ if $\{x^i\}_{i\in \ZZ}=\{a^i\}_{i\in \ZZ}$, and $X=M'$ if $\{x^i\}_{i\in \ZZ}=\{b^i\}_{i\in \ZZ}$. 

In particular, $(x^m, x^{m+1},X)\cap (x^j, x^{j+1},X)$ and $(x^m, x^{m+1},X)\cup (x^j, x^{j+1},X)$ are contractible.  \end{lemma}\bpr
We show first  the inclusion $\subset$.  Assume that  $\sigma \in (x^m, x^{m+1},X)\cap (x^j, x^{j+1},X)$ and $m>j$. Then $X,x^{m+1},x^m, $ $x^{j+1},x^j$  are all semistable and by table \eqref{left right M} we have   

\begin{gather} \begin{array}{c} \phi ( x^m) < \phi (x^{m+1}) \\ \phi ( x^m) < \phi\left (X\right) \\ \phi ( x^{m+1}) < \phi\left (X\right)+1\end{array}   \quad \begin{array}{c} \phi( x^j) < \phi (x^{j+1}) \\ \phi ( x^j) < \phi (X) \\ \phi ( x^{j+1} ) < \phi\left (X\right)+1\end{array}.  \end{gather}
By $m>j$ it follows $\hom^1(x^{m+1}, x^j)\neq 0$, hence $\phi(x^{m+1})\leq \phi(x^j)+1$ (see Remark \ref{inc dec seq} (b)). On the other hand from the inequalities above we have  $\phi(x^j)+1<\phi(x^{j+1})+1$ and  by  Remark \ref{inc dec seq} (a) we obtain $\phi(x^{j+1})+1\leq \phi(x^m)+1$. Thus we obtain  $\phi(x^{m+1})< \phi(x^m)+1$ and $\subset$ follows.

We show now $\supset$. The condition defining the set on the right-hand side is the same as  $\sigma \in (x^m, x^{m+1},X)$      and  $ \phi(x^m) > \phi(x^{m+1}[-1])$ (see table \eqref{left right M}). From Proposition \ref{phi_1>phi_2} (a) it follows that $\mc A \cap \mc T_{exc}\subset \sigma^{ss}$, where $\mc A$ is the extension closure of $(x^m, x^{m+1}[-1])$, hence   By  Remark \ref{ext closure ab} we have $\{x^{j+1},x^j\}\subset \sigma^{ss}$.  The inequality  $0 <\phi(x^{m+1}) - \phi(x^m) <1$ and \eqref{phase formula}  show that $Z(x^{m+1})$, $Z(x^m)$ are not collinear, hence by  Corollary \ref{from equal phases to} (b) we get $\phi(x^{j+1}) \neq \phi(x^j)$. Now by Ramark \ref{inc dec seq} (a) and the incidence  $\sigma \in (x^m, x^{m+1},X)$   we get: $  \begin{array}{c} \phi ( x^j) < \phi(x^{j+1}) \\ \phi ( x^j) \leq \phi ( x^m)< \phi (X) \\ \phi ( x^{j+1})\leq \phi ( x^{m+1}) < \phi (X)+1\end{array}     $. In table \eqref{left right M} we see that $\sigma \in (x^j, x^{j+1}, X)$ and  the inclusion $\supset$ is shown.

The proved equality  implies that   $(x^m, x^{m+1}, X) \cap  (x^j, x^{j+1}, X)$ is contractible (see the arguments for the proof that \eqref{Theta_E n=2 3} is contractible in Proposition \ref{mutations}). Since $(x^m, x^{m+1}, X)$ and $  (x^j, x^{j+1}, X)$ are contractible, by Remark \ref{VK} we see that $(x^m, x^{m+1}, X) \cup  (x^j, x^{j+1}, X)$ is contractible as well.
\epr

\begin{coro} \label{(_,_,X)} The subsets $(\_,\_,M)$ and $(\_,\_,M')$ of $\st(D^b(Q ))$ are  contractible. 
\end{coro}
\bpr Recalling \eqref{(_,_,M)} and using the notations of the previous lemma, we see that we have to show that $ \bigcup_{j\in \ZZ} (x^{j}, x^{j+1},X)$ is contractible. It is shown in   Lemma \ref{(_,_,X)0}, that  for a given $m\in \ZZ$ the intersection $(x^m, x^{m+1},X)\cap (x^j, x^{j+1},X)$  is contractible and  it is  the same for all $j<m$. Now by induction and using   Remark \ref{VK} one shows that $\bigcup_{k=0}^n(x^{m-k}, x^{m-k+1},X)$ is contractible for any $n\in \NN$ and any $m\in \ZZ$.  Using  again Remark \ref{VK} we deduce that $ \bigcup_{j\in \ZZ} (x^{j}, x^{j+1},X)$ is contractible.  The corollary follows.
\epr

 The steps in  the proof that   $(M,\_,\_)$ and $(M',\_,\_)$  are analogous. We show first:

\begin{lemma} \label{(X,_,_)0} Let $\{x^i\}_{i\in \ZZ}$ be  either the sequence $\{a^i\}_{i\in \ZZ}$ or   the sequence $\{b^i\}_{i\in \ZZ}$ . If  $m<j$, then:  \begin{gather}\label{(X,_,_)a}  (X,x^m, x^{m+1})\cap (X,x^j, x^{j+1}) =\left \{ \sigma:\begin{array}{c} X\in \sigma^{ss}\\ x^m\in \sigma^{ss} \\ x^{m+1}  \in \sigma^{ss}\end{array},  \  \begin{array}{c} \phi(X)<\phi(x^m) \\ \phi(X) +1 < \phi(x^{m+1}) \\  0 <\phi(x^{m+1}) - \phi(x^m) <1  \end{array}  \right \}  \end{gather}
where $X=M'$ if $\{x^i\}_{i\in \ZZ}=\{a^i\}_{i\in \ZZ}$, and $X=M$ if $\{x^i\}_{i\in \ZZ}=\{b^i\}_{i\in \ZZ}$. 

In particular, $ (X,x^m, x^{m+1})\cap (X,x^j, x^{j+1})$ and  $(X,x^m, x^{m+1})\cup (X,x^j, x^{j+1})$ are contractible.  \end{lemma}
\bpr  By table \eqref{left right M} we see that  the condition defining the set on the right-hand side of \eqref{(X,_,_)a}  is the same as  $\sigma \in (X,x^m, x^{m+1})$  and $\phi(x^m)>\phi(x^{m+1}[-1])$.  

  The inclusion $\subset$   follows from table \eqref{left right M},  $\hom^1(x^{j+1},x^m) \neq 0$ and Remark \ref{inc dec seq} (a) as follows $\phi(x^{m+1})\leq \phi(x^j)<\phi(x^{j+1})\leq \phi(x^m)+1$.  

 To show the converse  inclusion $\supset$ in  \eqref{(X,_,_)a},  assume that  $\sigma \in (X,x^m, x^{m+1})$  and $\phi(x^m)>\phi(x^{m+1}[-1])$.     From Proposition \ref{phi_1>phi_2} (b) and   Remark \ref{ext closure ab} it follows that $x^j, x^{j+1} \in \sigma^{ss}$. Since we have  $0 <\phi(x^{m+1}) - \phi(x^m) <1$, it follows that $Z(x^m)$, $Z(x^{m+1})$ are not collinear, therefore by Corollary \ref{from equal phases to} (b) and   Remark \ref{inc dec seq} (a) we obtain $\phi(x^j)<\phi(x^{j+1})$.  Since $j>m$, by Remark \ref{inc dec seq} (a) we obtain also $\phi(X)<\phi(x^m)\leq  \phi(x^j)$, $ \phi(X) +1 < \phi(x^{m+1}) \leq \phi(x^{j+1})$, hence $\sigma \in  (X,x^j, x^{j+1})$.
 
The proved equality  implies that   $ (X,x^m, x^{m+1})\cap (X,x^j, x^{j+1})$  is contractible (see the arguments for the proof that \eqref{Theta_E n=2 3} is contractible in Proposition \ref{mutations}). Since  $ (X,x^m, x^{m+1})$ and  $(X,x^j, x^{j+1})$ are contractible, by Remark \ref{VK} we see that  $ (X,x^m, x^{m+1})\cup (X,x^j, x^{j+1})$ is contractible as well. 
\epr 

\begin{coro} \label{(X,_,_)} The subsets $(M,\_,\_), (M',\_,\_)\subset \st(D^b(Q ))$ are  contractible. 
\end{coro} 
\bpr  Recalling \eqref{(M',_,_)} and using the notations of the previous lemma, we see that we have to show that $ \bigcup_{j\in \ZZ} (X,x^{j}, x^{j+1})$ is contractible.  From Lemma \ref{(X,_,_)0}  we know  that for a given $m\in \ZZ$ the intersection $(X,x^m, x^{m+1})\cap (X,x^j, x^{j+1})$ is contractible and it is  the same for all $j>m$.  Now by induction and using   Remark \ref{VK} one shows that $\bigcup_{k=0}^n(X,x^{m+k}, x^{m+k+1})$ is contractible for any $n\in \NN$ and any $m\in \ZZ$.  Using  again Remark \ref{VK} we deduce that $ \bigcup_{j\in \ZZ} (X,x^{j}, x^{j+1})$ is contractible.  The corollary follows.
\epr

\section{The subsets \texorpdfstring{$\mk{T}_{a}^{st}$}{\space} and \texorpdfstring{$\mk{T}_{b}^{st}$}{\space} are contractible} \label{T_a T_b cont}

We start by showing  some empty intersections:  
\begin{lemma} \label{T12lemma1} The unions $\bigcup_{p\in \ZZ}(a^p,b^{p+1},a^{p+1})$  and $\bigcup_{p\in \ZZ}(b^p,a^p,b^{p+1})$ are   disjoint. Furthermore, we have: 
\begin{gather} \label{T_12c}  p\neq q \ \Rightarrow \ (a^p,b^{p+1},a^{p+1}) \cap (a^q,a^{q+1},M)=(a^p,b^{p+1},a^{p+1}) \cap (M',a^q,a^{q+1})=\emptyset\\ \label{T_43c} p\neq q \ \Rightarrow \ (b^p,a^p,b^{p+1}) \cap (b^q,b^{q+1},M')=(b^p,a^p,b^{p+1})  \cap (M,b^q,b^{q+1})=\emptyset.     \end{gather}   \end{lemma}
\bpr 
If $\sigma \in  (a^p,b^{p+1},a^{p+1})$, then these exceptional objects are semistable and  by table \eqref{no M}  we have  $\phi(a^p)+1 < \phi(a^{p+1})$.  Now by Lemma \ref{from phases of big distance}   we see that $a^{j}$ with $j\not \in \{p,p+1\}$ can not  be semistable, therefore $\sigma \not \in  (a^q,b^{q+1},a^{q+1})$, $\sigma \not \in  (a^q,a^{q+1},M)$, and  $\sigma \not \in   (M',a^q,a^{q+1})$  for $q\neq p$.

If $\sigma \in   (b^p,a^p,b^{p+1}) $, then $b^p,a^p,b^{p+1}$ are semistable and  by table \eqref{no M}  we have  $\phi(b^p)+1 < \phi(b^{p+1})$.  Now by Lemma \ref{from phases of big distance}  we see that  $b^{j}$ with $j\not \in \{p,p+1\}$   can not  be semistable, therefore $\sigma \not \in  (b^q,a^q,b^{q+1})$, $\sigma \not \in  (b^q,b^{q+1},M')$, and  $\sigma \not \in  (M,b^q,b^{q+1})$ for $q\neq p$.
\epr
Now we attach the  pairwise non-intersecting contractible blocks  $\{ (a^p,b^{p+1},a^{p+1}) \}_{p\in\ZZ}$ to $(\_,\_,M)$ and to  $(M',\_,\_)$ 

\begin{lemma} \label{T12lemma3} For any $p\in \ZZ$  the sets $ (a^p,b^{p+1},a^{p+1})  \cap  (\_,\_,M)$;    $  (a^p,b^{p+1},a^{p+1})  \cap  (M',\_,\_)$;   $ (b^p,a^p,b^{p+1}) \cap  (\_,\_,M')$; and    $  (b^p,a^p,b^{p+1})  \cap  (M,\_,\_)$ are non-empty and contractible.
\end{lemma}
\bpr  From    \eqref{(M',_,_)}, \eqref{(_,_,M)}  and \eqref{T_12c}  it follows that: $(a^p,b^{p+1},a^{p+1}) \cap (\_,\_,M) = (a^p,b^{p+1},a^{p+1})\cap (a^{p}, a^{p+1}, M) $ and $(a^p,b^{p+1},a^{p+1})\cap (M',\_,\_) =(a^p,b^{p+1},a^{p+1})\cap (M',a^p,a^{p+1})$.
 
 From    \eqref{(M',_,_)}, \eqref{(_,_,M)}  and \eqref{T_43c}  it follows that: $(b^p,a^p,b^{p+1}) \cap  (\_,\_,M') = (b^p,a^p,b^{p+1})\cap (b^{p}, b^{p+1}, M') $ and $ (b^p,a^p,b^{p+1})  \cap  (M,\_,\_) =(b^p,a^p,b^{p+1})\cap (M,b^p,b^{p+1})$.

 From Proposition \ref{mutations} it follows that $(a^p,b^{p+1},a^{p+1})\cap (a^{p}, a^{p+1}, M) $, $(a^p,b^{p+1},a^{p+1})\cap (M',a^p,a^{p+1})$, $(b^p,a^p,b^{p+1})\cap (b^{p}, b^{p+1}, M') $,  and $(b^p,a^p,b^{p+1})\cap (M,b^p,b^{p+1})$  are contractible. 

  The lemma follows. 
\epr  
Let us denote: 
\begin{gather} \label{T12Z} Z= (M',\_,\_) \cup  \bigcup_{p\in \ZZ}(a^p,b^{p+1},a^{p+1}).  \end{gather}
Corollary \ref{(X,_,_)} and  Lemmas \ref{T12lemma1}, \ref{T12lemma3}  imply  (recall Remark \ref{VK}) that $Z$ is contractible. From \eqref{T_12new} and \eqref{(_,_,M)} we see that: \begin{gather} \label{T12=Zcup} \mk{T}_{a}^{st} =Z \cup (\_,\_,M)= Z  \cup  \bigcup_{m\in \ZZ} (a^{m}, a^{m+1}, M).\ \end{gather} 
We start to glue the contractible summands in formula \eqref{T12=Zcup}. The first step is:
\begin{lemma} \label{T12Zcap(E_1)} The set $ (a^{m}, a^{m+1}, M) \cap Z$ consists of the stability conditions $\sigma$ for which  $ a^{m},a^{m+1}, M $ are semistable and:   \begin{gather}   \label{T12Zcap(E_1)ineq} \begin{array}{c}  \phi(a^{m})-1<\phi(a^{m+1}[-1])< \phi(a^{m})\\ \phi(a^{m})-1<\phi(M[-1])< \phi(a^{m}) \\ \arg_{(\phi(a^{m})-1,\phi(a^{m}))}(Z(a^{m})- Z(a^{m+1}))> \phi(M)-1  \end{array}    \ \mbox{or} \  \begin{array}{c}\phi(M)<\phi(a^{m+1}) \\ \phi\left ( a^{m}\right) < \phi\left (a^{m+1}\right) \\ \phi\left ( a^{m}\right) < \phi\left (M\right) \\ \phi\left ( a^{m+1}\right) < \phi\left (M\right)+1\end{array}.   \end{gather}
It follows that  $ (a^{m}, a^{m+1}, M) \cap Z$ and  $ (a^{m}, a^{m+1}, M) \cup Z$  are  contractible. 
\end{lemma}
\bpr In  \eqref{T_12c} we  have that  $(a^{m}, a^{m+1}, M)\cap (a^{j}, b^{j+1} , a^{j+1})=\emptyset$ for $j\neq m$.  Therefore (recall \eqref{T12Z}) 
 \begin{gather}\label{T12Zcap(E_1)0}  (a^{m}, a^{m+1}, M) \cap Z = (a^{m}, a^{m+1}, M) \cap ( (a^{m},b^{m+1} , a^{m+1})\cup (M',\_,\_)).\end{gather}

We show first the inclusion $\subset$. Assume that  $\sigma \in  (a^{m}, a^{m+1}, M) \cap Z$. Then $a^{m}, a^{m+1}, M$ are semistable and    from table \eqref{left right M} we see that 
\begin{gather}\label{T12Zcap(E_1)1} \begin{array}{c}  \phi\left ( a^{m}\right) < \phi\left (a^{m+1}\right) \\ \phi\left ( a^{m}\right) < \phi\left (M\right) \\ \phi\left ( a^{m+1}\right) < \phi\left (M\right)+1\end{array}. \end{gather} Taking into account  \eqref{T12Zcap(E_1)0} we consider two  cases.

\ul{If $\sigma \in  (a^{m},b^{m+1} , a^{m+1})$,} then $b^{m+1} \in \sigma^{ss}$ and in table \eqref{no M} we see that $\phi(b^{m+1})<\phi(a^{m+1})$. From $\hom(M,b^{m+1})\neq 0$ (see \eqref{nonvanishing1}) it follows that $\phi(M)\leq \phi(b^{m+1})<\phi(a^{m+1})$ and we obtain the second system of inequalities in \eqref{T12Zcap(E_1)ineq}.

\ul{If $\sigma \in  (M', a^{j}, a^{j+1} )$,} then $M', a^{j}, a^{j+1} \in \sigma^{ss}$ and in table \eqref{left right M} we see that $\phi(M')+1<\phi(a^{j+1})$ and  $\phi(a^{j})<\phi(a^{j+1})$. From $\hom^1(M,M')\neq 0$ it follows that $\phi(M)\leq \phi(M')+1<\phi(a^{j+1})$. Since $\phi(a^{m})<\phi(M)$,   it follows from Remark \ref{inc dec seq} (a) that  $m\leq j$. 

If $m=j$,  then $\phi(M)<\phi(a^{m+1})$ and we obtain the second system of inequalities. 

If $m<j$,   then we show  that the  first system of inequalities in \eqref{T12Zcap(E_1)ineq} holds.   Now $\phi(M)<\phi(a^{j+1})$ and  $\hom^1(a^{j+1}, a^{m})\neq 0$,  hence $\phi(a^{m+1})\leq \phi(a^{j})<\phi(a^{j+1})\leq \phi(a^{m})+1$ and $\phi(M)<\phi(a^{j+1})\leq \phi(a^{m})+1$.  We have also  $M'\in \sigma^{ss}$ and by $\hom^1(M,M')\neq 0$ it follows that    $\phi(M)\leq \phi(M')+1$. 
From $\hom(M',a^{m})\neq 0$ it follows  $\phi(M')\leq a^{m}$.   These  arguments   together with \eqref{T12Zcap(E_1)1} imply:  
\begin{gather} \label{T12Zcap(E_1)3} \begin{array}{c}  \phi(a^{m})-1<\phi(a^{m+1}[-1])< \phi(a^{m});\\ \phi(a^{m})-1<\phi(M[-1])\leq \phi(M') \leq  \phi(a^{m}) ; \ \ \ \ \ \phi(M[-1]) <  \phi(a^{m}).  \end{array} \end{gather}

  In \eqref{Z(delta)}  we have   $Z(a^{m})-Z(a^{m+1})=Z(\delta)$, therefore it remains to show that:
\begin{gather} \label{T12Zcap(E_1)4} \arg_{(\phi(a^{m})-1,\phi(a^{m}))}(Z(\delta))> \phi(M)-1.    \end{gather}   
From the second row of \eqref{T12Zcap(E_1)3} and \eqref{phase formula} we see that $Z(\delta)$ and  $Z(M[-1])$ both lie  in the half-plane\footnote{The notation  $Z(a^{m})^c_-$ is explained in \eqref{complement of a line}.}  $Z(a^{m})^c_-$. In \eqref{[delta]}    we have  aslo $Z(M')=Z(\delta) + Z(M[-1])$, therefore the  vector $Z(M')$   is  in $Z(a^{m})^c_-$ as well, hence by $Z(M')=Z(\delta) + Z(M[-1])$ it follows that  the inequality   \eqref{T12Zcap(E_1)4}  is equivalent to  $\phi(M')>\phi(M[-1]) $. Therefore it remains  to show that   $\phi(M') \neq \phi(M[-1]) $. Indeed, on one hand $\phi(M[-1]) = \phi(M')$ implies $ \arg_{(\phi(a^{m})-1,\phi(a^{m}))}(Z(\delta))= \phi(M')$. On the other hand,
   $\sigma \in  (M', a^{j}, a^{j+1} )$,   $m<j$ and   \eqref{T12Zcap(E_1)1} imply    $\phi(M')+1<  \phi(a^{j+1})\leq \phi(a^{m})+1 < \phi(M)+1\leq\phi(M')+2$.  Thus, we see that   $\phi(M[-1]) = \phi(M')$ implies  $Z(a^{j+1}) \in Z(\delta)^c_-$.  However, from the first inequality in  \eqref{T12Zcap(E_1)3}  and Corollary \ref{from equal phases to} (d)  it follows that $Z(a^{j+1}) \in Z(\delta)^c_+$, which is a contradiction,  and \eqref{T12Zcap(E_1)4} follows.

So far we showed that $\sigma \in (a^{m}, a^{m+1}, M) \cap Z$ implies  \eqref{T12Zcap(E_1)ineq}. We show now converse inclusion. 

We assume first that the second system of inequalities  in \eqref{T12Zcap(E_1)ineq}  holds. In particular $\sigma \in (a^{m}, a^{m+1}, M)$. By the inequality $\phi(M)<\phi(a^{m+1}) $ we can apply  Proposition  \ref{phi_1>phi_2} (b), hence the triangle \eqref{short filtration 2}  implies that $b^{m+1} \in \sigma^{ss}$, $\phi(M)\leq \phi(b^{m+1})\leq \phi(a^{m+1})$, and $Z(M)+Z(a^{m+1})=Z(b^{m+1})$. We have in \eqref{T12Zcap(E_1)ineq} also $\phi(a^{m+1})-1<\phi(M)<\phi(a^{m+1})$  and it follows that  $\phi(M)<\phi(b^{m+1})<\phi(a^{m+1})$.   If the inequality   $\phi(a^{m+1})>  \phi(a^{m})+1$ holds, then due to   $\phi(M)<\phi(b^{m+1})<\phi(a^{m+1})$ and  $\phi(a^{m})<\phi(M)$  hold we obtain that  $\sigma \in (a^{m},b^{m+1},a^{m+1})\subset Z$ (see table \eqref{no M}).  

Thus, we can assume  that  $\phi(a^{m+1}[-1])\leq  \phi(a^{m})$ and combining with   the inequalities   $\phi(M[-1])<\phi(a^{m+1}[-1])<\phi(M)$,  $\phi(a^{m})<\phi(M)$ (given in \eqref{T12Zcap(E_1)ineq}) we get $\phi(M[-1])< \phi(a^{m+1}[-1])\leq \phi(a^{m})<\phi(M)$. Now  it is easy to show with the help of  Corollaries \ref{nonvanishings} and \ref{one nonvanishing degree} that  $(a^{m}, a^{m+1}[-1], M[-1])$  is a $\sigma$-triple (see Definition \ref{sigma triple}).
  Combining  the triangles  \eqref{short filtration 1} and \eqref{short filtration 2} we get the following sequence:    \be \label{T12Zcap(E_1)the filtation} 
\begin{diagram}[size=1em] 
0 & \rTo      &        &       &   M[-1] & \rTo      &             &        & b^{m+1}[-1] &   \rTo      &           &       & M' \\
  & \luDashto &        & \ldTo &         &\luDashto  &             &  \ldTo &           &   \luDashto &           & \ldTo &       \\
  &           &  M[-1] &       &         &           &  a^{m+1}[-1]  &        &           &             & a^{m} &
\end{diagram}.
\ee
The conditions of Lemma \ref{three comp factors} (b) are satisfied with the triple  $(a^{m}, a^{m+1}[-1], M[-1])$ and the diagram above. Therefore $M'\in \sigma^{ss}$ and  $\phi(M')<\phi(a^{m})$. 

If   $\phi(a^{m+1}[-1])=\phi(a^{m})$, then it follows that  $\phi(M')+1<\phi(a^{m+1})$, and recalling that we have also $\phi(a^{m})<\phi(a^{m+1})$ we see that $\sigma \in  (M', a^{m}, a^{m+1} )\subset Z$ (see table \eqref{left right M}).

  Therefore we can assume that   $\phi(M[-1])< \phi(a^{m+1}[-1]) <  \phi(a^{m})<\phi(M)$. We will show in this case that $\sigma \in (M',a^{j}, a^{j+1})$ for some big enough $j$. From Proposition \ref{all exc semistable} and  Remark \ref{ext closure ab} it follows that $\{a^{j+1}\}_{j\in \ZZ}\subset \sigma^{ss}$.  From  $\phi(a^{m})<\phi(a^{m+1}) <  \phi(a^{m})+1 $  and Corollary \ref{from equal phases to} (b) and (d) we see that $\phi(a^{j})<\phi(a^{j+1})$ and $Z(a^j)\in Z(\delta)_+^c$ for each $j$ (recall also Remark \ref{inc dec seq} (a)). 
		We will show that for big enough $j$ we have  $\phi(M')+1< \phi(a^j) $ and then from table \eqref{left right M} we obtain  $\sigma \in (M',a^{j}, a^{j+1})\subset Z$.
		
		Now we have  $ \phi(a^{m}[-1])<\phi(M[-1])< \phi(a^{m+1}[-1]) <  \phi(a^{m})$. Recalling that $Z(\delta)=Z(a^{m})+Z(a^{m+1}[-1])$, we see that we can choose $t\in \RR$ so that $t< \phi(a^{m})<\phi(M)< \phi(a^{m+1})<t+1 $ and $Z(\delta)=\abs{Z(\delta)} \exp(\ri \pi t)$.    Since $\hom^1(a^j,a^{m})\neq 0$, $\hom(a^{m},a^j)\neq 0$   for $j>m+1$ and  by Corollary \ref{from equal phases to} (b), we have  $\phi(a^{m})<\phi(a^j)<  \phi(a^{m})+1 $ for  $j>m+1$. These inequalities together with the incidence $Z(a^j)\in Z(\delta)_+^c$ imply that $\arg_{(t,t+1)}(Z(a^j))=\phi(a^j)$  for $j>m+1$ (see Remark \ref{arg remark}  (c)). Now the formula \eqref{noncolinear ab2} in Corollary \ref{noncolinear ab}  gives us the following equality:
		\begin{gather}\label{T12Zcap(E_1)the limit} \lim_{j\rightarrow  \infty} \phi(a^j)= t+1.\end{gather}
We showed that  $\phi(M')<\phi(a^{m})$ (see below \eqref{T12Zcap(E_1)the filtation}) and we have also $\phi(a^{m})<\phi(M)$. Using $\hom^1(M,M')\neq 0$ we see that $\phi(a^{m}[-1])<\phi(M')<\phi(a^{m})$. We showed also that   $t< \phi(a^{m})<\phi(M)< \phi(a^{m+1})<t+1 $. Since $Z(M)+Z(M')=Z(\delta)=\abs{Z(\delta)} \exp(\ri \pi t)$, it follows that $Z(M')\in  Z(\delta)^c_-$ and $\phi(M')<t$. By \eqref{T12Zcap(E_1)the limit} we get the desired  $\phi(a^j)>\phi(M')+1$ for big  $j$. 

So, we showed that the second system of inequalities in \eqref{T12Zcap(E_1)ineq} implies that $\sigma\in (a^{m}, a^{m+1}, M) \cap Z$.
We show now that the  first system  in \eqref{T12Zcap(E_1)ineq} implies  $\sigma\in (a^{m}, a^{m+1}, M) \cap Z$ as well. 
Assume that $a^{m}, a^{m+1}, M \in \sigma^{ss}$ and that these inequalities hold. They contain  the inequalities defining $(a^{m}, a^{m+1}, M)$ (see table \eqref{left right M}), therefore we obtain $\sigma \in (a^{m}, a^{m+1}, M)$ immediately.   Furthermore the first two inequalities show that  $(a^{m}, a^{m+1}[-1], M[-1])$ is a $\sigma$-triple. 
The conditions of Lemma \ref{three comp factors} (a) are satisfied with the triple  $(a^{m}, a^{m+1}[-1], M[-1])$ and the diagram \eqref{T12Zcap(E_1)the filtation}. Therefore $M'\in \sigma^{ss}$ and $\phi(M')<\phi(a^{m})$. By $\hom^1(M,M')\neq 0$ we can write 
$\phi(a^{m})-1< \phi(M[-1])\leq \phi(M') < \phi(a^{m})$, hence by  \eqref{phase formula} we see that  $Z(\delta), Z(M[-1]), Z(M') \in Z(a^{m})^c_-$. Let us denote $t=\arg_{(\phi(a^{m})-1,\phi(a^{m}))}(Z(\delta))$. The third inequality in \eqref{T12Zcap(E_1)ineq} is the same as $t> \phi(M)-1 $. Combining  these arguments with the equality   $Z(M')=Z(\delta)+Z(M[-1])$ we write:  
\begin{gather} \label{T12Zcap(E_1) help} \phi(a^{m}[-1])< \phi(M[-1])< \phi(M') <t< \phi(a^{m}).  \end{gather}
We will show that  $\sigma \in  (M', a^{j}, a^{j+1} )$  for some big enough $j$.  We have  $\phi(a^{m})<\phi(a^{m+1}) <  \phi(a^{m})+1 $ (the first inequality in \eqref{T12Zcap(E_1)ineq}), which by     Proposition \ref{all exc semistable} and  Remark \ref{ext closure ab} implies  that $\{a^{j+1}\}_{j\in \ZZ}\subset \sigma^{ss}$,   and by Corollary \ref{from equal phases to} (b), (d) implies  that $\phi(a^{j})<\phi(a^{j+1})$ and $Z(a^j)\in Z(\delta)_+^c$ for each $j$.   Since $\hom^1(a^j,a^{m})\neq 0$, $\hom(a^{m},a^j)\neq 0$   for $j>m+1$, we have  $\phi(a^{m})<\phi(a^j)<  \phi(a^{m})+1 $ for  $j>m+1$. These inequalities together with   the incidences $Z(a^j)\in Z(\delta)_+^c$,  $\phi(a^m)\in(t,t+1)$  imply that $\arg_{(t,t+1)}(Z(a^j))=\phi(a^j)$  for $j>m+1$ (see Remark \ref{arg remark}  (c)). Now  the formula \eqref{noncolinear ab2} in Corollary \ref{noncolinear ab}  leads to  \eqref{T12Zcap(E_1)the limit} again.   Therefore by \eqref{T12Zcap(E_1) help} we see that $\phi(a^j)>\phi(M')+1$ for big enough $j$. It follows that $\sigma \in  (M', a^{j}, a^{j+1} )\subset Z$  (see table \eqref{left right M}) . 

The first part of the lemma is shown.
 It is easy now to show that the  intersection is contractible. The  intersection in quiestion  is the same as $ (a^{m}, a^{m+1}[-1], M[-1]) \cap Z$. Let us denote $\mc E= (a^{m}, a^{m+1}[-1], M[-1])$.  We have a homeomorphism ${f_{\mc E}}_{\vert \Theta_{\mc E}}:\Theta_{\mc E} \rightarrow f_{\mc E}(\Theta_{\mc E})$ (see \eqref{A and homeo}, \eqref{the map}). The proved description of $Z\cap \Theta_{\mc E}$ by the inequalities \eqref{T12Zcap(E_1)ineq}  shows that $f_{\mc E}(Z\cap \Theta_{\mc E}) $ is union of two sets.  The first set after permutation of the coordinates in $\RR^6$ is  the same as the set considered in  Corollary \ref{contractible 1b}, hence it is also contractible.  The second is obviously contractible.   Furthermore, one easily shows that  the intersection of these two sets  is
$ \RR_{>0}^3 \times  \left \{  \phi_0-1 < \phi_2 <\phi_1 < \phi_0 \right \},$ which is contractible as well. Now by Remark \ref{VK} it follows  that $f_{\mc E}(Z\cap \Theta_{\mc E}) $ is contractible, therefore $Z\cap \Theta_{\mc E}$ is contractible as well. Recalling that $Z$ and $ (a^{m}, a^{m+1}, M) $ are contractible and applying again Remark \ref{VK} we deduce that $ (a^{m}, a^{m+1}, M) \cup Z$  is contractible.
 The lemma is proved. 
\epr 

\begin{coro} \label{T12Zcap union(E_1)} The set $\mk{T}_{a}^{st}$ is contractible.  
\end{coro}
\bpr Recall that   $\mk{T}_{a}^{st}= Z \cup \bigcup_{j\in \NN}  (a^{j}, a^{j+1}, M) $  (see \eqref{T12=Zcup}). We will show  that  $Z \cup \bigcup_{j=0}^n  (a^{m-j+1}, a^{m-j}, M) $ is contractible  for each $m\in \ZZ$ and each $n\in \NN$. Then the corollary follows from Remark \ref{VK}. 

 Assume that for some $n\in \NN$ the set $Z \cup \bigcup_{j=0}^n  (a^{m-j+1}, a^{m-j}, M) $ is contractible  for each $m\in \ZZ$. We have shown this statement  for $n=0$   in Lemma \ref{T12Zcap(E_1)}, and now we make induction assumption. 
Take any $m\in \NN$ and consider  $Z \cup \bigcup_{j=0}^{n+1} (a^{m-j+1}, a^{m-j}, M)$ $= \left (Z \cup  \bigcup_{j=1}^{n+1}  (a^{m-j+1}, a^{m-j}, M)\right ) \cup  (a^{m}, a^{m+1}, M)$. By the induction assumption $Z \cup  \bigcup_{j=1}^{n+1}  (a^{m-j+1}, a^{m-j}, M) $ and $(a^{m}, a^{m+1}, M)$ are contractible. We will show now that the intersection of these sets  is contractible as well and then by Remark \ref{VK} we obtain that the union $Z \cup \bigcup_{j=0}^{n+1}  (a^{m-j+1}, a^{m-j}, M)$ is contractible. Indeed,  we have 
\begin{gather}  \left (Z \cup  \bigcup_{j=1}^{n+1}  (a^{m-j+1}, a^{m-j}, M) \right ) \cap  (a^{m}, a^{m+1}, M) = \nonumber \\[-2mm] \label{T12Zcap union(E_1)1}  \\[-2mm]\left  ( (a^{m}, a^{m+1}, M) \cap Z \right ) \cup  \left ( (a^{m}, a^{m+1}, M) \cap  \bigcup_{j=1}^{n+1}   (a^{m-j+1}, a^{m-j}, M)  \right ).\nonumber \end{gather}
Using Lemmas \ref{T12Zcap(E_1)} and \ref{(_,_,X)0} we deduce that the considered intersection consists of  the stability conditions for which 
 $a^{m}, a^{m+1}, M$ are semi-stable and some of the two systems of inequalities in \eqref{T12Zcap(E_1)ineq} or the system    $\begin{array}{c} \phi(a^{m})<\phi(a^{m+1})<\phi(a^{m})+1 \\ \phi(a^{m})<\phi(M) \\ \phi(a^{m+1})<\phi(M)+1\end{array}$ holds. Since the first system in \eqref{T12Zcap(E_1)ineq}  implies the last system we deduce that the intersection \eqref{T12Zcap union(E_1)1} is described by the inequalities:
 \begin{gather}    \begin{array}{c} \phi(a^{m})<\phi(a^{m+1})<\phi(a^{m})+1 \\ \phi(a^{m})<\phi(M) \\ \phi(a^{m+1})<\phi(M)+1\end{array}   \ \mbox{or} \  \begin{array}{c} \phi\left ( a^{m}\right) < \phi\left (a^{m+1}\right) \\ \phi\left ( a^{m}\right) < \phi\left (M\right) \\ \phi(M)< \phi\left ( a^{m+1}\right) < \phi\left (M\right)+1\end{array}.   \end{gather}
Now analogous  arguments as in the last paragraph of the proof of Lemma \ref{T12Zcap(E_1)} show that the intersection \eqref{T12Zcap union(E_1)1} is contractible. The corollary follows.
\epr

We pass to the proof that $ \mk{T}_b^{st}$ is contractible. Let us denote 
\begin{gather} \label{T43Z} W= (M,\_,\_) \cup  \bigcup_{p\in \ZZ}(b^p,a^{p},a^{p+1}).  \end{gather}
Corollary \ref{(X,_,_)} and  Lemmas \ref{T12lemma1}, \ref{T12lemma3}  imply  (recall Remark \ref{VK}) that $W$ is contractible. From \eqref{T_43new} and \eqref{(_,_,M)} we see that: \begin{gather} \label{T43=Zcup} \mk{T}_{b}^{st} =W \cup (\_,\_,M')= W  \cup  \bigcup_{m\in \ZZ} (b^{m}, b^{m+1}, M').\ \end{gather} 
 The proof  of the next Lemma \ref{T43Zcap(E_1)} is analogous to the proof of Lemma \ref{T12Zcap(E_1)}):
\begin{lemma} \label{T43Zcap(E_1)} The set $ (b^{m}, b^{m+1}, M') \cap W$ consists of the stability conditions $\sigma$ for which  $ b^{m},b^{m+1}, M' $ are semistable and   \begin{gather}   \label{T43Zcap(E_1)ineq} \begin{array}{c}  \phi(b^{m})-1<\phi(b^{m+1}[-1])< \phi(b^{m})\\ \phi(b^{m})-1<\phi(M'[-1])< \phi(b^{m}) \\ \arg_{(\phi(b^{m})-1,\phi(a^{m}))}(Z(b^{m})- Z(b^{m+1}))> \phi(M')-1  \end{array}    \ \mbox{or} \  \begin{array}{c}\phi(M')<\phi(b^{m+1}) \\ \phi\left ( b^{m}\right) < \phi\left (b^{m+1}\right) \\ \phi\left ( b^{m}\right) < \phi\left (M'\right) \\ \phi\left ( b^{m+1}\right) < \phi\left (M'\right)+1\end{array}.   \end{gather}
It follows that  $ (b^{m}, b^{m+1}, M') \cap W$ and  $ (b^{m}, b^{m+1}, M') \cup W$  are  contractible. 
\end{lemma}
\bpr In  \eqref{T_43c} we  have that  $(b^{m}, b^{m+1}, M')\cap (b^{j}, a^{j} , b^{j+1})=\emptyset$ for $j\neq m$.  Therefore (recall \eqref{T43Z}) 
 \begin{gather}\label{T43Zcap(E_1)0}  (b^{m}, b^{m+1}, M') \cap W = (b^{m}, b^{m+1}, M') \cap ( (b^{m},a^{m} , b^{m+1})\cup (M,\_,\_)).\end{gather}

We show first the inclusion $\subset$. Assume that  $\sigma \in  (b^{m}, b^{m+1}, M') \cap W$. Then $b^{m}, b^{m+1}, M'$ are semistable and    from table \eqref{left right M} we see that 
\begin{gather}\label{T43Zcap(E_1)1} \begin{array}{c}  \phi\left ( b^{m}\right) < \phi\left (b^{m+1}\right) \\ \phi\left ( b^{m}\right) < \phi\left (M'\right) \\ \phi\left ( b^{m+1}\right) < \phi\left (M'\right)+1\end{array}. \end{gather} Taking into account  \eqref{T43Zcap(E_1)0} we consider two  cases.

\ul{If $\sigma \in  (b^{m},a^{m} , b^{m+1})$,} then $a^{m} \in \sigma^{ss}$ and  $\phi(a^{m})<\phi(b^{m+1})$ (see table \eqref{no M}). From $\hom(M',a^{m})\neq 0$ (see \eqref{nonvanishing1}) it follows $\phi(M')\leq \phi(a^{m})<\phi(b^{m+1})$ and we get the second system  in \eqref{T43Zcap(E_1)ineq}.

\ul{If $\sigma \in  (M, b^{j}, b^{j+1} )$,} then $M, b^{j}, b^{j+1}\in \sigma^{ss}$ and in table \eqref{left right M} we see that $\phi(M)+1<\phi(b^{j+1})$ and  $\phi(b^{j})<\phi(b^{j+1})$. From $\hom^1(M',M)\neq 0$ it follows that $\phi(M')\leq \phi(M)+1<\phi(b^{j+1})$. Since we have $\phi(b^{m})<\phi(M')$,   Remark \ref{inc dec seq} (a) implies that  $m\leq j$. 

If $m=j$,  then $\phi(M')<\phi(b^{m+1})$ and we obtain the second system of inequalities. 

If $m<j$,   then we will show  that the  first system of inequalities in \eqref{T43Zcap(E_1)ineq} holds.   Now $\phi(M')<\phi(b^{j+1})$ and  $\hom^1(b^{j+1}, b^{m})\neq 0$,  hence $\phi(b^{m+1})\leq \phi(b^{j})<\phi(b^{j+1})\leq \phi(b^{m})+1$, $\phi(M')<\phi(b^{j+1})\leq \phi(b^{m})+1$.  We have also  $M\in \sigma^{ss}$ and by $\hom^1(M',M)\neq 0$ and  $\hom(M,b^{m})\neq 0$  it follows that    $\phi(M')\leq \phi(M)+1$ and  $\phi(M)\leq b^{m}$.   These  arguments   together with \eqref{T43Zcap(E_1)1} imply  
\begin{gather} \label{T43Zcap(E_1)3} \begin{array}{c}  \phi(b^{m})-1<\phi(b^{m+1}[-1])< \phi(b^{m})\\ \phi(b^{m})-1<\phi(M'[-1])\leq \phi(M) \leq  \phi(b^{m}) ; \ \ \ \ \ \phi(M'[-1]) <  \phi(b^{m})  \end{array}. \end{gather}

  Due to \eqref{Z(delta)}, to show the first system in \eqref{T43Zcap(E_1)ineq} it remains to derive the following  inequality:
\begin{gather} \label{T43Zcap(E_1)4} \arg_{(\phi(b^{m})-1,\phi(b^{m}))}(Z(\delta))> \phi(M')-1.    \end{gather}   
From \eqref{T43Zcap(E_1)3} we see that $Z(\delta)$ and  $Z(M'[-1])$ both lie  in the half-plane\footnote{The notation  $Z(b^{m})^c_-$ is explained in \eqref{complement of a line}.}  $Z(b^{m})^c_-$. In \eqref{[delta]}    we have  aslo $Z(M)=Z(\delta) + Z(M'[-1])$, therefore the  vector $Z(M)$   is  in $Z(b^{m})^c_-$ as well. Now the equality  $Z(M)=Z(\delta) + Z(M'[-1])$ implies  that    \eqref{T43Zcap(E_1)4}  is equivalent to  $\phi(M)>\phi(M'[-1]) $. Hence we  have to show that   $\phi(M) \neq \phi(M'[-1]) $. Indeed, on one hand $\phi(M'[-1]) = \phi(M)$ implies $ \arg_{(\phi(b^{m})-1,\phi(b^{m}))}(Z(\delta))= \phi(M)$. On the other hand,
   $\sigma \in  (M, b^{j}, b^{j+1} )$,   $m<j$ and   \eqref{T43Zcap(E_1)1} imply    $\phi(M)+1<  \phi(b^{j+1})\leq \phi(b^{m})+1 < \phi(M')+1\leq\phi(M)+2$.  Thus, we see that   $\phi(M'[-1]) = \phi(M)$ implies  $Z(b^{j+1}) \in Z(\delta)^c_-$.  However, from the first inequality in  \eqref{T43Zcap(E_1)3}  and Corollary \ref{from equal phases to} (d)  it follows that $Z(b^{j+1}) \in Z(\delta)^c_+$, which is a contradiction,  and \eqref{T43Zcap(E_1)4} follows.

So far we showed the inclusion $\subset$. We show now the inverse inclusion $\supset$. 

We assume first that the second system of inequalities in \eqref{T43Zcap(E_1)ineq}  holds. In particular $\sigma \in (b^{m}, b^{m+1}, M')$. By the inequality $\phi(M')<\phi(b^{m+1}) $ we can apply  Proposition  \ref{phi_1>phi_2} (b), hence the short exact sequence \eqref{short filtration 1}  implies that $a^{m} \in \sigma^{ss}$ and $Z(M')+Z(b^{m+1})=Z(a^{m})$. We have also $\phi(b^{m+1})-1<\phi(M')<\phi(b^{m+1})$, and it follows that  $\phi(M')<\phi(a^{m})<\phi(b^{m+1})$. If  the inequality  $\phi(b^{m+1})>  \phi(b^{m})+1$ holds, then recalling that   $\phi(b^{m})<\phi(M')$ we obtain that  $\sigma \in (b^{m},a^{m},b^{m+1})\subset W$ (see table \eqref{no M}). 

 Therefore we reduce to the inequality  $\phi(b^{m+1}[-1])\leq  \phi(b^{m})$. Combining with    $\phi(M'[-1])<\phi(b^{m+1}[-1])<\phi(M')$ and  $\phi(b^{m})<\phi(M')$,      we can write   $\phi(M'[-1])< \phi(b^{m+1}[-1])\leq \phi(b^{m})<\phi(M')$  and then $(b^{m}, b^{m+1}[-1], M'[-1])$  is a $\sigma$-triple (see Definition \ref{sigma triple}).
  Combining  the triangles  \eqref{short filtration 1} and \eqref{short filtration 2} we obtain the following sequence of triangles in $\mc T$:    \be \label{T43Zcap(E_1)the filtation} 
\begin{diagram}[size=1em] 
0 & \rTo      &        &       &   M'[-1] & \rTo      &             &        & a^{m}[-1] &   \rTo      &           &       & M \\
  & \luDashto &        & \ldTo &         &\luDashto  &             &  \ldTo &           &   \luDashto &           & \ldTo &       \\
  &           &  M'[-1] &       &         &           &  b^{m+1}[-1]  &        &           &             & b^{m} &
\end{diagram}.
\ee
The conditions of Lemma \ref{three comp factors} (b) are satisfied with the triple  $(b^{m}, b^{m+1}[-1], M'[-1])$ and the diagram above. Therefore $M\in \sigma^{ss}$ and  $\phi(M)<\phi(b^{m})$. 

If   $\phi(b^{m+1}[-1])=\phi(b^{m})$, then we have also $\phi(M)+1<\phi(b^{m+1})$, and recalling that we have also $\phi(b^{m})<\phi(b^{m+1})$ we see that $\sigma \in  (M, b^{m}, b^{m+1} )\subset W$ (see table \eqref{left right M}).

  Therefore we can assume that   $\phi(M'[-1])< \phi(b^{m+1}[-1]) <  \phi(b^{m})<\phi(M')$. We will that $\sigma \in (M,b^{j}, b^{j+1})$ for some big $j$ in this case.  Proposition \ref{all exc semistable} and  Remark \ref{ext closure ab} ensure  that $\{b^{j+1}\}_{j\in \ZZ}\subset \sigma^{ss}$.  From  $\phi(b^{m})<\phi(b^{m+1}) <  \phi(b^{m})+1 $  and Corollary \ref{from equal phases to} (b) and (d) we see that $\phi(b^{j})<\phi(b^{j+1})$ and $Z(b^j)\in Z(\delta)_+^c$ for each $j$. 
		Now to show that  $\sigma \in (M,b^{j}, b^{j+1})\subset W$  it is enough to derive  $\phi(M)+1< \phi(b^j) $  for big enough $j$ (see  table \eqref{left right M}).
		
		Since we have  $ \phi(b^{m}[-1])<\phi(M'[-1])< \phi(b^{m+1}[-1]) <  \phi(b^{m})$ and $Z(\delta)=Z(b^{m})+Z(b^{m+1}[-1])$, we see that we can choose $t\in \RR$ so that $t< \phi(b^{m})<\phi(M')< \phi(b^{m+1})<t+1 $ and $Z(\delta)=\abs{Z(\delta)} \exp(\ri \pi t)$.    Since $\hom^1(b^j,b^{m})\neq 0$, $\hom(b^{m},b^j)\neq 0$   for $j>m+1$ and  by Corollary \ref{from equal phases to} (b), we have  $\phi(b^{m})<\phi(b^j)<  \phi(b^{m})+1 $ for  $j>m+1$. These inequalities together with the incidence $Z(b^j)\in Z(\delta)_+^c$ imply that $\arg_{(t,t+1)}(Z(b^j))=\phi(b^j)$  for $j>m+1$ (see Remark \ref{arg remark}  (c)). The formula \eqref{noncolinear ab2} in Corollary \ref{noncolinear ab}  gives us the following:
		\begin{gather}\label{T43Zcap(E_1)the limit} \lim_{j\rightarrow  \infty} \phi(b^j)= t+1.\end{gather}
We showed that  $\phi(M)<\phi(b^{m})$ (see below \eqref{T43Zcap(E_1)the filtation}) and we have also $\phi(b^{m})<\phi(M')$. From $\hom^1(M',M)\neq 0$ we derive  $\phi(b^{m}[-1])<\phi(M)<\phi(b^{m})$. We showed also that   $t< \phi(b^{m})<\phi(M')< \phi(b^{m+1})<t+1 $. Since $Z(M)+Z(M')=Z(\delta)=\abs{Z(\delta)} \exp(\ri \pi t)$, it follows that $Z(M)\in  Z(\delta)^c_-$ and $\phi(M)<t$. Now \eqref{T43Zcap(E_1)the limit} ensures that  $\phi(b^j)>\phi(M)+1$ for big enough $j$.  
So far  we showed that the second system of inequalities in \eqref{T43Zcap(E_1)ineq} implies  $\sigma\in (b^{m}, b^{m+1}, M') \cap W$.

We pass to  the  first system of inequalities in \eqref{T43Zcap(E_1)ineq}.  
So assume that $b^{m}, b^{m+1}, M' \in \sigma^{ss}$ and that these inequalities hold. They contain  the inequalities defining $(b^{m}, b^{m+1}, M')$ (see table \eqref{left right M}),  hence $\sigma \in (b^{m}, b^{m+1}, M')$.   Furthermore,  the first two inequalities show that  $(b^{m}, b^{m+1}[-1], M'[-1])$ is a $\sigma$-triple and that the conditions of Lemma \ref{three comp factors} (a) are satisfied with this triple and the diagram \eqref{T43Zcap(E_1)the filtation}. Therefore $M\in \sigma^{ss}$ and $\phi(M)<\phi(b^{m})$. By $\hom^1(M',M)\neq 0$ we can write 
$\phi(b^{m})-1< \phi(M'[-1])\leq \phi(M) < \phi(b^{m})$ (we use also \eqref{T43Zcap(E_1)ineq}), hence $Z(\delta), Z(M'[-1]), Z(M) \in Z(b^{m})^c_-$. Let us denote $t=\arg_{(\phi(b^{m})-1,\phi(b^{m}))}(Z(\delta))$. The third inequality in \eqref{T43Zcap(E_1)ineq} is the same as $t> \phi(M')-1 $. Combining  these arguments  with the equality   $Z(M)=Z(\delta)+Z(M'[-1])$ we deduce that: 
\begin{gather} \label{T43Zcap(E_1) help} \phi(b^{m}[-1])< \phi(M'[-1])< \phi(M) <t< \phi(b^{m}).  \end{gather}
We will show that  $\sigma \in  (M, b^{j}, b^{j+1} )$  for some big enough $j$.  We have  $\phi(b^{m})<\phi(b^{m+1}) <  \phi(b^{m})+1 $ (the first inequality in \eqref{T43Zcap(E_1)ineq}), which by     Proposition \ref{all exc semistable} and  Remark \ref{ext closure ab} implies  that $\{b^{j+1}\}_{j\in \ZZ}\subset \sigma^{ss}$,   and by Corollary \ref{from equal phases to} (b), (d) implies  that $\phi(b^{j})<\phi(b^{j+1})$ and $Z(b^j)\in Z(\delta)_+^c$ for each $j$.   Using Remark \ref{inc dec seq} one easily shows that   $\phi(b^{m})<\phi(b^j)<  \phi(b^{m})+1 $ for  $j>m+1$. These inequalities together with   the incidences $Z(b^j)\in Z(\delta)_+^c$,  $\phi(b^m)\in(t,t+1)$  imply that $\arg_{(t,t+1)}(Z(b^j))=\phi(b^j)$  for $j>m+1$ (see Remark \ref{arg remark}  (c)). The formula \eqref{noncolinear ab2} in Corollary \ref{noncolinear ab}  leads to  \eqref{T43Zcap(E_1)the limit} again.   Now \eqref{T43Zcap(E_1) help} implies  that $\phi(b^j)>\phi(M)+1$ for big  $j$, hence  $\sigma \in  (M, b^{j}, b^{j+1} )\subset W$  (see table \eqref{left right M}).

The arguments showing that  $ (b^{m}, b^{m+1}, M') \cap W$ and  $ (b^{m}, b^{m+1}, M') \cup W$  are contractible  are as in the last paragraph of the proof of Lemma \ref{T12Zcap(E_1)}. 
 The lemma is proved. 
\epr 

\begin{coro} \label{T43Zcap union(E_1)} The set $\mk{T}_{b}^{st}$ is contractible.  
\end{coro}
\bpr Recall that   $\mk{T}_{b}^{st}= W \cup \bigcup_{j\in \NN}  (b^{j}, b^{j+1}, M') $  (see \eqref{T43=Zcup}). Using  Lemmas \ref{T43Zcap(E_1)} and \ref{(_,_,X)0} one shows by induction that   $W \cup \bigcup_{j=0}^n  (b^{m-j+1}, b^{m-j}, M') $ is contractible  for each $m\in \ZZ$ and each $n\in \NN$ (see the proof of Corollary \ref{T12Zcap union(E_1)} for details). Then the corollary follows from Remark \ref{VK}. 
\epr

\section{Connecting \texorpdfstring{$\mk{T}_{a}^{st}$}{\space} and  \texorpdfstring{$ \mk{T}_{b}^{st}$}{\space} by \texorpdfstring{$(\_,M,\_)$}{\space} and \texorpdfstring{$(\_,M',\_)$}{\space}} \label{connecting}

   Due to the union  \eqref{st with T}, to prove Theorem \ref{main theo} it remains to connect  the contractible non-intersecting pieces  $\mk{T}_{a}^{st}$, $\mk{T}_{b}^{st}$ by $(\_,M,\_)$ and $(\_,M',\_)$, and to show that in this procedure the contractibility is preserved.  We describe first the building blocks of  $(\_,M,\_)$ and $(\_,M',\_)$ by Proposition \ref{lemma for f_E(Theta_E)}:

From the list of triples $\mk{T}$ given in Corollary \ref{exceptional colleections} we see that (see also \eqref{(A,_,_)}):\begin{gather} \label{middle M ab}  (\_,M\_) =\bigcup_{q\in \ZZ}(a^q,M,b^{q+1}) \qquad  
 (\_,M'\_)=\bigcup_{q\in \ZZ}(b^q,M',a^{q}) . \end{gather}

We apply Proposition \ref{lemma for f_E(Theta_E)} to the triples  $(a^p,M,b^{p+1})$ and $(b^q,M',a^{q})$. Using Corollaries \ref{nonvanishings} and \ref{one nonvanishing degree}   one shows that  the coefficients $\alpha, \beta, \gamma$ defined in \eqref{alpha,beta,gamma} are $\alpha=0$, $\beta=\gamma=-1$ in both the cases. Thus   we obtain  the   formulas  in the first and the second column of table \eqref{middle M} for the  contractible subsets $(a^p,M,b^{p+1}) \subset \st(D^b(\mc T))$   and  $(b^q,M',a^{q})\subset \st(D^b(\mc T))$, respectively:
\begin{gather} \label{middle M} \begin{array}{| c   | c  |}   \hline  
(a^p,M,b^{p+1})   & (b^q,M',a^{q})  \\
  \hline
 \left \{a^p,M,b^{p+1}  \in \sigma^{ss} :                      \begin{array}{c} \phi\left ( a^{p}\right) < \phi\left (M\right)+1 \\ \phi\left ( a^{p}\right) < \phi\left (b^{p+1}\right) \\ \phi\left ( M\right) < \phi\left (b^{p+1}\right)\end{array}  \right \}                            &       \left \{b^q,M',a^{q} \in \sigma^{ss} :                      \begin{array}{c} \phi\left ( b^{q}\right) < \phi\left (M'\right) +1\\ \phi\left ( b^{q}\right) < \phi\left (a^{q}\right) \\ \phi\left ( M'\right) < \phi\left (a^{q}\right)\end{array}  \right \}   \\  \hline
	\end{array}.
\end{gather}
\begin{remark}  \label{middle M rem}   $(a^p,M,b^{p+1}[-1])$,  $(b^q,M',a^{q}[-1])$ are Ext-exceptional triples (satisfy {\rm (a)} in Def. \ref{sigma triple}). 
	\end{remark}
In some steps of this section, when we need to show that certain exceptional  objects are semi-stable, the tools in Section \ref{general remarks} are not efficient enough. For these cases we prove Lemmas \ref{nonsemistable a} and \ref{nonsemistable b} below.   The relation   $\bd  R & \rDotsto & (S,E)  \ed $  between a $\sigma$-regular object  $R$ and an exceptional pair generated by it (introduced in \cite{DK1}) is utilized in the proof of these lemmas.   
\begin{lemma} \label{nonsemistable a} Let $a^m \not \in \sigma^{ss}$ and $t=\phi_-(a^m)$, then  one of the following holds:

 {\rm \textbf{(a)}} $a^j\in \sigma^{ss}$   for some $j<m-1$  and  $t=\phi(a^j)+1$;  {\rm \textbf{(b)}}  $a^j\in \sigma^{ss}$  for some  $m<j$  and     $t=\phi(a^j)$; 

{\rm \textbf{(c)}}  $b^j\in \sigma^{ss}$  for some $j<m$ and $t=\phi(b^j)+1$; 
 {\rm \textbf{(d)}}  $b^j\in \sigma^{ss}$  for some $m<j$ and $t=\phi(b^j)$;

 {\rm \textbf{(e)}} $M\in \sigma^{ss}$ and  $t=\phi(M)+1$. 
\end{lemma}
\bpr  Recall that any $X\in \{E_i^j: j\in \NN, 1\leq i \leq 4 \}$ is a trivially coupling object (see \cite[after Lemma 10.28]{DK1}). Since $a^m[k]\in  \{E_i^j: j\in \NN, 1\leq i \leq 4 \}$, where $k\in \{ 0,-1 \}$,  from $a^m \not \in \sigma^{ss}$ and \cite[Lemma 6.3]{DK1} it follows that $a^m[k]$ is a $\sigma$-regular object, hence $a^m$ is a $\sigma$-regular object.  Therefore  we have  $\begin{diagram} R & \rDotsto & (S,E) \end{diagram}$ for some exceptional pair $(S,E)$ (see \cite[Section 5]{DK1}).  We will need the following two properties of the exceptional object $S$. The first is $S\in \sigma^{ss}$,  $\phi(S)=\phi_-(a^m)$ (see \cite[formula (42) after Definition 5.2]{DK1}).  The second property is $\hom(a^m,S)\neq 0$, which follows from \cite[(c) after formula (19)]{DK1} and the way $S$ was chosen (see \cite[Definition 5.2]{DK1}).  Recall  that there exists at most one nonzero element in the family $\{ \hom^k(a^m,X) \}_{k\in \ZZ}$ for any   $X\in \mc T_{exc}$  (Corollary \ref{one nonvanishing degree}). By Remark \ref{exceptional objects} we have   $S\in \{a^j[k], b^j[k]: j\in \ZZ, k\in \ZZ\}\cup\{M[k],M'[k] ; k \in \ZZ\}$.
   Obviously $S\not = a^m[k]$ (since $a^m\not \in \sigma^{ss}$ and $S \in \sigma^{ss}$). 
	
	Now we will use the property    $\hom(a^m,S)\neq 0$ and Corollary \ref{nonvanishings} to prove the lemma.  By  $\hom^*(a^m, M')= 0$ (see  \eqref{nonvanishing1}) we exclude also  the case $S=M'[k]$.    It remains to consider the following cases (one of them must appear):

If $S=a^j[k]$ for some $j \neq m$ and $k\in \ZZ$, then by \eqref{nonvanishing5} we see that either   $j<m-1$ and  $k=1$, or    $m<j$ and  $k=0$. 

If $S=b^j[k]$ for some $j \in \ZZ$ and $k\in \ZZ$, then by \eqref{nonvanishing4}  we see that either    $j<m$ and  $k=1$, or   $m<j$ and  $k=0$. 

If $S=M[k]$ for some  $k\in \ZZ$, then by \eqref{nonvanishing2}  we get $k=1$. 
The lemma follows.
\epr

\begin{lemma} \label{nonsemistable b} Let $b^m \not \in \sigma^{ss}$ and $t=\phi_-(b^m)$, then  one of the following holds:

 {\rm \textbf{(a)}} $a^j\in \sigma^{ss}$   for some $j<m-1$  and  $t=\phi(a^j)+1$;  {\rm \textbf{(b)}}  $a^j\in \sigma^{ss}$  for some  $m\leq j$  and     $t=\phi(a^j)$; 

{\rm \textbf{(c)}}  $b^j\in \sigma^{ss}$  for some $j<m-1$ and $t=\phi(b^j)+1$; 
 {\rm \textbf{(d)}}  $b^j\in \sigma^{ss}$  for some $m<j$ and $t=\phi(b^j)$;

 {\rm \textbf{(e)}} $M'\in \sigma^{ss}$ and  $t=\phi(M')+1$. 
\end{lemma}
\bpr By the same arguments as in the proof of  Lemma \ref{nonsemistable a}  one shows that  $\hom(b^m,S)\neq 0$ and $\phi(S)=t$ for some  $S \in \sigma^{ss}\cap \left ( \{a^j[k], M, M': j\in \ZZ, k\in \ZZ\}\cup\{b^j[k] ; k \in \ZZ, j\in \ZZ,j\neq m\}\right )$.     Now we  will use  Corollaries \ref{nonvanishings} and \ref{one nonvanishing degree}.  
 By  $\hom^*(b^m, M)=0$ (see  \eqref{nonvanishing2}) we exclude   the case $S=M[k]$.    It remains to consider the following cases (one of them must appear):

If $S=a^j[k]$ for some $j \in \ZZ$ and $k\in \ZZ$, then by \eqref{nonvanishing3}  we see that either    $j<m-1$ and  $k=1$, or   $m\leq j$ and  $k=0$. 

If $S=b^j[k]$ for some $j \neq m$ and $k\in \ZZ$, then by \eqref{nonvanishing6} we see that either   $j<m-1$ and  $k=1$, or    $m<j$ and  $k=0$.

If $S=M'[k]$ for some  $k\in \ZZ$, then by \eqref{nonvanishing2}  we get $k=1$. 
The lemma follows.
\epr
Lemmas \ref{semi-stability of a} and \ref{semi-stability of b}  put together  the arguments which ensure semi-stability, necessary later  in the  analysis of   the intersections  $(a^p,M,b^{p+1}) \cap \mk{T}_{a/b}^{st}$  and  $(a^p,M,b^{p+1}) \cap (\mk{T}_{a}^{st} \cup (a^q,M,b^{q+1})\cup \mk{T}_{b}^{st})$. 

\begin{lemma} \label{semi-stability of a} Let $\sigma \in (a^p,M,b^{p+1})$ and let the following inequality hold:  \begin{gather} \label{semi-stability of a ineq}  \begin{array}{c}  \phi\left (b^{p+1}\right)-1< \phi\left ( M\right)<\phi\left (b^{p+1}\right)\end{array}.   \end{gather}
 Then we have the following:

{\rm \textbf{(a)}} $a^{p+1}\in \sigma^{ss}$ and   $\phi(b^{p+1})-1<\phi(a^{p+1})-1<\phi(M)$. 

{\rm \textbf{(b)}} If in addition to  \eqref{semi-stability of a ineq} we have $\phi(a^p)<\phi(M)$, then $\sigma \in (a^p,a^{p+1},M)$.

{\rm \textbf{ (c)}} If in addition to \eqref{semi-stability of a ineq} we have  \begin{gather} \label{semi-stability of a ineq1} \phi\left (b^{p+1}\right)-1< \phi\left ( a^{p}\right)<\phi\left (b^{p+1}\right), \end{gather}  then  $M'\in \sigma^{ss}$ and  $\phi(b^{p+1})-1<\phi(M')=\arg_{(\phi\left (b^{p+1}\right)-1,\phi\left (b^{p+1}\right))}(Z(a^p)-Z(b^{p+1}))<\phi(a^{p})$.

{\rm \textbf{ (d)}} If    \eqref{semi-stability of a ineq}, \eqref{semi-stability of a ineq1} hold and     $\phi(M')<\phi(M)$, then $\sigma \in (a^j,a^{j+1},M)$ for some $j\in \ZZ$.
\end{lemma}
\bpr  (a)  We  apply Proposition \ref{phi_1>phi_2} (b) to the triple $(a^p,M,b^{p+1})$ and since $a^{p+1}[-1]$ is in the extension closure of $M, b^{p+1}[-1]$ (by  \eqref{short filtration 2}) it follows that  $a^{p+1}\in \sigma^{ss}$.  The inequality  $\phi(b^{p+1})-1<\phi(a^{p+1})-1<\phi(M)$ follows from the given inequality \eqref{semi-stability of a ineq} and $Z(a^{p+1}[-1])=Z(M)+Z(b^{p+1}[-1])$. 

(b) From the given inequalities it follows  that $\phi(a^p)<\phi(b^{p+1})$. We have also  $\phi(b^{p+1})<\phi(a^{p+1})<\phi(M)+1$ from (a). Therefore we obtain the inequalities  $\phi(a^{p})<\phi(a^{p+1})$, $\phi(a^p)<\phi(M)$,  $\phi(a^{p+1})<\phi(M)+1$, which means that $\sigma \in (a^p,a^{p+1},M)$ (see table \eqref{left right M}). 

    (c)   Follows from  Lemma \ref{two comp factors} applied  to the Ext-triple $(a^p,M,b^{p+1}[-1])$ and the triangle \eqref{short filtration 1}.

 (d) Now by the given inequalities and (c) we have $\phi(b^{p+1})-1<\phi(M')<\phi(M)<\phi(b^{p+1})$. Recalling that  $Z(\delta)=Z(M')+Z(M)$, we see that we can choose  $t\in \RR$ with   $Z(\delta)=\abs{Z(\delta)}\exp(\ri \pi t)$ and $\phi(M')<t<\phi(M)<\phi(b^{p+1})<t+1$.  If $\phi(a^p)<\phi(M)$, then (d) follows from (b). 

So let $\phi(M)\leq \phi(a^p)$. Since we have also $\phi(a^p)<\phi(b^{p+1})$, we obtain  $t<\phi(M)\leq \phi(a^p)<\phi(b^{p+1})<t+1$. Now Corollary  \ref{noncolinear ab}   shows that  $\{Z(a^j), Z(b^j)\}_{j\in \ZZ}\subset Z(\delta)^c_+$ and  that  \eqref{noncolinear ab2}, \eqref{noncolinear ab1} hold for both the sequences  $\{Z(a^j)\}_{j\in \ZZ}$ and $\{Z(b^j)\}_{j\in \ZZ}$.  From (a) we see that $\phi(b^{p+1})<\phi(a^{p+1})<\phi(M)+1$, hence $t< \phi(a^{p+1})<t+2$, which combined with $Z(a^{p+1})\in Z(\delta)^c_+$ implies that $\phi(a^{p+1})<t+1$. Thus we obtain the inequalities  $t<\phi(M)\leq \phi(a^p)<\phi(b^{p+1})<\phi(a^{p+1})<t+1$. 

From \eqref{noncolinear ab2} we see that there exists  $N\in \ZZ$, $N<p$ such that  $t<\arg_{(t,t+1)}(Z(a^j))<\phi(M)$ for  $j<N$. We will show below  that  $a^j\in \sigma^{ss}$ for $j<N$. Then (d) follows. Indeed,  assume that   $a^j\in \sigma^{ss}$ for each $j<N$.  Then by \eqref{nonvanishing5} and Corollary \ref{from equal phases to}  (a) it follows that  $\phi(a^{p+1})-1<\phi(a^j) <\phi(a^{p+1})$ for  $j<N$, therefore  $t-1<\phi(a^j) <t+1$, which combined with $Z(a^j)\in Z(\delta)^c_+$ implies that  $\arg_{(t,t+1)}(Z(a^j))=\phi(a^j)$. Putting the last equality  in  \eqref{noncolinear ab1}  and in  $\arg_{(t,t+1)}(Z(a^j))<\phi(M)$  we get $\phi(a^{j-1})<\phi(a^j)<\phi(M)$    which by table \eqref{left right M} implies that $\sigma \in (a^{j-1},a^j,  M)$. 

 Suppose $a^j \not \in \sigma^{ss}$ for some $j<N$.  From Remark \ref{ext closure ab} we know  that $a^j$ is in the extension closure of $a^p$, $a^{p+1}[-1]$. It follows that  $a^j\in\mc P[\phi(a^{p+1})-1,\phi(a^p)]$ and then $ \phi(a^{p+1})-1\leq \phi_-(a^j)$ (recall the paragraph after \eqref{sigma^{ss}}). We will use Lemma \ref{nonsemistable a} and show that each of the five cases given there leads to a contradiction. We fist derive  \eqref{semi-stability of a ineq2}.  The inequalities  $ \phi(a^p)-1<\phi(a^{p+1})-1<\phi(M) 
\leq \phi(a^p)$ can be used due to the  previous steps.   Therefore we have $a^j\in\mc P[\phi(a^{p+1})-1,\phi(a^p)] \subset  \mc P(\phi(a^{p})-1,\phi(a^p)]$. Using $\phi(a^p)\in (t,t+1)$,  $Z(a^j)\in Z(\delta)^c_+$, and Remark \ref{arg remark}  (c)   we get: $\arg_{(\phi(a^{p})-1,\phi(a^p)]}(Z(a^j))=\arg_{(t,t+1)}(Z(a^j))$.  Now by  Remark \ref{arg remark}  (a) we get  $\phi_-(a^j) < \arg_{(t,t+1)}(Z(a^j))$ and by our choice of $N$ we have  $\arg_{(t,t+1)}(Z(a^j))<\phi(M)$. We combine these facts in the following inequalities:  
\begin{gather} \label{semi-stability of a ineq2} \phi(a^{p+1})-1\leq \phi_-(a^j) < \arg_{(t,t+1)}(Z(a^j))<\phi(M)\leq \phi(a^p)<\phi(a^{p+1}). \end{gather}

 One of the cases in Lemma \ref{nonsemistable a} must appear.   In case (a)   we have   $\phi_-(a^j) =\phi(a^k)+1 $ for some $k<j-1$, hence by   \eqref{semi-stability of a ineq2} it follows $\hom^1(a^p,a^k)=0$, which contradicts    \eqref{nonvanishing5} and  $j<N<p$. 

In case (b):    $\phi_-(a^j) =\phi(a^k) $ for some $k>j$. It follows that $\phi(a^k)=\arg_{(t,t+1)}(Z(a^k))$  (see Remark \ref{arg remark} (c)), hence   by \eqref{semi-stability of a ineq2} and \eqref{phase formula} we get $\arg_{(t,t+1)}(Z(a^k))<\arg_{(t,t+1)}(Z(a^j))$, which contradicts  \eqref{noncolinear ab1}.

In cases (c) and (d) we have $\phi_-(a^j) =\phi(b^k) $ of $\phi(b^k)+1$ for some $k\in \ZZ$, and then  \eqref{semi-stability of a ineq2} implies $\hom(M,b^k)=0$, which contradicts \eqref{nonvanishing1}. 

Case  (e) in Lemma \ref{nonsemistable a} and \eqref{semi-stability of a ineq2} imply that $\phi(M)+1 <\phi(M)$ and  we proved the lemma. 
\epr

\begin{lemma} \label{semi-stability of b} Let $\sigma \in (a^p,M,b^{p+1})$ and let the following inequality hold:  \begin{gather} \label{semi-stability of b ineq}  \begin{array}{c}  \phi\left (a^{p}\right)-1< \phi\left ( M\right)<\phi\left (a^{p}\right)\end{array}.   \end{gather}
 Then we have the following:

{\rm \textbf{(a)}} $b^{p}\in \sigma^{ss}$ and   $\phi(M)<\phi(b^{p})<\phi(a^{p})$. 

{\rm \textbf{(b)}} If in addition to  \eqref{semi-stability of b ineq} we have $\phi(M)+1<\phi(b^{p+1})$, then $\sigma \in (M, b^p, b^{p+1})$.

{\rm \textbf{ (c)}} If in addition to \eqref{semi-stability of b ineq} we have  \begin{gather} \label{semi-stability of b ineq1}  \phi\left (a^{p}\right)-1< \phi\left ( b^{p+1}\right)-1<\phi\left (a^{p}\right), \end{gather}  then  $M'\in \sigma^{ss}$ and  $\phi(b^{p+1})-1<\phi(M')=\arg_{(\phi\left (a^{p}\right)-1,\phi\left (a^{p}\right))}(Z(a^p)-Z(b^{p+1}))<\phi(a^{p})$.

{\rm \textbf{ (d)}} If    \eqref{semi-stability of b ineq}, \eqref{semi-stability of b ineq1} hold and     $\phi(M)<\phi(M')$, then $\sigma \in (b^j,b^{j+1},M')$ for some $j\in \ZZ$ or $\sigma \in  (M, b^p, b^{p+1})$.

{\rm \textbf{ (e)}} If    \eqref{semi-stability of b ineq}, \eqref{semi-stability of b ineq1} hold and     $\phi(M)=\phi(M')$, then $\sigma \in (a^j,M,b^{j+1})$ for each $j<p$.
\end{lemma}
\bpr  \textbf{(a)}  We  apply Proposition \ref{phi_1>phi_2} (a) to the triple $(a^p,M,b^{p+1})$ and since $b^{p}$ is in the extension closure of $M, a^p$ (by  \eqref{short filtration 2}) it follows that  $b^{p}\in \sigma^{ss}$ and $\phi(M)\leq \phi(b^{p})\leq \phi(a^{p})$.  The inequality  $\phi(M)<\phi(b^{p})<\phi(a^{p})$ follows from the given inequality \eqref{semi-stability of b ineq} and $Z(b^{p})=Z(M)+Z(a^{p})$. 

\textbf{(b)} From the given inequalities we have $\phi(a^p)<\phi(M)+1<\phi(b^{p+1})$. In (a) we showed that   $\phi(M)<\phi(b^{p})<\phi(a^{p})$. Therefore we obtain the inequalities  $\phi(M)<\phi(b^{p})$, $\phi(M)+1<\phi(b^{p+1})$,  $\phi(b^{p})<\phi(b^{p+1})$, which means that $\sigma \in (M,b^p,b^{p+1})$ (see table \eqref{left right M}). 

    \textbf{(c)}   Follows from  Lemma \ref{two comp factors} applied  to the Ext-triple $(a^p,M,b^{p+1}[-1])$ and the triangle \eqref{short filtration 1}.

\textbf{(d)} Now by the given inequalities and (c) we have $\phi(a^{p})-1<\phi(M)<\phi(M')<\phi(a^{p})$. Recalling that  $Z(\delta)=Z(M')+Z(M)$, we see that we can choose  $t\in \RR$ with   $Z(\delta)=\abs{Z(\delta)}\exp(\ri \pi t)$ and $\phi(M)<t<\phi(M')<\phi(a^{p})<\phi(M)+1$.  If $\phi(M)+1<\phi(b^{p+1})$, then we apply (b). 

So, let $\phi(b^{p+1})\leq\phi(M)+1$. Since we have also $\phi(a^p)<\phi(b^{p+1})$, we obtain  $t<\phi(M')<\phi(a^p)<\phi(b^{p+1})\leq\phi(M)+1 <t+1$. Now Corollary  \ref{noncolinear ab}   ensures that  $\{Z(a^j), Z(b^j)\}_{j\in \ZZ}\subset Z(\delta)^c_+$ and  that \eqref{noncolinear ab2}, \eqref{noncolinear ab1} hold for both the sequences  $\{Z(a^j)\}_{j\in \ZZ}$ and $\{Z(b^j)\}_{j\in \ZZ}$.  From (a) we see that $\phi(M)<\phi(b^{p})<\phi(a^{p})$, hence $t-1< \phi(b^{p})<t+1$, which combined with $Z(b^{p})\in Z(\delta)^c_+$ implies that $t<\phi(b^{p})$.  Hence we obtain the inequalities  \begin{gather} \label{help ineq} t<\phi(b^p)< \phi(a^p)<\phi(b^{p+1})<t+1; \quad t<\phi(M')<\phi(a^p)<\phi(b^{p+1}) <t+1 .  \end{gather}

From \eqref{noncolinear ab2} and $t<\phi(M')$ it follows that there exists  $N\in \ZZ$, $N<p$ such that  $t<\arg_{(t,t+1)}(Z(b^j))<\phi(M')$ for  $j<N$. We will show below  that  $b^j\in \sigma^{ss}$ for $j<N$. Then (d) follows. Indeed,  assume that  $b^j\in \sigma^{ss}$ for  each $j<N$.  Then by \eqref{nonvanishing6} and Corollary \ref{from equal phases to} (a) it follows that  $\phi(b^{p+1})-1<\phi(b^j) <\phi(b^{p+1})$ for  $j<N$, and by \eqref{help ineq} we get  $t-1<\phi(b^j) <t+1$, which combined with $Z(b^j)\in Z(\delta)^c_+$ implies that  $\arg_{(t,t+1)}(Z(b^j))=\phi(b^j)$. Putting the last equality in  \eqref{noncolinear ab1} and in  $\arg_{(t,t+1)}(Z(b^j))<\phi(M')$ we obtain $\phi(b^{j-1})<\phi(b^j)<\phi(M')$, which implies $\sigma \in (b^{j-1},b^j,  M')$. 

 Suppose $b^j \not \in \sigma^{ss}$ for some $j<N$.  We apply  Lemma \ref{nonsemistable b}  and show that each of the five cases given there leads to a contradiction. We show first \eqref{semi-stability of b ineq2}.  From Remark \ref{ext closure ab} we know  that $b^j$ is in the extension closure of $b^p$, $b^{p+1}[-1]$ (recall that  $N<p$) and   we  have  $ \phi(b^p)-1<\phi(b^{p+1})-1 <t
< \phi(b^p)$ in \eqref{help ineq}.   It follows that  $b^j\in\mc P[\phi(b^{p+1})-1,\phi(b^p)] \subset  \mc P(\phi(b^{p})-1,\phi(b^p)]$. Using $\phi(b^p)\in (t,t+1)$,  $Z(b^j)\in Z(\delta)^c_+$,  Remark \ref{arg remark} (c) and (a), we deduce that $\arg_{(\phi(b^{p})-1,\phi(b^p)]}(Z(b^j))=\arg_{(t,t+1)}(Z(b^j))>\phi_-(b^j)$. The incidence   $b^j\in\mc P[\phi(b^{p+1})-1,\phi(b^p)]$ implies  $\phi(b^{p+1})-1\leq \phi_-(b^j)$, and we get: 
\begin{gather} \label{semi-stability of b ineq2} \phi(b^{p+1})-1\leq \phi_-(b^j) < \arg_{(t,t+1)}(Z(b^j))<\phi(M') <\phi(b^{p+1}). \end{gather}

 One of the cases in Lemma \ref{nonsemistable b} must appear. In cases (a) and (b) we have $\phi_-(b^j) =\phi(a^k) $ of $\phi(a^k)+1$ for some $k\in \ZZ$, and then  \eqref{semi-stability of b ineq2} implies $\hom(M',a^k)=0$, which contradicts \eqref{nonvanishing1}. 

In case (c) we have   $\phi_-(b^j) =\phi(b^k)+1 $ for some $k<j-1$, and \eqref{semi-stability of b ineq2} implies that $\hom^1(b^{p+1},b^k)=0$, which contradicts \eqref{nonvanishing6} and $k<j-1<p-1$. 

In case (d) we have   $\phi_-(b^j) =\phi(b^k) $ for some $k>j$. From   $Z(b^k) \in Z(\delta)^c_+$,  \eqref{semi-stability of b ineq2}, and  $\phi(b^{p+1})\in (t,t+1)$  it follows that  $\phi(b^k)=\arg_{(t,t+1)}(Z(b^k))$. Hence \eqref{semi-stability of b ineq2} and \eqref{phase formula} imply  $\arg_{(t,t+1)}(Z(b^k))<\arg_{(t,t+1)}(Z(b^j))$, which contradicts  $k>j$ and  \eqref{noncolinear ab1}.

Case  (e) in Lemma \ref{nonsemistable b} and \eqref{semi-stability of b ineq2} imply that $\phi(M')+1 <\phi(M')$.  We  proved completely part  (d) of the lemma. 

\textbf{(e)}  Now by the given inequalities  we have $\phi(a^{p})-1<\phi(M)=\phi(M')<\phi(a^{p})$. Recalling that  $Z(\delta)=Z(M')+Z(M)$, we see that   $t=\phi(M)=\phi(M')$ satisfies    $Z(\delta)=\abs{Z(\delta)}\exp(\ri \pi t)$ and $t<\phi(a^{p})<t+1$.   From (a) we get $t<\phi(b^{p})<\phi(a^{p})<t+1$.   Now we can apply Corollary \ref{noncolinear ab full1}, which  besides   $\{Z(a^j), Z(b^j)\}_{j\in \ZZ}\subset Z(\delta)^c_+$ and  formulas    \eqref{noncolinear ab1}, \eqref{noncolinear ab2}  gives us the inequalities \eqref{j<mleqi}. 

 We extend  the inequality $t<\phi(b^{p})<\phi(a^{p})<t+1$  to \eqref{semi-stability of b ineq3} as follows.    We already have  that $a^p, b^p, b^{p+1}\in \sigma^{ss}$.   In \eqref{semi-stability of b ineq1} is given that  $\phi(a^{p})<\phi(b^{p+1})$. 
  From $\hom^1(b^{p+1},M')$ (see \eqref{nonvanishing2}) it follows $\phi(b^{p+1})\leq t+1$ and from $Z(b^{p+1})\in Z(\delta)^c_+$ we see that  $\phi(b^{p+1})< t+1=\phi(M)+1$.  We have also  $\phi(M)<\phi(b^{p+1})$ (due to $\sigma \in (a^p,M,b^{p+1})$). Therefore $\phi(b^{p+1})-1<\phi(M)< \phi(b^{p+1})$ and from Lemma \ref{semi-stability of a} (a) we get  $a^{p+1}\in \sigma^{ss}$ and   $\phi(b^{p+1})-1<\phi(a^{p+1})-1<\phi(M)$. Thus, we derive:
\begin{gather} \label{semi-stability of b ineq3} \phi(a^{p})-1<\phi(b^{p+1})-1<\phi(a^{p+1})-1<t<\phi(b^{p})<\phi(a^{p})<\phi(b^{p+1})<\phi(a^{p+1})<t+1.\end{gather}

We will show below that $a^j$ and $b^j$ are semi-stable  for each $j<p$. We claim that this implies    $\sigma \in (a^{j}, M, b^{j+1})$ for $j<p$. Indeed, assume that $a^j, b^j \in \sigma^{ss}$ for each  $j<p$.  Then by  \eqref{nonvanishing5}, \eqref{nonvanishing6} we get $\phi(a^{p+1})-1\leq \phi(a^j)\leq \phi(a^{p+1})$ and $\phi(b^{p+1})-1\leq \phi(b^j)\leq \phi(b^{p+1})$, which combined with  $t<\phi(b^{p+1})<\phi(a^{p+1})<t+1$ and $Z(a^j), Z(b^j)\in Z(\delta)^c_+$ implies that $\phi(a^j),\phi(b^j) \in (t,t+1)$, in particular   $\phi(a^j)=\arg_{(t,t+1)}(Z(a^j))$ and  $\phi(b^j)=\arg_{(t,t+1)}(Z(b^j))$ for each $j<p$. The last two equalities hold also for $j=p$ by \eqref{semi-stability of b ineq3}. Putting these equalities in  \eqref{j<mleqi} we get that $\phi(a^{j})<\phi(b^{j+1})$ for each $j<p$. Thus,  we obtain   $\phi(M)<\phi(a^{j})<\phi(b^{j+1})<\phi(M)+1$ for each $j<p$, which by table \eqref{middle M} gives  $\sigma \in (a^{j}, M, b^{j+1})$.

Suppose that $b^j\not \in \sigma^{ss}$ for some $j<p$.   Remark \ref{ext closure ab} asserts that  $b^j$ is in the extension closure of $b^p$, $b^{p+1}[-1]$, therefore $b^j\in\mc P[\phi(b^{p+1})-1,\phi(b^p)]$, and hence $\phi(b^{p+1})-1\leq \phi_-(b^j)$, $\phi_+(b^j)\leq \phi(b^p)$. Due to \eqref{semi-stability of b ineq3} we can write  $b^j\in\mc P[\phi(b^{p+1})-1,\phi(b^p)] \subset  \mc P(\phi(a^{p})-1,\phi(a^p)]$. Using $\phi(a^p)\in (t,t+1)$,  $Z(b^j)\in Z(\delta)^c_+$ and Remark \ref{arg remark} (c) we conclude that $\arg_{(\phi(a^{p})-1,\phi(a^p)]}(Z(b^j))=\arg_{(t,t+1)}(Z(b^j))$.  Now  using Remark \ref{arg remark}  (a),    we  obtain: 
\begin{gather} \label{semi-stability of b ineq4} \phi(b^{p+1})-1\leq \phi_-(b^j) < \arg_{(t,t+1)}(Z(b^j))<\phi_+(b^j)\leq  \phi(b^p) <\phi(b^{p+1}). \end{gather}
 We use Lemma \ref{nonsemistable b}  and  show that each of the five cases given there leads to a contradiction. 

Case (a) ensures $\phi_-(b^j) = \phi(a^k)+1$  for some $k<j-1$ and \eqref{semi-stability of b ineq4} implies that $\hom^1(b^{p+1},a^k)=0$, which contradicts \eqref{nonvanishing3} (now $k<p$).

Case  (b) ensures  $\phi_-(b^j) =\phi(a^k) $  for some $k\geq j$, and then  \eqref{semi-stability of b ineq4}  and  $Z(a^k) \in Z(\delta)^c_+$ 
 imply  $\arg_{(t,t+1)}(Z(a^k))=\phi(a^k)$, hence by  \eqref{semi-stability of b ineq4} and \eqref{phase formula} we get $\arg_{(t,t+1)}(Z(a^k))<\arg_{(t,t+1)}(Z(b^j))$, which contradicts \eqref{j<mleqi} and  $k\geq j$. 

Case (c) ensures   $\phi_-(b^j) =\phi(b^k)+1 $ for some $k<j-1<p-1$, and \eqref{semi-stability of b ineq4} implies that $\hom^1(b^{p+1},b^k)=0$, which contradicts \eqref{nonvanishing6}. 

In case (d) we have   $\phi_-(b^j) =\phi(b^k) $ for some $k>j$. It follows by $Z(b^k)\in Z(\delta)^c_+$ and \eqref{semi-stability of b ineq4} that  $\phi(b^k)=\arg_{(t,t+1)}(Z(b^k))$, and then \eqref{semi-stability of b ineq4} gives $\arg_{(t,t+1)}(Z(b^k))<\arg_{(t,t+1)}(Z(b^j))$, which contradicts   \eqref{noncolinear ab1}.

In case  (e)  using \eqref{semi-stability of b ineq4} we obtain $\phi(M')+1 < \phi(b^{p+1})$, which contradicts \eqref{nonvanishing2}.

Suppose that $a^j\not \in \sigma^{ss}$ for some $j<p$. Since  $a^j$ is in the extension closure of $a^p$, $a^{p+1}[-1]$ (see  Remark \ref{ext closure ab} ), therefore $a^j\in\mc P[\phi(a^{p+1})-1,\phi(a^p)]$, and hence $\phi_\pm(a^j) \in [\phi(a^{p+1})-1,\phi(a^p)]$. Due to \eqref{semi-stability of b ineq3} we have  $a^j\in\mc P[\phi(a^{p+1})-1,\phi(a^p)] \subset  \mc P(\phi(b^{p+1})-1,\phi(b^{p+1})]$ and  Remark \ref{arg remark} (c) shows that that $\arg_{(\phi(b^{p+1})-1,\phi(b^{p+1})]}(Z(a^j))=\arg_{(t,t+1)}(Z(a^j))$.  Now   Remark \ref{arg remark}  (a) completes the following: 
\begin{gather} \label{semi-stability of b ineq5} \phi(a^{p+1})-1\leq \phi_-(a^j) < \arg_{(t,t+1)}(Z(a^j))<\phi_+(a^j)\leq  \phi(a^p) <\phi(a^{p+1}). \end{gather}

 We use Lemma \ref{nonsemistable a}  to get a contradiction.  One of the five  cases given there must appear.

In case (a) of  Lemma \ref{nonsemistable a}   we have   $\phi_-(a^j) =\phi(a^k)+1 $ for some $k<j-1<p-1$, and \eqref{semi-stability of b ineq5} implies $\hom^1(a^{p+1},a^k)=0$, which contradicts \eqref{nonvanishing5}. 

Case (b)  ensures  $\phi_-(a^j) =\phi(a^k) $ for some $k>j$. It follows that $\phi(a^k)=\arg_{(t,t+1)}(Z(a^k))$  (see Remark \ref{arg remark} (c)), hence   by \eqref{semi-stability of b ineq5} we get $\arg_{(t,t+1)}(Z(a^k))<\arg_{(t,t+1)}(Z(a^j))$, which contradicts  \eqref{noncolinear ab1}.

In case (c) we have  $\phi_-(a^j) = \phi(b^k)+1$  for some $k<j$ and \eqref{semi-stability of b ineq5} implies that $\hom^1(a^{p+1},b^k)=0$, which contradicts \eqref{nonvanishing4} (now $k<p$).

Case  (d) ensures  $\phi_-(a^j) =\phi(b^k) $  for some $ j<k$, and then  $\arg_{(t,t+1)}(Z(b^k))=\phi(b^k)$ (see Remark \ref{arg remark} (c)), 
 hence by  \eqref{semi-stability of b ineq5} we get $\arg_{(t,t+1)}(Z(b^k))<\arg_{(t,t+1)}(Z(a^j))$, which contradicts \eqref{j<mleqi}. 

In case  (e) we have   $\phi_-(a^j) =\phi(M)+1$, and   \eqref{semi-stability of b ineq5} implies  $\hom^1(a^{p+1},M)=0$, which contradicts \eqref{nonvanishing2}.
The lemma is proved.
\epr
Next we glue $(a^p,M,b^{p+1})$ and  $\mk{T}_{a}^{st}$.
\begin{lemma} \label{middle M cap left M'} For any $p\in \ZZ$ the set  $(a^p,M,b^{p+1}) \cap \mk{T}_{a}^{st}$ consists of the stability conditions $\sigma$ for which  $a^p,M,b^{p+1}$ are semistable and:   \begin{gather} \label{middle M cap left M' sys}\begin{array}{c}  \phi\left (b^{p+1}\right)-1< \phi\left ( M\right)<\phi\left (b^{p+1}\right)\\  \phi\left (b^{p+1}\right)-1< \phi\left ( a^{p}\right)<\phi\left (b^{p+1}\right) \\ \arg_{(\phi\left (b^{p+1}\right)-1,\phi\left (b^{p+1}\right))}(Z(a^p)-Z(b^{p+1}))<\phi(M) \end{array}    \ \mbox{or} \  \begin{array}{c}  \phi\left (a^p\right)< \phi\left ( M\right)\\ \phi\left (b^{p+1}\right)-1< \phi\left ( M\right)<\phi\left (b^{p+1}\right)\end{array}.   \end{gather}
It follows that $(a^p,M,b^{p+1}) \cap \mk{T}_{a}^{st}$ and  $   (a^p,M,b^{p+1}) \cup \mk{T}_{a}^{st}$ are  contractible. 
\end{lemma}
\bpr We start with the inclusion $\subset$. Assume that $\sigma \in (a^p,M,b^{p+1})$. Then $a^p,M,b^{p+1}$ are semi-stable and by table \eqref{middle M} we get  \begin{gather} \label{middle M cap left M'1} \begin{array}{c} \phi\left ( a^{p}\right) < \phi\left (M\right)+1 \\ \phi\left ( a^{p}\right) < \phi\left (b^{p+1}\right) \\ \phi\left ( M\right) < \phi\left (b^{p+1}\right)\end{array} \end{gather}
Recalling \eqref{T_12new},  we see that we have to consider three cases.

\ul{If $\sigma \in (M',a^j, a^{j+1})$}, then $M',a^j, a^{j+1}$ are semi-stable and   from table \eqref{left right M} we see that $\phi(M')+1<\phi(a^{j+1})$.  Since we have also $\hom^1(b^{p+1},M')$, $\hom^1(a^{j+1},M)\neq 0$(see Corollary \ref{nonvanishings}),  we obtain $ \phi\left (b^{p+1}\right) \leq \phi(M')+1<\phi(a^{j+1})\leq \phi(M)+1$, which combined with \eqref{middle M cap left M'1} implies 
\begin{gather}\label{middle M cap left M'2} \phi\left (b^{p+1}\right)-1< \phi\left ( M\right)<\phi\left (b^{p+1}\right) \qquad \phi(M')< \phi(M). \end{gather} These non-vanishings and  inequalities give also  $\phi\left ( a^{p}\right) < \phi\left (b^{p+1}\right) \leq \phi(M')+1<\phi(a^{j+1})$.  Using Remark \ref{inc dec seq} (a)  we deduce that $p\leq j$. 

 We  show now that $\phi(b^{p+1})<\phi(a^p)+1$. If $j=p$, then we immediately obtain this  by $\hom^1(b^{p+1}, M')\neq 0$ and $ \phi(M')<\phi(a^p)$(see table \eqref{left right M}). If $j>p$, then $\hom(b^{p+1},a^j)\neq 0$ and $\hom^1(a^{j+1},a^p)\neq 0$ (see Corollary \ref{nonvanishings}) and we can write $\phi(b^{p+1})\leq\phi(a^{j})<\phi(a^{j+1})\leq \phi(a^{p})+1 $.

 To obtain the first system of inequalities in \eqref{middle M cap left M' sys} it remains to show the third inequality. 
From the triangle \eqref{short filtration 1}  it follows that $\phi(b^{p+1})-1\leq \phi(M')\leq \phi(a^p)$  and  $Z(M')=Z(a^p)-Z(b^{p+1})$, now  $\phi(M')=\arg_{(\phi\left (b^{p+1}\right)-1,\phi\left (b^{p+1}\right))}(Z(a^p)-Z(b^{p+1}))<\phi(M)$ follows from the already proved $\phi\left (b^{p+1}\right)-1< \phi\left ( a^{p}\right)<\phi\left (b^{p+1}\right)$ and \eqref{middle M cap left M'2}.

\ul{If $\sigma\in  (a^m,a^{m+1},M)$}, then $a^m,a^{m+1}$ are semistable as well and in table \eqref{left right M} we see that $\phi(a^m)<\phi(M)$, which together with the third inequality in \eqref{middle M cap left M'1} imply that $\phi(a^m)<\phi(b^{p+1})$ and hence $\hom(b^{p+1},a^m)=0$. By \eqref{nonvanishing3} we deduce that $p\geq m$. 

If $p=m$, then  we get immediately $\phi(a^p)<\phi(M)$. In table \eqref{left right M} we have $\phi(a^{p+1})<\phi(M)+1$ and in Corollary \ref{nonvanishings} we have $\hom(b^{p+1},a^{p+1})\neq 0$, hence $\phi(b^{p+1})<\phi(M)+1$ and we obtain the second system of inequalities in \eqref{middle M cap left M' sys}.

If $p>m$, then $\hom^1(b^{p+1},a^m)\neq 0$ and from the inequalities $\phi(a^m)<\phi(M)$, $\phi(a^m)<\phi(a^{m+1})$ (due to $\sigma \in (a^m,a^{m+1},M)$) it follows $\phi(b^{p+1})<\phi(M)+1$ and  $\phi(b^{p+1})\leq \phi(a^m)+1<\phi(a^{m+1})+1\leq \phi(a^{p})+1$. Recalling \eqref{middle M cap left M'1} we see that we obtained  the first two equalities in \eqref{middle M cap left M' sys}. Hence by Lemma \ref{semi-stability of a} (c) we get $M'\in \sigma^{ss}$ and $\phi(M')=\arg_{(\phi\left (b^{p+1}\right)-1,\phi\left (b^{p+1}\right))}(Z(a^p)-Z(b^{p+1}))$. From $\hom(M',a^m)\neq 0$ it follows $\phi(M')\leq \phi(a^m)<\phi(M)$ and we obtain the complete first system of inequalitites in \eqref{middle M cap left M' sys}. 

\ul{If $\sigma\in  (a^m, b^{m+1},a^{m+1})$}, then $a^m, b^{m+1},a^{m+1}\in \sigma^{ss}$ and in table \eqref{no M} we see that $\phi(a^m)+1<\phi(a^{m+1})$, hence  Lemma \ref{from phases of big distance} and $a^p \in \sigma^{ss}$ imply that $p=m$ or $p=m+1$. If $p=m+1$, then by \eqref{middle M cap left M'1} we obtain $\phi(a^m)+1<\phi(a^{m+1})<\phi(b^{m+2})$, and hence $\hom^1(b^{m+2},a^m)=0$, which contradicts \eqref{nonvanishing3}. Thus, it remains to consider the case $m=p$. Now we have $\phi(a^p)+1<\phi(a^{p+1})$ and $\phi(b^{p+1})<\phi(a^{p+1})$(see table \eqref{no M}), which  together with $\hom^1(a^{p+1},M)\neq 0$ imply $\phi(a^p)<\phi(M)$ and $\phi(b^{p+1})<\phi(M)+1$,  hence we obtain the second system  in \eqref{middle M cap left M' sys}. The inclusion $\subset$ is shown. 

 We show now the converse  inclusion $\supset$. Assume that   $a^p,M,b^{p+1}$ are semi-stable and that one of the two systems of inequalities in \eqref{middle M cap left M' sys} holds. In both the cases the given inequalities imply the inequalities \eqref{middle M cap left M'1}, therefore  $\sigma \in (a^p,M,b^{p+1})$.  If the second system in \eqref{middle M cap left M' sys} holds, then by Lemma \ref{semi-stability of a} (b) we get $\sigma \in (a^p,a^{p+1},M)\subset \mk{T}_{a}^{st}$. If   the first  system in \eqref{middle M cap left M' sys} holds, then by Lemma \ref{semi-stability of a} (c) and (d) we get $\sigma \in (a^j,a^{j+1},M)\subset \mk{T}_{a}^{st}$ for some $j\in \ZZ$, and the inclusion $\supset$ is proved as well. 

As in the last paragraph of the proof of Lemma \ref{T12Zcap(E_1)} one shows that the   two systems of  inequalities in \eqref{middle M cap left M' sys} correspond to two contractible sets (the first is contractible by Corollary  \ref{contractible 1a}),  and it is easy to show that their intersection is homeomorphic to $\RR_{>0}^3 \times \{\phi_2-1<\phi_0<\phi_1 < \phi_2\}$, which is also contractible.  Remark \ref{VK} shows that  $ (a^p,M,b^{p+1}) \cap \mk{T}_{a}^{st}$ is contractible.   Since  $  (a^p,M,b^{p+1}) $ and  $ \mk{T}_{a}^{st}$ are both contractible (Proposition \ref{lemma for f_E(Theta_E)} and Corollary  \ref{T12Zcap union(E_1)}), Remark \ref{VK} shows that  $   (a^p,M,b^{p+1}) \cup \mk{T}_{a}^{st}$ is contractible as well. 
\epr
\begin{lemma} \label{middle M cap left M} For any $p\in \ZZ$ the set  $(a^p,M,b^{p+1}) \cap \mk{T}_{b}^{st}$ consists of the stability conditions $\sigma$ for which  $a^p,M,b^{p+1}$ are semistable and:   \begin{gather} \label{middle M cap left M sys} \begin{array}{c}\phi\left (a^{p}\right)-1< \phi\left ( M\right)<\phi\left (a^{p}\right)\\
 \phi\left (a^{p}\right)-1< \phi\left ( b^{p+1}\right)-1<\phi\left (a^{p}\right) \\ \arg_{(\phi\left (a^{p}\right)-1,\phi\left (a^{p}\right))}(Z(a^p)-Z(b^{p+1}))>\phi(M) \end{array}    \ \mbox{or} \  \begin{array}{c}  \phi\left (a^{p}\right)-1< \phi\left ( M\right)<\phi\left (a^{p}\right)\\ \phi(M)+1<\phi(b^{p+1})\end{array}.   \end{gather}
It follows that $(a^p,M,b^{p+1}) \cap \mk{T}_{b}^{st}$ and  $   (a^p,M,b^{p+1}) \cup \mk{T}_{b}^{st}$ are  contractible. 
\end{lemma}
\bpr We start with the inclusion $\subset$. Assume that $\sigma \in (a^p,M,b^{p+1})$. Then $a^p,M,b^{p+1}$ are semi-stable and by table \eqref{middle M} we get  \begin{gather} \label{middle M cap left M1} \begin{array}{c} \phi\left ( a^{p}\right) < \phi\left (M\right)+1 \\ \phi\left ( a^{p}\right) < \phi\left (b^{p+1}\right) \\ \phi\left ( M\right) < \phi\left (b^{p+1}\right)\end{array} \end{gather}
Recalling \eqref{T_43new},  we see that we have to consider three cases.

\ul{If $\sigma \in (M,b^j, b^{j+1})$}, then $M,b^j, b^{j+1}$ are semi-stable and   from table \eqref{left right M} we see that $\phi(M)<\phi(b^{j})$ and   $\phi(M)+1<\phi(b^{j+1})$, hence $\phi(a^p)<\phi(b^{j+1})$ and $\hom(b^{j+1},a^p)=0$. From \eqref{nonvanishing3} it follows that  $p\leq j$. If $j=p$, then  $\phi(M)+1<\phi(b^{p+1})$ and  by   $\hom(b^{p},a^{p})\neq 0$ (see  \eqref{nonvanishing3}) we get $\phi(M)<\phi(a^{p})$, which implies the second system in \eqref{middle M cap left M sys}. It remains to consider the case  $p< j$.

 In this case $\hom^1(b^{j+1}, a^p)\neq 0$ (see \eqref{nonvanishing3}) and we obtain $\phi(M)+1<\phi(b^{j+1})\leq \phi(a^{p})+1$, which combined with \eqref{middle M cap left M1} implies 
$\phi\left (a^{p}\right)-1< \phi\left ( M\right)<\phi\left (a^{p}\right)$. 
On the other hand, we have  $\phi(b^{j})<\phi(b^{j+1})$ (see table \eqref{left right M}), and  by $p<j$ we can write $\phi(b^{p+1})\leq \phi(b^{j})<\phi(b^{j+1})\leq  \phi(a^{p})+1$, which combined with \eqref{middle M cap left M1} implies  $ \phi\left (a^{p}\right)-1< \phi(b^{p+1})-1<\phi\left (a^{p}\right)$.  Now we can use Lemma \ref{semi-stability of b} (c) to deduce that $M'\in \sigma^{ss}$ and $\phi(M')=\arg_{(\phi\left (a^{p}\right)-1,\phi\left (a^{p}\right))}(Z(a^p)-Z(b^{p+1}))$. From $\hom^1(b^{j+1},M')\neq 0$ and $\phi(M)+1<\phi(b^{j+1})$ it follows that $\phi(M)<\phi(M')$ and   the first system in \eqref{middle M cap left M sys} follows.

\ul{If $\sigma\in  (b^m,b^{m+1},M')$}, then $b^m,b^{m+1}, M'$ are  semistable  and in table \eqref{left right M} we see that $\phi(b^m)<\phi(M')$. By $\hom(M',a^p)\neq 0$ and  $\hom(M,b^m)\neq 0$   (see \eqref{nonvanishing1}) we get:

\begin{gather} \label{middle M cap left M11} \phi(M)\leq \phi(b^m)<\phi(M')\leq \phi(a^p). \end{gather} 
Whence $\phi(M)<\phi(a^p)$ and combining with \eqref{middle M cap left M1} we derive  $ \phi\left (a^{p}\right)-1< \phi(M)<\phi\left (a^{p}\right)$. On the other hand, in \eqref{middle M cap left M11} we have also  $\phi(b^m)<\phi(a^p)$, and hence  $\hom(a^p,b^m)=0$,  threfore  by  \eqref{nonvanishing4} we see   that $p\geq m$.   In \eqref{middle M cap left M11} we have also $\phi(M)<\phi(M')$. Taking into account Lemma \ref{semi-stability of b} (c), we see that if we show that  $ \phi\left (a^{p}\right)-1< \phi\left (b^{p+1}\right)-1<\phi\left (a^{p}\right)$, then the first system in \eqref{middle M cap left M sys} follows. Since we have $\phi\left (a^{p}\right)< \phi\left (b^{p+1}\right)$ (see \eqref{middle M cap left M1}), it remains to show that $\phi\left (b^{p+1}\right)<\phi\left (a^{p}\right)+1$. If $p=m$, then from table \eqref{left right M} we obtain $\phi(b^{p+1})<\phi(M')+1$ and the inequality in question follows from $\phi(M')\leq \phi(a^p)$. If $m<p$, then $\hom^1(b^{p+1},b^m)\neq 0$ and  we get   $\phi(b^{p+1})\leq \phi(b^{m})+1<\phi\left (a^{p}\right)+1$ (see \eqref{middle M cap left M11}).

\ul{If $\sigma\in  (b^m, a^{m},b^{m+1})$}, then $b^m, a^{m},b^{m+1} \in \sigma^{ss}$ and in table \eqref{no M} we see that $\phi(b^m)+1<\phi(b^{m+1})$,  hence Lemma \ref{from phases of big distance} and $b^{p+1} \in \sigma^{ss}$ imply that $p=m$ or $p=m-1$. If $p=m-1$, then by \eqref{middle M cap left M1} we obtain $\phi(a^{m-1})+1<\phi(b^m)+1<\phi(b^{m+1})$, and hence $\hom^1(b^{m+1},a^{m-1})=0$, which contradicts \eqref{nonvanishing3}. Therefore we have  $m=p$. Now we have $\phi(b^p)+1<\phi(b^{p+1})$ and $\phi(b^{p})<\phi(a^{p})$ (see table \eqref{no M}), which  together with $\hom(M,b^p)\neq 0$ imply $\phi(M)+1<\phi(b^{p+1})$ and $\phi(M)<\phi(a^{p})$,  hence the second system  in \eqref{middle M cap left M sys} follows. Thus we showed the inclusion $\subset$. 

 We show now the inverse inclusion $\supset$. Assume that   $a^p,M,b^{p+1}$ are semi-stable and that one of the two systems of inequalities in \eqref{middle M cap left M sys} holds. In both the cases the given inequalities imply the inequalities \eqref{middle M cap left M1}, therefore  $\sigma \in (a^p,M,b^{p+1})$.  If the second system in \eqref{middle M cap left M sys} holds, then by Lemma \eqref{semi-stability of b} (b) we get $\sigma \in (M, b^p,b^{p+1}) \subset \mk{T}_{b}^{st}$. If   the first  system in \eqref{middle M cap left M sys} holds, then the desired  $\sigma \in \mk{T}_{b}^{st}$ follows  from Lemma \eqref{semi-stability of b} (c) and (d).  The inclusion $\supset$ is proved as well. 

In Corollary \ref{T43Zcap union(E_1)} was shown that $\mk{T}_{b}^{st}$ is contractible.  The proof that  $(a^p,M,b^{p+1}) \cap \mk{T}_{b}^{st}$ and  $   (a^p,M,b^{p+1}) \cup \mk{T}_{b}^{st}$ are contractible  
 is as in the last paragraph of Lemma \ref{middle M cap left M'}.  The two  systems  in \eqref{middle M cap left M sys} correspond to   contractible subsets of $(a^p,M,b^{p+1}) \cap \mk{T}_{b}^{st}$ (the first is contractible by  Corollary  \ref{contractible 1b}).  
The intersection of these subsets is homeomorphic to $\RR_{>0}^3 \times \{\phi_0-1<\phi_1<\phi_2-1 < \phi_0\}$, which is also contractible. Now we apply  Remark \ref{VK} twice and the lemma  follows. 
\epr
\begin{coro} \label{middle M cup left right M coro} For any $p\in \ZZ $  the set $\mk{T}_{a}^{st} \cup (a^p,M,b^{p+1})\cup \mk{T}_{b}^{st} $ is contractible. \end{coro}
\bpr In Lemma \ref{middle M cap left M'} we showed that $\mk{T}_{a}^{st} \cup (a^p,M,b^{p+1})$  is contractible. Since  $\mk{T}_{a}^{st} \cap  \mk{T}_{b}^{st} =\emptyset $ (see Subsection \ref{empty intersection}), it follows that  $(\mk{T}_{a}^{st} \cup (a^p,M,b^{p+1}))\cap \mk{T}_{b}^{st} = (a^p,M,b^{p+1})\cap \mk{T}_{b}^{st}$, which is contractible by Lemma \ref{middle M cap left M}. Now we apply Remark \ref{VK}. 
\epr

\begin{lemma} \label{middle M cap left right M}  For any $q<p$ the set  $(a^p,M,b^{p+1}) \cap (\mk{T}_{a}^{st} \cup (a^q,M,b^{q+1})\cup \mk{T}_{b}^{st}  )$  consists of the stability conditions $\sigma$ for which  $a^p,M,b^{p+1}$ are semistable and:   \begin{gather} \begin{array}{c}\phi\left (a^{p}\right)-1< \phi\left ( M\right)<\phi\left (a^{p}\right)\\
 \phi\left (a^{p}\right)-1< \phi\left ( b^{p+1}\right)-1<\phi\left (a^{p}\right)  \end{array}    \ \mbox{or} \  \begin{array}{c}  \phi\left (a^{p}\right)-1< \phi\left ( M\right)<\phi\left (a^{p}\right)\\ \phi(M)+1<\phi(b^{p+1})\end{array} \nonumber  \\[-2mm]  \label{middle M cap left right M sys} \\[-2mm]  \mbox{or}  \begin{array}{c}  \phi\left (a^p\right)< \phi\left ( M\right)\\ \phi\left (b^{p+1}\right)-1< \phi\left ( M\right)<\phi\left (b^{p+1}\right)\end{array}   \ \mbox{or} \  \begin{array}{c}  \phi\left (b^{p+1}\right)-1< \phi\left ( M\right)<\phi\left (b^{p+1}\right)\\  \phi\left (b^{p+1}\right)-1< \phi\left ( a^{p}\right)<\phi\left (b^{p+1}\right) \\ \arg_{(\phi\left (b^{p+1}\right)-1,\phi\left (b^{p+1}\right))}(Z(a^p)-Z(b^{p+1}))<\phi(M)  \nonumber \end{array}.
 \end{gather}   
It follows that  $(a^p,M,b^{p+1}) \cap (\mk{T}_{a}^{st} \cup (a^q,M,b^{q+1})\cup \mk{T}_{b}^{st}  )$  and  $(a^p,M,b^{p+1}) \cup (\mk{T}_{a}^{st} \cup (a^q,M,b^{q+1})\cup \mk{T}_{b}^{st}  )$   are  contractible. 
\end{lemma}
\bpr   We start with the inclusion $\subset$. Assume that $\sigma \in (a^p,M,b^{p+1})$. Then $a^p,M,b^{p+1}$ are semi-stable and by table \eqref{middle M} we get  \begin{gather} \label{middle M cap left right M1} \begin{array}{c} \phi\left ( a^{p}\right) < \phi\left (M\right)+1 \\ \phi\left ( a^{p}\right) < \phi\left (b^{p+1}\right) \\ \phi\left ( M\right) < \phi\left (b^{p+1}\right)\end{array}. \end{gather}

\ul{If $\sigma \in  (a^q,M,b^{q+1})$ and $q<p$,} then $a^q,b^{q+1} \in \sigma^{ss}$ and   $\phi(M)<\phi(b^{q+1})$,  $\phi\left ( a^{q}\right) < \phi\left (b^{q+1}\right) $ a well.  By \eqref{nonvanishing3}   we have  $\hom(b^{q+1}, a^p)\neq 0$ and  $\hom^1(b^{p+1}, a^q)\neq 0$,  therefore $\phi(M)<\phi(b^{q+1})\leq \phi(a^{p})$ and $\phi(b^{p+1})\leq \phi(a^{q})+1<\phi(b^{q+1})+1\leq \phi\left (a^{p}\right)+1$. Combining with \eqref{middle M cap left right M1} we obtain the  system in the first row and first column in \eqref{middle M cap left right M sys}. 

\ul{If $\sigma \in \mk{T}_{a}^{st}$}, then by  Lemma \ref{middle M cap left M'} some of the systems on the second row of \eqref{middle M cap left right M sys} follows. 

\ul{If $\sigma \in \mk{T}_{b}^{st}$}, then by  Lemma \ref{middle M cap left M}  some of the systems on the first row of \eqref{middle M cap left right M sys} follows (\eqref{middle M cap left M sys} implies \eqref{middle M cap left right M sys}).  So we showed the inclusion $\subset$.

We show now the inclusion $\supset$. So let $a^p,M,b^{p+1}$ be semi-stable.  If some of the systems on the second row of \eqref{middle M cap left right M sys} holds,  then by   \ref{middle M cap left M'} it follows that $\sigma \in (a^p,M,b^{p+1})\cap \mk{T}_{a}^{st}$. If the  system in the first row and second column of \eqref{middle M cap left right M sys} holds,  then Lemma \ref{middle M cap left M} ensures that $\sigma \in (a^p,M,b^{p+1})\cap \mk{T}_{b}^{st}$.  

Thus, it remains to consider the first system in \eqref{middle M cap left right M sys}. We assume till the end of the proof that 
\begin{gather} \label{middle M cap left right M2}
\begin{array}{c}\phi\left (a^{p}\right)-1< \phi\left ( M\right)<\phi\left (a^{p}\right)\\
 \phi\left (a^{p}\right)-1< \phi\left ( b^{p+1}\right)-1<\phi\left (a^{p}\right)  \end{array}.
\end{gather}
Lemma \ref{semi-stability of b} (c) ensures  that 
 \begin{gather} \label{middle M cap left right M3} M'\in \sigma^{ss};  \qquad \phi(b^{p+1})-1<\phi(M')=\arg_{(\phi\left (a^{p}\right)-1,\phi\left (a^{p}\right))}(Z(a^p)-Z(b^{p+1}))<\phi(a^{p}). \end{gather}
Now we consider three cases.

If $\phi(M')>\phi(M)$, then  \eqref{middle M cap left right M2} and \eqref{middle M cap left right M3} yield the  first system in \eqref{middle M cap left M sys}  is satisfied and then Lemma \ref{middle M cap left M}   says that $\sigma \in (a^p,M,b^{p+1})\cap \mk{T}_{b}^{st}$.

If $\phi(M')<\phi(M)$,  then by $\hom^1(b^{p+1}, M')\neq 0$ it follows that $\phi(b^{p+1})-1<\phi(M)$. Combining  this inequality with  \eqref{middle M cap left right M2} one easily shows that: 
\begin{gather} \label{middle M cap left right M6}\begin{array}{c}\phi\left (b^{p+1}\right)-1< \phi\left ( M\right)<\phi\left (b^{p+1}\right)\\
 \phi\left (b^{p+1}\right)-1< \phi\left ( a^{p}\right)<\phi\left (b^{p+1}\right)  \end{array}. \end{gather}
Having obtained \eqref{middle M cap left right M6} we can use Lemma \ref{semi-stability of a}  (c)  and due to $\phi(M')<\phi(M)$ we derive  the first system  in  \eqref{middle M cap left M' sys}.   Thus Lemma \ref{middle M cap left M'}   ensures that $\sigma \in (a^p,M,b^{p+1})\cap \mk{T}_{a}^{st}$.

Finally, if $\phi(M)=\phi(M')$, then due to \eqref{middle M cap left right M2} we can apply Lemma \ref{semi-stability of b} \textbf{(e)}, which says that $\sigma \in (a^p,M,b^{p+1})\cap (a^q,M,b^{q+1})$ (recall that $q<p$). So far we showed the first part of the lemma.

We explain  now, using  the obtained representation through the systems of inequalities \eqref{middle M cap left right M sys}, that     $(a^p,M,b^{p+1}) \cap (\mk{T}_{a}^{st} \cup (a^q,M,b^{q+1})\cup \mk{T}_{b}^{st}  )$ is contractible.   The four systems correspond to four open  subsets of  $(a^p,M,b^{p+1}) \cap (\mk{T}_{a}^{st} \cup (a^q,M,b^{q+1})\cup \mk{T}_{b}^{st}  )$ (see the last paragraph of the proof of Lemma \ref{T12Zcap(E_1)}). We denote these subsets by $S_{11} , S_{12}, S_{21}, S_{22}$, where $S_{ij}$ corresponds to the system in the $i$-th row  and $j$-th column of \eqref{middle M cap left right M sys}.  The proved part of the lemma is the equality  $(a^p,M,b^{p+1}) \cap (\mk{T}_{a}^{st} \cup (a^q,M,b^{q+1})\cup \mk{T}_{b}^{st}  )=\bigcup_{1\leq i,j \leq 2} S_{ij}$. The subset $S_{22}$ is contractible by Corollary  \ref{contractible 1a}. The subsets $S_{11} , S_{12}, S_{21}$ are contractible since they are homeomorphic to convex subsets of $\RR^6$. For example $S_{11}$ is homeomorphic to $$\RR_{>0}^3\times \left \{ (\phi_0,\phi_1,\phi_2)  \in \RR^3:  \begin{array}{c}\phi_0-1< \phi_1<\phi_0\\
 \phi_0-1< \phi_2-1<\phi_0  \end{array}\right \}. $$

One easily shows that $S_{11}\cap S_{12}$ is homeomorphic to  $\RR_{>0}^3 \times \{\phi_0-1<\phi_1<\phi_2-1 < \phi_0\}$, hence it is contractible, and by Remark \ref{VK} we deduce that $S_{11}\cup S_{12}$ is contractible. Note that in $S_{12}$ we have $\phi(M)+1<\phi(b^{p+1})$ and in $S_{22}$ we have  $\phi(M)+1>\phi(b^{p+1})$ , therefore $S_{12}\cap S_{22}=\emptyset$. Hence $S_{22}\cap (S_{11}\cup S_{12})=S_{22}\cap S_{11}$. One easily shows that $S_{22}\cap S_{11}$ is  homeomorphic to:
\begin{gather}   \RR_{>0}^3 \times \left \{ (\phi_0,\phi_1,\phi_2) \in \RR^3: \begin{array}{c} r_i>0\\ \phi_2-1 < \phi_1 < \phi_0 < \phi_2 \\   \arg_{(\phi_2-1,\phi_2)}(r_0 \exp(\ri \pi \phi_0) - r_2 \exp(\ri \pi \phi_2))<\phi_1\end{array} \right \},\end{gather} 
 which by Corollary \ref{contractible 2a} is contractible as well. Thus, we see that $S_{22}\cap (S_{11}\cup S_{12})$ is contractible, therefore by Remark \ref{VK} we see that  $S_{22}\cup S_{11}\cup S_{12}$ is contractible.  In $S_{11}$ and  $S_{12}$   we have $\phi(M)<\phi(a^p)$ and in $S_{21}$ we have  $\phi(M)>\phi(a^p)$, therefore $S_{21}\cap(S_{22}\cup S_{11}\cup S_{12})=S_{21}\cap S_{22}$. On the other hand, one easily shows (by drawing a picture) that the intersection $S_{21}\cap S_{22}$ is homeomorphic to  $\RR_{>0}^3 \times \{\phi_2-1<\phi_0<\phi_1 < \phi_2\}$,  which is contractible as well, and hence $S_{21}\cap(S_{22}\cup S_{11}\cup S_{12})$ is contractible.  Applying Remark \ref{VK} again ensures that     $S_{21}\cup S_{22}\cup S_{11}\cup S_{12} = (a^p,M,b^{p+1}) \cap (\mk{T}_{a}^{st} \cup (a^q,M,b^{q+1})\cup \mk{T}_{b}^{st}  )$ is contractible. In  Corollary  \ref{middle M cup left right M coro} is shown that  $\mk{T}_{a}^{st} \cup (a^q,M,b^{q+1})\cup \mk{T}_{b}^{st}$ is contractible and  with  one more  reference to Remark \ref{VK} we prove the lemma. 
\epr

\begin{coro} \label{middle M cup left right M} The set $(\_,M,\_)\cup \mk{T}_{a}^{st}\cup \mk{T}_{b}^{st}$ is contractible. 
\end{coro}
\bpr Recall that $ (\_,M\_) =\bigcup_{q\in \ZZ}(a^q,M,b^{q+1})$ (see  \eqref{middle M ab}).  We will prove  that for each $p\in \ZZ$ and for each $k\geq 1$  the set \eqref{middle M cup left right M1} below  is contractible, and the corollary follows from Remark \ref{VK}: 
\begin{gather} \label{middle M cup left right M1} \bigcup_{i=0}^k  (a^{p-i},M,b^{p+1-i})\cup (\mk{T}_{a}^{st} \cup \mk{T}_{b}^{st}  ).\end{gather}

In the previous lemma was shown  that for $k=1$ and any $p\in \ZZ$ the set \eqref{middle M cup left right M1} is contractible.  Assume that for some $k\geq 1$  this set is contractible for each $p\in \ZZ$. Take now any $p\in \ZZ$. We have 
\begin{gather}\label{middle M cup left right M2}\bigcup_{i=0}^{k+1}  (a^{p-i},M,b^{p+1-i})\cup (\mk{T}_{a}^{st} \cup \mk{T}_{b}^{st}  )= (a^{p},M,b^{p+1})\cup \left (\bigcup_{i=1}^{k+1}  (a^{p-i},M,b^{p+1-i})\cup (\mk{T}_{a}^{st} \cup \mk{T}_{b}^{st}  ) \right ). \end{gather}
 Proposition \ref{lemma for f_E(Theta_E)}  and the induction assumption say that  the two components on RHS of \eqref{middle M cup left right M2} are contractible. Since the intersection analyzed  in  Lemma \ref{middle M cap left right M}  does not depend on $q$, we can write: 
$$(a^{p},M,b^{p+1})\cap\left ( \bigcup_{i=1}^{k+1}  (a^{p-i},M,b^{p+1-i})\cup (\mk{T}_{a}^{st} \cup \mk{T}_{b}^{st}  ) \right )=(a^{p},M,b^{p+1})\cap\left (   (a^{p-1},M,b^{p})\cup (\mk{T}_{a}^{st} \cup \mk{T}_{b}^{st} )  \right ),  $$
which by Lemma \ref{middle M cap left right M} is contractible. Now  Remark \ref{VK} ensures that that \eqref{middle M cup left right M2} is contractible. 
\epr

The next step is to glue   $(\_,M,\_)\cup \mk{T}_{a}^{st}\cup \mk{T}_{b}^{st}$ and $(b^p,M',a^p)$. This is done in several substeps: Lemmas \ref{middle M' cap left M}, \ref{middle M' cap left M'}, \ref{semi-stability of a'}, \ref{semi-stability of b'}, which lead to Corollary \ref{middle M' cap left right middle M}.    In  the next two lemmas we prove inclusions in only one direction not equality of sets.

\begin{lemma} \label{middle M' cap left M} Let  $p\in \ZZ$.  If  $\sigma \in (b^p,M',a^{p}) \cap \mk{T}_{b}^{st}$, then  $b^p,M',a^{p}$ are semistable and:   \begin{gather} \label{middle M' cap left M sys} \begin{array}{c}\phi\left (a^{p}\right)-1< \phi\left ( M'\right)<\phi\left (a^{p}\right)\\
 \phi\left (a^{p}\right)-1< \phi\left ( b^{p}\right)<\phi\left (a^{p}\right) \end{array}    \ \mbox{or} \  \begin{array}{c}  \phi\left (a^{p}\right)-1< \phi\left ( M'\right)<\phi\left (a^{p}\right)\\ \phi(b^{p})<\phi(M')\end{array}.   \end{gather}
\end{lemma}
\bpr  In table \eqref{middle M} we see that $b^p,M',a^{p}$ are semi-stable and: 
\begin{gather} \label{middle M' cap left M1}  \begin{array}{c} \phi\left ( b^{p}\right) < \phi\left (M'\right) +1\\ \phi\left ( b^{p}\right) < \phi\left (a^{p}\right) \\ \phi\left ( M'\right) < \phi\left (a^{p}\right)    
	\end{array} \end{gather}

Recalling \eqref{T_43new},  we see that we have to consider three cases.

\ul{If $\sigma \in (M,b^j, b^{j+1})$}, then $M,b^j, b^{j+1}$ are semi-stable and   from table \eqref{left right M} we see that $\phi(M)<\phi(b^{j})$ and   $\phi(M)+1<\phi(b^{j+1})$.  By $\hom^1(a^p,M)\neq 0$  and  $\hom^1(b^{j+1},M')\neq 0$ (see \eqref{nonvanishing2}) we can write  $\phi(a^p)\leq \phi(M)+1 <\phi(b^{j+1})\leq \phi(M')+1$, therefore (see also \eqref{middle M' cap left M1})  we get
\begin{gather} \label{middle M' cap left M2} \phi\left (a^{p}\right)-1< \phi\left ( M'\right)<\phi\left (a^{p}\right).\end{gather}

Since $\phi(b^p)<\phi(a^p)\leq \phi(M)+1<\phi(b^{j+1})$,  due to   \eqref{nonvanishing6} the inequality   $p\leq j$ must hold. 

If $j=p$, then the inequality  $\phi(M)<\phi(b^{p})$ (coming from  $\sigma \in (M,b^j, b^{j+1})$)   implies   $\phi(a^p)-1\leq \phi(M) <\phi(b^{p})$ and  combining with \eqref{middle M' cap left M1} and \eqref{middle M' cap left M2}  we obtain the first system in \eqref{middle M' cap left M sys}. 

 If $p<j$, then we have   $\hom(a^p, b^j)\neq 0$ (see \eqref{nonvanishing4}) and $\hom^1(b^{j+1}, b^p)\neq 0$, hence $\phi(a^p)\leq \phi(b^j)<\phi(b^{j+1})\leq  \phi(b^p)+1$ and again the  first system in \eqref{middle M' cap left M sys} follows.

\ul{If $\sigma\in  (b^m,b^{m+1},M')$}, then $b^m,b^{m+1}, M'$ are  semistable  and in table \eqref{left right M} we see that $\phi(b^m)<\phi(M')$,  therefore $\phi(b^m)<\phi(M')<\phi(a^p)$ and $\hom(a^p,b^m)=0$. From \eqref{nonvanishing4} we deduce that $m\leq p$.  

If $m=p$, then the incidence $\sigma\in  (b^m,b^{m+1},M')$  gives $\phi(b^p)<\phi(M')$ and $\phi(b^{p+1})-1<\phi(M')$ (see table \eqref{left right M}), and from $\hom(a^p,b^{p+1})\neq 0$ we obtain $\phi(a^{p})-1<\phi(M')$, therefore the second system in  \eqref{middle M' cap left M sys} holds. 
 
Let $m<p$. Then we have $\phi(b^m)<\phi(b^{m+1})$ and $\phi(b^m)<\phi(M')$ (see table \eqref{left right M}).  Using  $\hom^1(a^p,b^m)\neq 0$ (see \eqref{nonvanishing4})  we deduce  $\phi(a^p)\leq \phi(b^m)+1<\phi(b^{m+1})+1\leq \phi(b^{p})+1$ and $\phi(a^p)\leq \phi(b^m)+1<\phi(M')+1$, which combined with \eqref{middle M' cap left M1} produces the first system in \eqref{middle M' cap left M sys}. 

\ul{If $\sigma\in  (b^m, a^{m},b^{m+1})$}, then $b^m, a^{m},b^{m+1} \in \sigma^{ss}$ and in table \eqref{no M} we see that $\phi(b^m)+1<\phi(b^{m+1})$,  hence  Lemma \ref{from phases of big distance} and $b^{p} \in \sigma^{ss}$ imply  $p=m$ or $p=m+1$. If $p=m+1$, then by \eqref{middle M' cap left M1} we obtain $\phi(b^m)+1<\phi(b^{m+1})<\phi(a^{m+1})$, and hence $\hom^1(a^{m+1},b^{m})=0$, which contradicts \eqref{nonvanishing4}. Therefore we have  $m=p$.  In table \eqref{no M}  we see that  $\phi(b^p)+1<\phi(b^{p+1})$ and $\phi(a^p)<\phi(b^{p+1})$. From $\hom^1(b^{p+1}, M')\neq 0$ it follows that $ \phi(b^p)+1<\phi(b^{p+1})\leq \phi(M')+1$ and $\phi(a^p)<\phi(b^{p+1})\leq \phi(M')+1$.  These inequalities together with \eqref{middle M' cap left M1} produce  the second system in \eqref{middle M' cap left M sys}.
\epr

\begin{lemma} \label{middle M' cap left M'} Let  $p\in \ZZ$.  If  $\sigma \in (b^p,M',a^{p}) \cap \mk{T}_{a}^{st}$, then  $b^p,M',a^{p}$ are semistable and:   \begin{gather} \label{middle M' cap left M' sys}\begin{array}{c}  \phi\left (b^{p}\right)-1< \phi\left ( M'\right)<\phi\left (b^{p}\right)\\  \phi\left (b^{p}\right)-1< \phi\left ( a^{p}\right)-1<\phi\left (b^{p}\right)  \end{array}    \ \mbox{or} \  \begin{array}{c}  \phi\left (b^{p}\right)-1< \phi\left ( M'\right)<\phi\left (b^{p}\right)\\  \phi\left ( M'\right)+1<\phi\left (a^p\right)\end{array}.   \end{gather}
\end{lemma}
\bpr  In table \eqref{middle M} we see that $b^p,M',a^{p}$ are semi-stable and:   \begin{gather} \label{middle M' cap left M'1} \begin{array}{c} \phi\left ( b^{p}\right) < \phi\left (M'\right) +1\\ \phi\left ( b^{p}\right) < \phi\left (a^{p}\right) \\ \phi\left ( M'\right) < \phi\left (a^{p}\right)    
	\end{array} \end{gather}
Recalling \eqref{T_12new},  we see that we have to consider the following  three cases.

\ul{If $\sigma \in (M',a^j, a^{j+1})$}, then $M',a^j, a^{j+1}$ are semi-stable and  $\phi(M')<\phi(a^{j})$, $\phi(M')+1<\phi(a^{j+1})$, $\phi(a^{j})<\phi(a^{j+1})$ (see table   \eqref{left right M}). On the other hand $\phi(b^p)<\phi(M')+1$, hence $\hom(a^{j+1},b^p)=0$ and \eqref{nonvanishing4} implies that  $p-1\leq j$.   If $p-1=j$, then we have  $\phi(M')+1<\phi(a^{p})$ and $\phi(M')<\phi(a^{p-1}) \leq \phi(b^p)$ (see also \eqref{nonvanishing4}) and combining with \eqref{middle M' cap left M'1} we derive the second system in \eqref{middle M' cap left M' sys}.  
Let $p\leq j$.  Then by \eqref{nonvanishing4} we have $\hom^1(a^{j+1},b^p)\neq 0$ and we can write $\phi(M')+1<\phi(a^{j+1})\leq\phi(b^p)+1$ and $\phi(a^p)\leq \phi(a^j)<\phi(a^{j+1})\leq \phi(b^p)+1$, therefore $\phi(M')<\phi(b^p)$ and $\phi(a^p)<\phi(b^p)+1$, which combined with \eqref{middle M' cap left M'1} amounts to the first system in \eqref{middle M' cap left M' sys}.

\ul{If $\sigma\in  (a^m,a^{m+1},M)$}, then $a^m,a^{m+1}, M$ are semistable as well and in table \eqref{left right M} we see that $\phi(a^m)<\phi(M)$, $\phi(a^{m+1})<\phi(M)+1$, $\phi(a^m)<\phi(a^{m+1})$. Since $\hom(M',a^m)\neq0$ and $\hom(M, b^p)\neq 0$, it follows 
that $\phi(M')\leq\phi(a^m)<\phi(M)\leq\phi( b^p)$ and hence (see also  \eqref{middle M' cap left M'1}):
\begin{gather} \label{middle M' cap left M'11} \phi\left (b^{p}\right)-1< \phi\left ( M'\right)<\phi\left (b^{p}\right) \end{gather}

On the other hand,  $\phi(a^m)<\phi(M)$ and   $\hom(M, b^p)\neq 0$ imply        that $\phi(a^m)<\phi( b^p)$ and $\hom(b^p,a^m)=0$. Now from \eqref{nonvanishing3} we  deduce that $m<p$. If $m=p-1$, then we have $\phi(a^{p})<\phi(M)+1\leq\phi( b^p)+1$, which together with  \eqref{middle M' cap left M'11} and \eqref{middle M' cap left M'1} amounts to the first system in  \eqref{middle M' cap left M' sys}.

If $m<p-1$, then $\hom^1(a^p, a^m)\neq 0$ and $\hom( a^{m+1},b^p)\neq 0$ (see \eqref{nonvanishing4}). Therefore we have $\phi(a^p)\leq \phi(a^m)+1<\phi(a^{m+1})+1\leq \phi(b^p)+1$ and  the first system in  \eqref{middle M' cap left M' sys}  follows again.

\ul{If $\sigma\in  (a^m, b^{m+1},a^{m+1})$}, then $a^m, b^{m+1},a^{m+1}\in \sigma^{ss}$ and in table \eqref{no M} we see that $\phi(a^m)+1<\phi(a^{m+1})$, hence Lemma \ref{from phases of big distance} and $a^p \in \sigma^{ss}$ imply $p=m$ or $p=m+1$. If $p=m$, then by \eqref{middle M' cap left M'1} we obtain $\phi(b^m)+1<\phi(a^m)+1<\phi(a^{m+1})$, and hence $\hom^1(a^{m+1},b^m)=0$, which contradicts \eqref{nonvanishing4}. Thus, it remains to consider the case $m=p-1$. Now we have $\phi(a^{p-1})+1<\phi(a^{p})$ and $\phi(a^{p-1})<\phi(b^{p})$(see table \eqref{no M}), which  together with $\hom(M',a^{p-1})\neq 0$ imply $\phi(M')+1<\phi(a^{p})$ and $\phi(M')<\phi(b^{p})$,  hence we obtain the second system of inequalities in \eqref{middle M' cap left M' sys}.
\epr

\begin{lemma} \label{semi-stability of a'} Let $\sigma \in (b^p,M',a^{p})$ and let the following inequality hold:  \begin{gather} \label{semi-stability of a' ineq}  \phi\left (b^{p}\right)-1< \phi\left ( M'\right)<\phi\left (b^{p}\right).   \end{gather}
 Then we have the following:

{\rm \textbf{(a)}} $a^{p-1}\in \sigma^{ss}$ and   $\phi(M')<\phi(a^{p-1})<\phi(b^{p})<\phi(a^{p})$. 

{\rm \textbf{(b)}} If in addition to  \eqref{semi-stability of a' ineq} we have $\phi\left ( M'\right)+1<\phi\left (a^p\right)$, then $\sigma \in (M',a^{p-1},a^{p})$.

{\rm \textbf{ (c)}} If in addition to \eqref{semi-stability of a' ineq} we have  $\phi\left (b^{p}\right)-1< \phi\left ( a^{p}\right)-1<\phi\left (b^{p}\right),$ then $\sigma \in (a^{p-1},M,b^p)$.
\end{lemma}
\bpr  (a)  We  apply Proposition \ref{phi_1>phi_2} (a) to the triple $ (b^p,M',a^{p})$ and since $a^{p-1}$ is in the extension closure of $M', b^{p}$ (by  \eqref{short filtration 1}) it follows that  $a^{p-1}\in \sigma^{ss}$, $\phi(M')\leq \phi(a^{p-1})\leq \phi(b^{p})$ .  The inequality   $\phi(M')<\phi(a^{p-1})<\phi(b^{p})$  follows from the given inequality \eqref{semi-stability of a' ineq}, formula \eqref{phase formula} and $Z(a^{p-1})=Z(M')+Z(b^{p})$.  The inequality  $\phi(b^{p})<\phi(a^{p})$ follows from $\sigma \in (b^p,M',a^{p})$ (see table \eqref{middle M}). 

(b) From the given inequalities  and (a) we have  $\phi(M')<\phi(a^{p-1})$,  $\phi\left ( M'\right)+1<\phi\left (a^p\right)$, and  $\phi(a^{p-1})<\phi(a^{p})$,  then  table \eqref{left right M} shows  that  $\sigma \in (M',a^{p-1},a^{p})$. 

    (c)    From  Lemma \ref{two comp factors} applied  to the Ext-triple $(b^p,M',a^{p}[-1]) $ and the triangle \eqref{short filtration 2} we obtain  $M\in \sigma^{ss}$ and  $\phi(a^{p})-1<\phi(M)<\phi(b^{p})$. In (a) we got $a^{p-1}\in \sigma^{ss}$ and $\phi(a^{p-1})<\phi(a^{p})$, therefore  $\phi(a^{p-1})<\phi(M)+1$. In (a) we have also			$\phi(a^{p-1})<\phi(b^p)$. 
			Looking at table \eqref{middle M} we see that  $\sigma \in  (a^{p-1},M,b^p)$.
\epr

\begin{lemma} \label{semi-stability of b'} Let $\sigma \in (b^p,M',a^{p})$ and let the following inequality hold:  \begin{gather} \label{semi-stability of b' ineq} \phi\left (a^{p}\right)-1< \phi\left ( M'\right)<\phi\left (a^{p}\right).   \end{gather}
 Then we have the following:

{\rm \textbf{(a)}} $b^{p+1}\in \sigma^{ss}$ and   $\phi\left (b^{p}\right)-1< \phi\left (a^{p}\right)-1<\phi\left (b^{p+1}\right)-1 <\phi\left ( M'\right)$. 

{\rm \textbf{(b)}} If in addition to  \eqref{semi-stability of b' ineq} we have $ \phi(b^{p})<\phi(M')$, then $\sigma \in (b^p,b^{p+1},M')$.

{\rm \textbf{ (c)}} If in addition to \eqref{semi-stability of b' ineq} we have $\phi\left (a^{p}\right)-1< \phi\left ( b^{p}\right)<\phi\left (a^{p}\right) $,  then  $\sigma \in  (a^{p},M,b^{p+1})$

\end{lemma}
\bpr  (a)  We  apply Proposition \ref{phi_1>phi_2} (b) to the triple $ (b^p,M',a^{p})$ and since $b^{p+1}[-1]$ is in the extension closure of $M', a^{p}[-1]$ (by  \eqref{short filtration 1}) it follows that  $b^{p+1}\in \sigma^{ss}$,  $ \phi\left (a^{p}\right)-1\leq \phi\left (b^{p+1}\right)-1 \leq \phi\left ( M'\right)$.  The inequality   $ \phi\left (a^{p}\right)-1<\phi\left (b^{p+1}\right)-1 <\phi\left ( M'\right)$  follows from the given inequality \eqref{semi-stability of b' ineq}, formula \eqref{phase formula}, and $Z(b^{p+1}[-1])=Z(M')+Z(a^{p}[-1])$.  The inequality  $\phi(b^{p})<\phi(a^{p})$ follows from $\sigma \in (b^p,M',a^{p})$.

(b) From the given inequalities  and (a) we have $\phi(b^{p})<\phi(b^{p+1})$, $\phi(b^{p})<\phi(M')$  and $\phi(b^{p+1})<\phi(M')+1$.  Now in table \eqref{left right M} we see that  $\sigma \in (b^p,b^{p+1},M')$. 

    (c)    From  Lemma \ref{two comp factors} applied  to the Ext-triple $(b^p,M',a^{p}[-1]) $ and the triangle \eqref{short filtration 2} we obtain  $M\in \sigma^{ss}$ and  $\phi(a^{p})-1<\phi(M)<\phi(b^{p})$. In (a) we showed that $b^{p+1}\in \sigma^{ss}$ and $\phi(a^{p})<\phi(b^{p+1})$,  $ \phi\left (b^{p}\right)<\phi\left (b^{p+1}\right)$.  Now all the conditions  determining  $(a^{p},M,b^{p+1})$ (given in table \eqref{middle M}) are satisfied, hence  $\sigma \in  (a^{p},M,b^{p+1})$.
\epr

\begin{coro} \label{middle M' cap left right middle M}  For any $p\in \ZZ$ the set  $(b^p,M',a^{p}) \cap (\mk{T}_{a}^{st} \cup (\_,M,\_)\cup \mk{T}_{b}^{st}  )$  consists of the stability conditions $\sigma$ for which  $ b^p,M',a^{p} $ are semistable and:   \begin{gather} \begin{array}{c}\phi\left (a^{p}\right)-1< \phi\left ( M'\right)<\phi\left (a^{p}\right)\\
 \phi\left (a^{p}\right)-1< \phi\left ( b^{p}\right)<\phi\left (a^{p}\right) \end{array}    \ \mbox{or} \  \begin{array}{c}  \phi\left (a^{p}\right)-1< \phi\left ( M'\right)<\phi\left (a^{p}\right)\\ \phi(b^{p})<\phi(M')\end{array}  \nonumber  \\[-2mm]  \label{middle M' cap left right middle M sys}  \\[-2mm]  \mbox{or}  \begin{array}{c}  \phi\left (b^{p}\right)-1< \phi\left ( M'\right)<\phi\left (b^{p}\right)\\  \phi\left (b^{p}\right)-1< \phi\left ( a^{p}\right)-1<\phi\left (b^{p}\right)  \end{array}    \ \mbox{or} \  \begin{array}{c}  \phi\left (b^{p}\right)-1< \phi\left ( M'\right)<\phi\left (b^{p}\right)\\  \phi\left ( M'\right)+1<\phi\left (a^p\right)\end{array}.
 \nonumber\end{gather}   
It follows that  $(b^p,M',a^{p}) \cap (\mk{T}_{a}^{st} \cup (\_,M,\_)\cup \mk{T}_{b}^{st}  )$   and   $(b^p,M',a^{p}) \cup (\mk{T}_{a}^{st} \cup (\_,M,\_)\cup \mk{T}_{b}^{st}  )$    are  contractible. 
\end{coro}
\bpr Due to Lemmas  \ref{middle M' cap left M} and \ref{middle M' cap left M'},  to show the inclusion $\subset$  it remains only to show that the incidence $\sigma \in (b^p,M',a^{p}) \cap  (\_,M,\_)$ implies some of the systems in \eqref{middle M' cap left right middle M sys}. Assume that  $\sigma \in (b^p,M',a^{p}) \cap  (a^q,M,b^{q+1})$ for some $q\in \ZZ$. From table \eqref{middle M} we see that $b^p,M',a^{p},a^q,M,b^{q+1}$ are semi-stable and: 
 \begin{gather}\label{middle M' cap left right middle M1} \begin{array}{c} \phi\left ( b^{p}\right) < \phi\left (M'\right) +1\\ \phi\left ( b^{p}\right) < \phi\left (a^{p}\right) \\ \phi\left ( M'\right) < \phi\left (a^{p}\right)    
	\end{array}    \ \mbox{and} \   \begin{array}{c} \phi\left ( a^{q}\right) < \phi\left (M\right)+1 \\ \phi\left ( a^{q}\right) < \phi\left (b^{q+1}\right) \\ \phi\left ( M\right) < \phi\left (b^{q+1}\right)\end{array}. \end{gather}
	If $p\leq q$, then the non-vanishings  $\hom(a^p,a^q)\neq 0$,  $\hom^1(b^{q+1},M')\neq 0$,  and    $\hom(M,b^p)\neq 0$ (see Corollary \ref{nonvanishings}) together with \eqref{middle M' cap left right middle M1} imply the following inequalities   $\phi(a^p)\leq \phi(a^q)<\phi(M)+1\leq \phi(b^p)+1$ and  $\phi(a^p)\leq \phi(a^q)<\phi(b^{q+1})\leq \phi(M')+1$, which combined with \eqref{middle M' cap left right middle M1} amount to the system in the first row and the first column  of \eqref{middle M' cap left right middle M sys}.

	If $q<p$, then  the non-vanishings   $\hom(M',a^q)\neq 0$, $\hom(b^{q+1},b^p)\neq 0$, and $\hom^1(a^p,M)\neq 0$ together with \eqref{middle M' cap left right middle M1} imply the inequalities $\phi(M')\leq \phi(a^q)<\phi(b^{q+1})\leq \phi(b^{p})$ and $\phi(a^p)\leq \phi(M)+1<\phi(b^{q+1})+1\leq \phi(b^{p})+1$.  The system in the second row and the first column  in \eqref{middle M' cap left right middle M sys} follows.	
	So far we showed the incusion $\subset$. 
	
	Assume that $b^p,M',a^{p}\subset \sigma^{ss}$ and that  \eqref{middle M' cap left right middle M sys} holds. Each of the systems in \eqref{middle M' cap left right middle M sys} contains in it the inequalities of $(b^p,M',a^{p})$ from table \eqref{middle M}, hence    $\sigma \in (b^p,M',a^{p})$.    Lemmas \ref{semi-stability of a'} and \ref{semi-stability of b'} ensure that  $\sigma \in (\mk{T}_{a}^{st} \cup (\_,M,\_)\cup \mk{T}_{b}^{st}  )$ as well and the first part of the corollary follows.  
	
	Now the arguments are analogous to those given in the end of the proof of Lemma \ref{middle M cap left right M}.

  The four systems  in \eqref{middle M' cap left right middle M sys} correspond to four open  subsets of  $(b^p,M',a^{p}) \cap (\mk{T}_{a}^{st} \cup (\_,M,\_)\cup \mk{T}_{b}^{st}  )$. We denote these subsets by $S_{11} , S_{12}, S_{21}, S_{22}$, where $S_{ij}$ corresponds to the system in the $i$-th row  and $j$-th. 	
The first part of the corollary and Remark \ref{VK} reduce the proof of the last statement  to proving that 
	 $\bigcup_{1\leq i,j \leq 2} S_{ij}$ is contractible . 
	
	All of  $S_{11} , S_{12}, S_{21}, S_{22}$ are contractible since they are homeomorphic to convex subsets of $\RR^6$.

One easily shows that:
\begin{itemize}
	\item $S_{11}\cap S_{12}$ is homeomorphic to  $\RR_{>0}^3 \times \{\phi_2-1<\phi_0<\phi_1 < \phi_2\}$
	\item $S_{21}\cap (S_{11}\cup S_{12})=S_{21}\cap  S_{11}$ is homeomorphic to  $\RR_{>0}^3 \times \{\phi_2-1<\phi_1<\phi_0 < \phi_2\}$
	\item  $S_{22}\cap (S_{11}\cup S_{12}\cup S_{21})=S_{22}\cap  S_{21}$ is homeomorphic to  $\RR_{>0}^3 \times \{\phi_0-1<\phi_1<\phi_2-1 < \phi_0\}$.
\end{itemize}
Since the obtained  subsets of $\RR^6$ are convex, in particular contractible, it follows by Remark \ref{VK} that  $\bigcup_{1\leq i,j \leq 2} S_{ij}$ is contractible. The corollary follows. 
\epr
We can prove now Theorem \ref{main theo}:
\begin{theorem}  $\st(D^b(Q ))$ is contractible.  
\end{theorem}
\bpr  Recall that $\st(\mc T)=\mk{T}_{a}^{st} \cup(\_,M',\_)\cup(\_,M,\_)\cup  \mk{T}_{b}^{st}  $ (see \eqref{st with T}). Recalling  \eqref{middle M ab} we  get:
\begin{gather} \label{one union} \st(D^b(Q ))=\mk{T}_{a}^{st} \cup(\_,M,\_)\cup  \mk{T}_{b}^{st} \cup  \bigcup_{k\in \ZZ}(b^{k},M',a^{k}). \end{gather}
Corollary  \ref{middle M cup left right M} says that   $ \mk{T}_{a}^{st}\cup(\_,M,\_)\cup \mk{T}_{b}^{st}$ is contractible and it remains to show that after adding $ \bigcup_{k\in \ZZ}(b^{k},M',a^{k})$ the result is still contractible. 

We first show that for any two integers $q>p$ we have:  \begin{gather} \label{one inclusion} (b^p,M',a^{p}) \cap (b^q,M',a^{q}) \subset  (b^p,M',a^{p}) \cap (\mk{T}_{a}^{st} \cup (\_,M,\_)\cup \mk{T}_{b}^{st}  ).
\end{gather}
Assume that $\sigma \in (b^p,M',a^{p}) \cap (b^q,M',a^{q})$. Then in table \eqref{middle M} we see that 
\begin{gather}\label{two middle M'} \begin{array}{c} \phi\left ( b^{p}\right) < \phi\left (M'\right) +1\\ \phi\left ( b^{p}\right) < \phi\left (a^{p}\right) \\ \phi\left ( M'\right) < \phi\left (a^{p}\right)    
	\end{array}    \ \mbox{and} \  \begin{array}{c} \phi\left ( b^{q}\right) < \phi\left (M'\right) +1\\ \phi\left ( b^{q}\right) < \phi\left (a^{q}\right) \\ \phi\left ( M'\right) < \phi\left (a^{q}\right)    
	\end{array}. \end{gather}
	 Since $p<q$, we have the  non-vanishings  $\hom(a^p,b^q)\neq 0$ and   $\hom^1(a^{q},b^{p})\neq 0$  (see \eqref{nonvanishing4}). We combine with \eqref{two middle M'} as follows   $\phi(a^p)\leq \phi(b^q)<\phi(a^q)\leq \phi(b^p)+1$ and  $\phi(a^p)\leq \phi(b^q)< \phi(M')+1$, hence  $\phi(a^p)-1< \phi(b^p)$  and  $\phi(a^p)-1<\phi(M')$. In \eqref{two middle M'} we have also $\phi(b^p)<\phi(a^p)$ and   $\phi(M')<\phi(a^p)$ and  the system in the first row and the first column  of \eqref{middle M' cap left right middle M sys} follows. Therefore  by Corollary \ref{middle M' cap left right middle M} we get  $\sigma \in  (b^p,M',a^{p}) \cap (\mk{T}_{a}^{st} \cup (\_,M,\_)\cup \mk{T}_{b}^{st}  )$ and we showed the inclusion \eqref{one inclusion}. This implies that for any $p\in \ZZ$ and any $n\geq 1$ holds the following equality: 
	\begin{gather}  (b^p,M',a^{p}) \cap\left ( \mk{T}_{a}^{st} \cup (\_,M,\_)\cup \mk{T}_{b}^{st} \bigcup_{k=1}^n(b^{p+k},M',a^{p+k}) \right )   = (b^p,M',a^{p}) \cap\left ( \mk{T}_{a}^{st} \cup (\_,M,\_)\cup \mk{T}_{b}^{st}\right )  .
\nonumber 
\end{gather}
In Corollary  \ref{middle M' cap left right middle M} we showed that  $(b^p,M',a^{p}) \cap (\mk{T}_{a}^{st} \cup (\_,M,\_)\cup \mk{T}_{b}^{st}  )$   and   $(b^p,M',a^{p}) \cup (\mk{T}_{a}^{st} \cup (\_,M,\_)\cup \mk{T}_{b}^{st}  )$    are  contractible (for any $p\in \ZZ$). Now using the equality above and Remark \ref{VK} one easily shows by induction that 
$  \mk{T}_{a}^{st} \cup (\_,M,\_)\cup \mk{T}_{b}^{st}\cup  \bigcup_{k=0}^n(b^{p+k},M',a^{p+k}) $ is contractible for any $p\in \ZZ$ and any $n\geq 1$. Applying Remark \ref{VK} again we deduce that the right-hand side of   \eqref{one union}  is contractible as well. Therefore  $\st(D^b(Q ))$ is contractible. 
\epr

\appendix

\section{Some contractible subsets of \texorpdfstring{$\RR^6$}{\space}} We prove here that some  subsets of  $\RR^6$, which we meet  in the proof of Theorem \ref{main theo}, are contractible. We start by the following subset
	\begin{lemma} \label{contractible 1} The   set $U_{>}$,  given below, is contractible: \begin{gather}  U_{>}= \left \{ (r_0,r_1,r_2,\phi_0,\phi_1,\phi_2) \in \RR^6 : \begin{array}{c} r_i>0\\ \phi_0 < \phi_1 < \phi_0+1 \\ \phi_0 < \phi_2 < \phi_0+1 \\ \arg_{(\phi_0,\phi_0+1)}(r_0 \exp(\ri \pi \phi_0) + r_1 \exp(\ri \pi \phi_1))>\phi_2\end{array} \right \}.\end{gather} The set $U_{<}$ defined by the same inequalities, except the last, where we take  $\arg_{(\phi_0,\phi_0+1)}(r_0 \exp(\ri \pi \phi_0) + r_1 \exp(\ri \pi \phi_1))<\phi_2$ is  contractible as well.  \end{lemma}
	
	\bpr    By drawing a picture one easily shows that:
	\begin{gather}  \label{contractible 11}  \forall   (r_0, r_1, r_2, \phi_0, \phi_1, \phi_2  ) \in U_{>}   \  \quad      \begin{array}{c}   r_1'\geq r_1  \\    r_2'>0<r_0' \leq r_0 \\ 
	 \phi_0<\phi_1\leq \phi_1'<\phi_0+1 \\
		\phi_0<\phi_2'\leq \phi_2<\phi_0+1 \end{array} \Rightarrow  (r_0', r_1', r_2' , \phi_0, \phi_1', \phi_2'  ) \in U_{>}.\end{gather}
	
	Let  	 $\gamma: \SS^n \rightarrow U $ be a continuous map  with  $n\geq 1$. Denote 
	\begin{gather} 0<r_0^{min} = \min\{ r_0(t): t\in \SS^n\}; \ \  0<r_1^{max}= \max\{ r_1(t): t\in \SS^n\};\nonumber  \\
	  0<u = \max\{ \phi_1(t)-\phi_0(t): t\in \SS^n\} <1; \ \  0<v = \min\{ \phi_2(t)-\phi_0(t): t\in \SS^n\}<1; \nonumber  \end{gather} 
	then by \eqref{contractible 11} for any $\delta >0$ and any $t\in \SS^n$, $s \in [0,1]$ the vector given below lies in $ U_{>} $:
	\begin{gather} F(t,s)=\left (\begin{array}{c} r_0(t)(1-s)+s r_0^{min}, r_1(t)(1-s)+s r_1^{max},  r_2(t)(1-s)+s \delta ,\\ 
	       \phi_0(t), \phi_0(t)+(1-s)(\phi_1(t)-\phi_0(t))+s u, \phi_0(t)+(1-s)(\phi_2(t)-\phi_0(t))+s v  \end{array}\right ).  \nonumber \end{gather}
	Hence we obtain a map $F:\SS^n \times [0,1]\rightarrow  U_{>}$, whose continuity is obvious.  This gives a homotopy from the map $\gamma$ to the following continuous  map: 
		\begin{gather} \gamma': \SS^n \rightarrow U_{>} \ \ \ \     \gamma'(t)=(r_0^{min}, r_1^{max}, \delta ,\phi_0(t), \phi_0(t)+ u,\phi_0(t)+ v )\end{gather}
Now we note that: 	
\begin{gather}  \label{contractible 12} \begin{array}{c} \forall (r_0,r_1,r_2,\phi_0,\phi_1,\phi_2) \in U   \ \  \forall \delta \in \RR \ \ \ (r_0,r_1,r_2,\phi_0+\delta ,\phi_1 + \delta ,\phi_2 + \delta) \in U_{>} \ \end{array} \end{gather}
Therefore for $t\in \SS^n$, $s \in [0,1]$  we have  \begin{gather} G(t,s) =\left (\begin{array}{c} r_0^{min}, r_1^{max}, \delta , \phi_0(t) +s (\phi_0(0)-\phi_0(t)),\\  \phi_0(t)+ u  +s (\phi_0(0)-\phi_0(t)),\phi_0(t)+ v  +s (\phi_0(0)-\phi_0(t)) \end{array} \right ) \in U_{>0}\nonumber \end{gather}
which gives a homotopy from $\gamma'$ to the constant map from $\SS^n$ to the point \\ $( r_0^{min}, r_1^{max}, \delta , \phi_0(0) ,  \phi_0(0)+u,\phi_0(0)+ v ) \in U_{>0}$. Thus, we showed that each continuous map $\gamma: \SS^n \rightarrow U_{>} $  with  $n\geq 1$ is homotopic to a constant map. 
If we show that $U_{>}$ is connected, then  Whitehead theorem ensures  that $U_{>}$ is contractible. Let $x= (r_0,r_1,r_2,\phi_0,\phi_1,\phi_2) \in U_{>} $ and  $x'= (r_0',r_1',r_2',\phi_0',\phi_1',\phi_2') \in U_{>} $. By \eqref{contractible 12}  we can move continuously  $x'$ in $U_{>}$ to  $x''=(r_0'',r_1'',r_2'',\phi_0,\phi_1'',\phi_2'')$ and now by \eqref{contractible 11} we can connect $x$, $x''$ by a continuous path in $U_{>}$.

The same idea shows that $U_{<}$ is contractible,  one must permute  $\leq \leftrightarrow \geq$, $\min \leftrightarrow \max$. 
The lemma is proved.
	\epr
	\begin{coro}  \label{contractible 1a} The   set $V$,  given below, is contractible: \begin{gather}  V= \left \{ (r_0,r_1,r_2,\phi_0,\phi_1,\phi_2) \in \RR^6 : \begin{array}{c} r_i>0\\ \phi_2-1 < \phi_0 < \phi_2 \\ \phi_2-1 < \phi_1 < \phi_2 \\ \arg_{(\phi_2-1,\phi_2)}(r_0 \exp(\ri \pi \phi_0) - r_2 \exp(\ri \pi \phi_2))>\phi_1\end{array} \right \}.\end{gather} After changing the last inequality to $\arg_{(\phi_2-1,\phi_2)}(r_0 \exp(\ri \pi \phi_0) - r_2 \exp(\ri \pi \phi_2))<\phi_1$ the  set remains contractible.  
	\end{coro}
\bpr  The assignment $(a_0,a_1,a_2,b_0,b_1,b_2)\mapsto (a_2,a_0,a_1,b_2-1,b_0,b_1)$  maps homeomorphically  the set $V$ to the set $U$ in Lemma \ref{contractible 1}. \epr

	\begin{coro}  \label{contractible 1b} The   set $V$,  given below, is contractible: \begin{gather}  V= \left \{ (r_0,r_1,r_2,\phi_0,\phi_1,\phi_2) \in \RR^6 : \begin{array}{c} r_i>0\\ \phi_0-1 < \phi_1 < \phi_0 \\ \phi_0-1 < \phi_2 < \phi_0 \\ \arg_{(\phi_0-1,\phi_0)}(r_0 \exp(\ri \pi \phi_0) + r_2 \exp(\ri \pi \phi_2))>\phi_1\end{array} \right \}.\end{gather}  
	\end{coro}
\bpr  The assignment $(a_0,a_1,a_2,b_0,b_1,b_2)\mapsto (a_0,a_2,a_1,-b_0,-b_2,-b_1)$  maps homeomorphically  the set $V$ to   the  set $U_{<}$  in Lemma \ref{contractible 1}(see \eqref{arg2}). \epr

	\begin{lemma} \label{contractible 2} The   set $U$,  given below, is contractible: \begin{gather}  U= \left \{ (r_0,r_1,r_2,\phi_0,\phi_1,\phi_2) \in \RR^6 : \begin{array}{c} r_i>0\\ \phi_2 < \phi_1 < \phi_0 < \phi_2+1 \\   \arg_{(\phi_2,\phi_2+1)}(r_0 \exp(\ri \pi \phi_0) + r_2 \exp(\ri \pi \phi_2))<\phi_1\end{array} \right \}.\end{gather}   \end{lemma}
	
	\bpr    By drawing a picture one easily checks  that:
	\begin{gather}  \label{contractible 21}  \forall   (r_0, r_1, r_2, \phi_0, \phi_1, \phi_2  ) \in U   \  \quad      \begin{array}{c}   r_2'\geq r_2  \\    r_1'>0<r_0' \leq r_0   \end{array} \Rightarrow  (r_0', r_1', r_2' , \phi_0, \phi_1, \phi_2  ) \in U.\end{gather}
	Let  	 $\gamma: \SS^n \rightarrow U $ be any continuous map  with  $n\geq 1$.       Denote 
	\begin{gather} 0<r_0^{min} = \min\{ r_0(t): t\in \SS^n\}; \ \  0<r_2^{max}= \max\{ r_2(t): t\in \SS^n\};\nonumber  \\[-2mm] 	 \label{contractible 24} \\[-2mm]
   0<u = \min\{ \phi_1(t)-\phi_2(t): t\in \SS^n\}<1.  \nonumber  \end{gather} 
			
			By drawing a picture one sees  that for big enough $A>r_2^{max}$ we have 
			\begin{gather} 	 \label{contractible 25} \forall \phi_2 \forall \phi_0 \ \  \  \phi_2< \phi_0 < \phi_2+1 \ \Rightarrow  \arg_{(\phi_2,\phi_2+1)}(r_0^{min} \exp(\ri \pi \phi_0) + A \exp(\ri \pi \phi_2))-\phi_2 <u. \end{gather}
This implication means that  for any $\delta >0$ the set $U'$, given below, is contained in $U$:
 \begin{gather} \label{contractible 22} U'= \left \{ (r_0^{min}, \delta ,A,\phi_0,\phi_1,\phi_2) \in \RR^6 : \begin{array}{c} \phi_2 < \phi_1 < \phi_0 < \phi_2+1 \\   u\leq \phi_1-\phi_2\end{array} \right \}\subset U\end{gather} 
where $A$, $r_0^{min}$, $u$ are fixed in  \eqref{contractible 24}, \eqref{contractible 25} and we chose any $\delta >0$.
			By \eqref{contractible 21} we see that  for  any $t\in \SS^n$, $s \in [0,1]$ we have:
	\begin{gather} F(t,s)=\left ( r_0(t)(1-s)+s r_0^{min}, r_1(t)(1-s)+s  \delta ,  r_2(t)(1-s)+s A , 
	       \phi_0(t), \phi_1(t), \phi_2(t)  \right ) \in U.  \nonumber \end{gather}

	Hence we obtain a  continuous map $F:\SS^n \times [0,1]\rightarrow  U$, which   is a  homotopy in $U$ from the map $\gamma$ to the following continuous  map: 
		\begin{gather} \gamma': \SS^n \rightarrow U \ \ \ \     \gamma'(t)=(r_0^{min}, \delta ,A,\phi_0(t), \phi_1(t),\phi_2(t)). \nonumber \end{gather}
Furthermore, by \eqref{contractible 24} we have $u\leq \phi_1(t)-\phi_2(t)$	 for $t\in \SS^n$, which means that $\im( \gamma')\subset U'$. Since $U'$ is contractible, there exists a homotopy  in $U'$ from $\gamma'$  to a constant map. Since $U'\subset U$(see \eqref{contractible 22}),  there exists a homotopy  in $U$ from $\gamma$  to a constant map. 
		
We show below that $U$ is connected, and then   by Whitehead theorem $U$ is contractible. 

  Let $x= (a_0,a_1,a_2,b_0,b_1,b_2) \in U $ and  $x'= (a_0',a_1',a_2',b_0',b_1',b_2') \in U $.  The formula  \eqref{contractible 12} holds again, and by using it we can move continuously  $x'$ in $U$ to a point $x''=(a_0',a_1',a_2',b_0'',b_1'',b_2)$. If we denote \begin{gather} 0<r_0^{min}=\min\{a_0,a_0'\}, 0<r_2^{max}=\max\{a_2,a_2'\}, 0<u = \min\{ b_1-b_2,  b_1''-b_2 \}<1 \nonumber \end{gather} then choose $A > r_2^{max}$ so that \eqref{contractible 25} holds with the chosen $u, r_0^{min}$, $r_2^{max}$, in particular  for any $\delta>0$  the corresponding  set $U'$ defined by \eqref{contractible 22} is  a subset of $U$.   By the properties \eqref{contractible 21} and by the choice of $u, r_0^{min}$, $A$, $\delta$  we can move  the points $x$ and $x''$, by changing only $a_0,a_1,a_2, a_0',a_1',a_2'$,   continuously in $U$ to points $y$, $y'$ in $U'$, respectively. Now the connectivity of $U$ follows from the connectivity of $U'$.
The lemma is proved.
	\epr	
	\begin{coro}   \label{contractible 2a} The   set $V$,  given below, is contractible: \begin{gather}  V= \left \{ (r_0,r_1,r_2,\phi_0,\phi_1,\phi_2) \in \RR^6 : \begin{array}{c} r_i>0\\ \phi_2-1 < \phi_1 < \phi_0 < \phi_2 \\   \arg_{(\phi_2-1,\phi_2)}(r_0 \exp(\ri \pi \phi_0) - r_2 \exp(\ri \pi \phi_2))<\phi_1\end{array} \right \}.\end{gather}  
	\end{coro}
\bpr  The assignment $(a_0,a_1,a_2,b_0,b_1,b_2)\mapsto (a_0,a_1,a_2,b_0,b_1,b_2-1)$  maps homeomorphically  the set $V$ to the set $U$ in Lemma \ref{contractible 2}. \epr	
\begin{remark} \label{VK} If we have two contractible open subsets $U$, $V$
 in a f.d. manifold $M$ and the interesection $U\cap V$ is contractible, then by Seifert-van Kampen  theorem, Mayer-Vietoris sequence, Hurewicz theorem  and Whitehead theorem it follows that $U \cup V$ is contractible.  

If $U=\bigcup_{i\in A} U_i$ is an union of open  subsets in  a f.d. manifold $M$ and for any finite subset $F \subset A$ we have that  $\bigcup_{i\in F} U_i$  is contractible, then using Witehead theorem one can easily show that  $U$ is contractible as well. \end{remark}


\begin{thebibliography}{99}

\bibitem{Bridg1}  T. Bridgeland, Stability conditions on triangulated categories, Annals of Math.  166  no. 2, 317-345 (2007).

\bibitem{Bridg2}  T. Bridgeland, Stability conditions on K3 surfaces, Duke Math. J. 141, no. 2, 241-291 (2008).

\bibitem{BS} T. Bridgeland and I. Smith, \textit{Quadratic differentials as
stability conditions}, arXiv:1302.7030

\bibitem{BBD} A.A.  Beilinson, J. Bernstein and P. Deligne, Faisceaux Pervers. Astґerisque 100, Soc. Math. de France (1982).

\bibitem{CP} J. Collins and A. Polishchuck, Gluing Stability Conditions, Adv. Theor. Math. Phys. Volume 14, Number 2 (2010), 563-608.


\bibitem{WCB1}  W. Crawley-Boevey, Exceptional sequences of representations of quivers,  in 'Representations of algebras', Proc. Ottawa 1992, eds V. Dlab and H. Lenzing, Canadian Math. Soc. Conf. Proc. 14 (Amer. Math. Soc., 1993), 117-124.

\bibitem{WCB2}  W. Crawley-Boevey, Lectures on Representations of Quivers, \url{http://www1.maths.leeds.ac.uk/~pmtwc/quivlecs.pdf}

\bibitem{DHKK} G. Dimitrov, F. Haiden, L. Katzarkov, M. Kontsevich, Dynamical Systems and Categories, arXiv:1307.8418.

\bibitem{DK1} G. Dimitrov,  L. Katzarkov, Non-semistable exceptional objects in hereditary categories, arXiv:1311.7125

\bibitem{GM}  I. Gelfand, Y. Manin, Methods of Homological algebra, 2nd edition, Springer.

\bibitem{Kac}  V.G. Kac: Infinite root systems, representations of graphs and invariant theory. Inventiones mathematicae 56, 57-92 (1980).

\bibitem{KS} M. Kontsevich and Y. Soibelman, Motivic
Donaldson-Thomas invariants: summary of results, \textit{Mirror symmetry and
tropical geometry}, 55--89, Contemp. Math., 527, Amer. Math. Soc.,
Providence, RI, 2010.


\bibitem{King} A. King, Moduli of representations of finite dimensional
algebras, \textit{Quart. J. Math. Oxford Ser.} 45 (1994), 515--530.

\bibitem{Lam}  T. Y. Lam,  A first course in non-commutative rings, Springer-Verlag.

\bibitem{Lenzing} H. Lenzing, Hereditary Categories Lectures 1,2, Advanced ICTP-school on Representation Theory and Related Topics
(9-27 January 2006), \url{https://www.imj-prg.fr/~bernhard.keller/ictp2006/lecturenotes/lenzing1.pdf}

\bibitem{Macri}   E.  Macr\`i, Stability conditions on curves, Math. Res. Lett. 14 (2007) 657-672. Also arXiv:0705.3794.

\bibitem{Woolf}  J. Woolf, Algebraic stability conditions and contractible stability spaces, arXiv:1407.5986.
\end{thebibliography}
\end{document}